\documentclass[A4paper, 12pt, french]{article}

\usepackage{cite}
\usepackage[portrait, top=2cm, bottom=2cm, left=2cm, right=2cm] {geometry}
\usepackage[french]{babel}
\usepackage{url}
\usepackage{array,epsfig}
\usepackage{amsmath}
\usepackage{amsfonts}
\usepackage{amssymb}
\usepackage{amsxtra}
\usepackage[thmmarks, amsmath, thref, amsthm]{ntheorem}
\usepackage{latexsym}
\usepackage{dsfont}
\usepackage{mathrsfs}
\usepackage{color}
\usepackage[all]{xy}
\usepackage{tensor}
\usepackage[pagebackref=true]{hyperref}

\usepackage{graphicx}
\usepackage{mathtools}
\usepackage{caption}
\usepackage[OT2, T1]{fontenc}
\usepackage{accents}
\usepackage[nottoc]{tocbibind}
\usepackage[page,toc,titletoc,title]{appendix}
\usepackage{sectsty}

\theoremstyle{plain}
\newtheorem{thm}{Th\'eor\`eme}[section]
\newtheorem{lem}[thm]{Lemme}

\newtheorem{prop}[thm]{Proposition}
\newtheorem{cor}[thm]{Corollaire}

\newtheorem{customthm}{Th\'eor\`eme}

\theoremstyle{definition}
\newtheorem{defn}[thm]{D\'efinition}
\newtheorem{rmk}[thm]{Remarque}
\newtheorem*{remercie}{Remerciements}

\theoremstyle{remark}

\numberwithin{equation}{section}

\newcommand{\tuple}[1]{\left(#1\right)}

\newcommand{\pair}[1]{\left<#1\right>}

\newcommand{\ol}[1]{\overline{#1}}

\renewcommand{\le}{\leqslant}
\renewcommand{\ge}{\geqslant}

\renewcommand{\H}{\mathrm{H}}

\newcommand{\Z}{\mathrm{Z}}
\newcommand{\diff}{\mathrm{d}}

\newcommand{\Cbb}{\mathbb{C}}
\newcommand{\Gbb}{\mathbb{G}}

\newcommand{\Pbb}{\mathbb{P}}

\newcommand{\Qbb}{\mathbb{Q}}
\newcommand{\Rbb}{\mathbb{R}}
\newcommand{\Zbb}{\mathbb{Z}}

\newcommand{\Escr}{\mathscr{E}}

\newcommand{\Xscr}{\mathscr{X}}

\newcommand{\afrak}{\mathfrak{a}}
\newcommand{\gfrak}{\mathfrak{g}}

\newcommand{\xfrak}{\mathfrak{x}}
\newcommand{\yfrak}{\mathfrak{y}}

\newcommand{\id}{\operatorname{id}}
\newcommand{\ab}{\operatorname{ab}}

\newcommand{\GL}{\operatorname{GL}}
\newcommand{\PGCD}{\operatorname{PGCD}}
\newcommand{\PPCM}{\operatorname{PPCM}}
\newcommand{\SL}{\operatorname{SL}}
\newcommand{\BM}{\operatorname{BM}}
\newcommand{\Spec}{\operatorname{Spec}}
\newcommand{\inv}{\operatorname{inv}}

\newcommand{\Gal}{\operatorname{Gal}}
\newcommand{\Ker}{\operatorname{Ker}}

\newcommand{\Aut}{\operatorname{Aut}}
\newcommand{\Out}{\operatorname{Out}}
\newcommand{\Img}{\operatorname{Im}}

\newcommand{\Br}{\operatorname{Br}}
\newcommand{\nr}{\operatorname{nr}}
\newcommand{\mr}{\operatorname{mr}}

\newcommand{\res}{\operatorname{res}}
\newcommand{\lien}{\operatorname{lien}}
\newcommand{\sh}{\operatorname{sh}}
\newcommand{\loc}{\operatorname{loc}}
\renewcommand{\int}{\operatorname{int}}
\newcommand{\Int}{\operatorname{Int}}

\newcommand{\Hom}{\operatorname{Hom}}

\newcommand{\Ind}{\operatorname{Ind}}

\DeclareSymbolFont{cyrletters}{OT2}{wncyr}{m}{n}
\DeclareMathSymbol{\Sha}{\mathalpha}{cyrletters}{"58}
\DeclareMathSymbol{\Be}{\mathalpha}{cyrletters}{"42}

\title{Sur les espaces homog\`enes de Borovoi--Kunyavski\u{\i}}
\author{Nguy$\tilde{\hat{\text{e}}}$n M$\d{\text{a}}$nh Linh}
\date{\today}

\begin{document}

\maketitle

\begin{abstract}
    Nous \'etablissons le principe de Hasse et l'approximation faible pour certains espaces homog\`enes de $\SL_n$ \`a stabilisateur g\'eom\'etrique nilpotent de classe 2, construits par Borovoi et Kunyavski\u{\i}. Ces espaces homog\`enes v\'erifient donc une conjecture de Colliot-Th\'el\`ene concernant l'obstruction de Brauer--Manin pour les vari\'et\'es g\'eom\'etriquement rationnellement connexes.
\end{abstract}

\tableofcontents

\bigskip
\section{Introduction}

\subsection{Int\'er\^et du probl\`eme} \label{subsection11}

Soit $X$ une vari\'et\'e lisse d\'efinie sur un corps de nombres $k$. On dira que $X$ est un contre-exemple au principe de Hasse si $X$ poss\`ede des points locaux dans tous les compl\'et\'es $k_v$ de $k$ mais $X(k) = \varnothing$. Manin \cite{manin1971} a introduit une m\'ethode g\'en\'erale pour l'\'etude des obstructions \`a l'existence de $k$-points sur $X$: \`a $X$ on associe le groupe de Brauer non ramifi\'e $\Br_{\nr} X$, et l'on consid\`ere l'{\em accouplement de Brauer--Manin}
	\begin{equation} \label{eqBrauerManin}
		\prod_{v} X(k_v) \times \Br_{\nr} X, \quad \pair{(P_v)_v,A}_{\BM} = \sum_{v} \inv_v(A(P_v)),
	\end{equation}
o\`u $v$ parcourt les places de $k$, et o\`u $\inv_v: \Br k_v \hookrightarrow \Qbb/\Zbb$ d\'esigne l'invariant local. Notant $\Br_0 X = \Img(\Br k \to \Br X)$, on voit par la loi de reciprocit\'e globale que l'accouplement \eqref{eqBrauerManin} se factorise par $(\Br_{\nr} X) / (\Br_0 X)$.

On plonge $X(k)$ diagonalement dans $\prod_v X(k_v)$. Alors l'adh\'erence de $X(k)$ (pour la topologie produit) est contenue dans l'ensemble (ferm\'e) $\left(\prod_v X(k_v)\right)^{\Br}$ des familles de points locaux de $X$ qui sont orthogonales \`a $\Br_{\nr} X$ par rapport \`a l'accouplement \eqref{eqBrauerManin}. Si $X(k) \neq \varnothing$, on dira que l'{\em obstruction de Brauer--Manin \`a l'approximation faible pour $X$ est la seule} si $X(k)$ est dense dans $\left(\prod_v X(k_v)\right)^{\Br}$. On dira que l'{\em obstruction de Brauer--Manin au principe de Hasse est la seule} pour une classe $\Xscr$ de vari\'et\'es lisses sur $k$ si pour toute $X \in \Xscr$, $\left(\prod_v X(k_v)\right)^{\Br} \neq \varnothing$ entra\^ine $X(k) \neq \varnothing$. Notons que ces deux propri\'et\'es sont des invariants birationnels stables. 

Une conjecture de Colliot-Th\'el\`ene \cite[Conjecture 14.1.2]{colliot2021brauer} pr\'edit que l'obstruction de Brauer--Manin au principe de Hasse et \`a l'approximation faible est la seule pour les espaces homog\`enes de groupes alg\'ebriques lin\'eaires. \`A la suite du travail de Demarche et Lucchini Arteche \cite[Th\'eor\`eme 1.1]{demarche2019reduction}, le cas g\'en\'eral se r\'eduit au cas o\`u le groupe ambiant est $\SL_n$ et o\`u le stabilisateur g\'eom\'etrique est fini. On sait que cette conjecture vaut si le stabilisateur est ab\'elien \cite[Theorem 2.2]{bor96}; dans ce cas, l'obstruction \`a l'existence de $k$-points sur $X$ est en fait contr\^ol\'ee par un morphisme $\Be(X) \to \Qbb/\Zbb$, o\`u $\Be(X)$ est le sous-groupe des \'el\'ements \og partout localement constants \fg{} de $(\Br_{\nr} X) / (\Br_0 X)$, c'est la {\em premi\`ere obstruction de Brauer--Manin} de $X$. Dans le cas o\`u le stabilisateur g\'eom\'etrique n'est pas ab\'elien, on ne sait pas si le groupe $\Be(X)$ est suffisant pour expliquer la d\'efaut du principe de Hasse. Borovoi et Kunyavski\u{\i} ont essay\'e de construire un espace homog\`ene o\`u l'obstruction par rapport \`a $\Be(X)$ ne suffit pas \cite{Bokun}, mais il s'av\`ere que cette construction ne fonctionne pas (voir \cite{BokunErratum}). Nous allons, dans ce texte, \'etudier leur construction (qu'on va appeler {\em espaces homog\`enes de Borovoi--Kunyavski\u{\i}}), et notamment nous montrons que ces espaces homog\`enes v\'erifient toujours le principe de Hasse ainsi que l'approximation faible; ils v\'erifient donc la conjecture de Colliot-Th\'el\`ene mentionn\'ee ci-dessus (en fait, il n'y pas d'obstruction de Brauer--Manin pour ceux-ci). 

Il est \`a noter qu'un r\'esultat r\'ecent de Harpaz et Wittenberg \cite[Th\'eor\`eme B]{Harpaz2020} dit que la conjecture de Colliot-Th\'el\`ene vaut si le stabilisateur g\'eom\'etrique est {\em hyper-r\'esoluble} (en tant que groupe fini muni d'une action ext\'erieure du groupe de Galois absolu). Cela n'est pas toujours le cas pour les espaces homog\`enes de Borovoi--Kunyavski\u{\i} (voir le corollaire \ref{corNotSuperSolvabe} ci-dessous), malgr\'e le fait que les stabilisateurs de ceux-ci sont nilpotents (un groupe {\em abstrait} fini et nilpotent est toujours hyper-r\'esoluble, mais ce n'est pas le cas pour les groupes alg\'ebriques finis).

\subsection{Notations et conventions} \label{subsection12}

On va fixer dans ce texte quelque notations.

Si $K$ est un corps parfait, $\ol{K}$ d\'esigne une cl\^oture alg\'ebrique de $K$ et $\Gamma_K = \Gal(\ol{K}/K)$ est le groupe de Galois absolu de $K$.

Soit $K$ un corps de caract\'eristique nulle. Si $F$ est un $K$-groupe fini ({\em i.e.} un sch\'ema en groupes fini sur $\Spec K$), on l'identifiera au groupe abstrait $F(\ol{K})$ muni de l'action naturelle de $\Gamma_K$. On note aussi $\hat{F} = \Hom(F,\ol{K}^\times)$ le $\Gamma_K$-module des caract\`eres de $F$.

Lorsque $F$ est un groupe abstrait ou un sch\'ema en groupes, on note $[F,F]$ son sous-groupe d\'eriv\'e, $F^{\ab} = F/[F,F]$ son ab\'elianis\'e, $Z(F)$ son centre, $\Aut(F)$ son groupe des automorphismes (en tant que groupe abstrait), $\Int(F) = F/Z(F)$ son groupe des automorphismes int\'erieurs, $\int(f) \in \Int(F)$ la conjugaison par un \'el\'ement $f \in F$, et $\Out(F) = \Aut(F)/\Int(F)$ son groupe des automorphismes ext\'erieurs. Si $F$ est un groupe fini et $\Gamma$ est un groupe profini, une action (resp. {\em action ext\'erieure}) de $\Gamma$ sur $F$ est un morphisme {\em continu} ({\em i.e.} localement constant) de $\Gamma$ dans $\Aut(F)$ (resp. dans $\Out(F)$). 

Rappelons la notion de groupe fini hyper-r\'esoluble au sens de Harpaz--Wittenberg \cite[D\'efinition 6.4]{Harpaz2020}. Soit $\kappa: \Gamma \to \Out(F)$ une action ext\'erieure d'un groupe profini $\Gamma$ sur un groupe fini $F$. Un sous-groupe distingu\'e $G$ de $F$ est dit stable sous $\kappa$ si pour tout $\sigma \in \Gamma$, il existe un relev\'e $\phi \in \Aut(F)$ de $\kappa(\sigma)$ tel que $\phi(G) \subseteq G$ (dans ce cas, comme $G$ est un sous-groupe distingu\'e de $F$, tous les relev\'es de $\kappa(\sigma)$ v\'erifient cette propri\'et\'e). On dit que $F$ (muni de l'action ext\'erieure $\kappa$) est hyper-r\'esoluble s'il existe un entier $n$ et une suite
    \begin{equation*}
        \{1\} = F_0 \subseteq F_1 
 \subseteq \cdots \subseteq F_n = F
    \end{equation*}
de sous-groupes distingu\'es de $F$, stables sous $\kappa$, tels que $F_i/F_{i-1}$ soit cyclique pour tout $i \in \{1,\ldots,n\}$. Il est clair que dans ce cas, pour tout sous-groupe distingu\'e $G$ de $F$ qui est stable sous $\kappa$, le quotient $F/G$ muni de l'action ext\'erieure de $\Gamma$ induite par $\kappa$ est \'egalement hyper-r\'esoluble\footnote{Attention, $\kappa$ n'induit pas forc\'ement une action ext\'erieure de $\Gamma$ sur $G$.}.

Les cohomologies utilis\'ees seront cohomologie de groupes profinis, cohomologie galoisienne et cohomologie \'etale. On dispose aussi de la notion des $\H^0$ et $\H^1$ non ab\'eliens. Le $\H^2$ non ab\'elien, d\'efini dans \cite[\S 1]{FSS}, sera mentionn\'e dans le paragraphe \ref{subsection31}.

Si $G$ est un groupe et $H$ est un sous-groupe distingu\'e de $G$, l'image d'un \'el\'ement $\sigma \in G$ dans $G/H$ sera not\'ee $\ol{\sigma}$. Pour tout $r \ge 1$, l'image dans $\H^r$ d'un $r$-cocycle $a$ de $\Z^r$ sera not\'ee $[a]$.

Si $G$ est un groupe abstrait et $P \subseteq G$, $\pair{P}$ d\'esigne le sous-groupe de $G$ (alg\'ebriquement) engendr\'e par $P$. Si de plus $G$ est suppos\'e topologique, l'adh\'erence $\ol{\pair{P}}$ est le sous-groupe de $G$ topologiquement engendr\'e par $P$. 

Soit $k$ un corps de nombres. L'ensemble des places de $k$ sera not\'e $\Omega_k$. Si $A$ est un $\Gamma_k$-module, $r$ est un entier et $v \in \Omega_k$, on notera $\loc_v: \H^r(k,A) \to \H^r(k_v,A)$ le morphisme de localisation, et $\alpha_v = \loc_v(\alpha)$ pour tout $\alpha \in \H^r(k,A)$. De plus, on d\'efinit les sous-groupes
	\begin{equation*}
	    \Sha^r_S(k,A):=\Ker\left(\H^r(k,A) \to \prod_{v \in \Omega_k \setminus S} \H^r(k_v,A)\right)
     \end{equation*}
pour tout sous-ensemble fini $S \subseteq \Omega_k$, et
    \begin{equation*}
        \Sha^r_\omega(k,A) := \bigcup_{S} \Sha^r_S(k,A), \quad \Sha^r(k,A) :=  \Sha^r_\varnothing(k,A).
    \end{equation*}

Soit $X$ une vari\'et\'e lisse et g\'eom\'etriquement int\`egre sur un corps $K$ de caract\'eristique nulle. Le groupe de Brauer de $X$ est toujours le groupe de Brauer--Grothendieck $\Br X:=\H^2(X,\Gbb_m)$. Le groupe de Brauer non ramifi\'e $\Br_{\nr} X$ est le groupe de Brauer de n'importe quelle vari\'et\'e propre et lisse contenant $X$ comme un ouvert dense; c'est naturellement un sous-groupe de $\Br X$, et c'est un invariant birationnel stable (voir \cite[\S 6.2]{colliot2021brauer}). Si $K$ est un corps, le groupe de Brauer de $K$ est alors $\Br K = \H^2(K,\Gbb_m)$. Notons que $\H^2(K,\mu_n) = (\Br K)[n]$ pour tout entier $n$. Dans le cas o\`u $K$ est un corps $p$-adique ou $\Rbb$ ou $\Cbb$, $\inv_K: \Br K \hookrightarrow \Qbb/\Zbb$ d\'esigne l'invariant local. Si $k$ un corps de nombres et $v$ est une place de $k$, on notera $\inv_v = \inv_{k_v}$.
 
Rappelons aussi l'intepr\'etation cohomologique suivante de l'approximation faible, dont la preuve se trouve par exemple dans \cite[\S 1.2]{harari_affbm}.

\begin{lem} \label{lemWA}
    Soient $k$ un corps de nombres et $F$ un $k$-groupe fini. On choisit un plongement $F \hookrightarrow \SL_n$ de $k$-groupes et l'on pose $X = F \backslash \SL_n$. Alors pour tout sous-ensemble fini $S \subseteq \Omega_k$, on a \'equivalence entre
    
    \begin{enumerate}
        \item $X(k)$ est dense dans $\prod_{v \in S} X(k_v)$;
        \item la restriction $\H^1(k,F) \to \prod_{v \in S} \H^1(k_v, F)$ est surjective.
    \end{enumerate}
\end{lem}

\subsection{\'Enonc\'es des principaux th\'eor\`emes} \label{subsection13} 

Soit $k$ un corps de nombres. On rappelle la construction de Borovoi--Kunyavski\u{\i}. On construit une extension centrale de $k$-groupes finis
    \begin{equation} \label{eqExtension0}
	0 \to Z \to F \to M \oplus M \to 0
    \end{equation}
avec $Z = Z(F) = [F,F]$ (donc $F$ est nilpotent de classe 2) et on consid\`ere les espaces homog\`enes $X$ de $\SL_n$ \`a stabilisateur g\'eom\'etrique $F$. Dans le \S \ref{section3}, cette construction sera discut\'ee en d\'etail. Pour le stabilisateur $F$ comme dans les th\'eor\`emes \ref{thmA} et \ref{thmB} ci-dessous, on verra qu'il n'y a pas d'obstruction de Brauer--Manin pour ces espaces homog\`enes $X$: leurs groupes de Brauer non ramifi\'es sont r\'eduits aux constantes (proposition \ref{propBorovoiKunyavskii}).

Voici les principaux r\'esultats de ce texte.

\begin{customthm} \label{thmA}
    Soient $A$ un group ab\'elien fini, $k$ un corps de nombres contenant $\mu_{\exp(A)}$, $L/k$ une extension finie galoisienne, $M = \Ind_{\Gamma_k}^{\Gamma_L} A$, $j: A \otimes A \hookrightarrow M \otimes M$ l'inclusion canonique, $Z$ le conoyau de $j$ et $\phi: M \otimes M \to Z$ la projection. On d\'efinit l'application biadditive
        \begin{equation*}
            \Phi: (M \oplus M) \times (M \oplus M) \to Z, \quad \Phi((x,y),(x',y')) = \phi(x \otimes y'),
        \end{equation*}
    vu comme $2$-cocycle normalis\'e sur $M \oplus M$ \`a coefficients dans $Z$ (muni de l'action triviale de $M \oplus M$) et soit $F$ le produit crois\'e $Z \times_\Phi (M \oplus M)$, muni de l'action coordonn\'ees par coordonn\'ees de $\Gamma_k$. Alors les espaces homog\`enes de $\SL_n$ de lien de Springer $\lien(F)$ (voir paragraphe \ref{subsection31} pour la d\'efinition du lien de Springer d'un espace homog\`ene) v\'erifient le principe de Hasse.
\end{customthm}

\begin{customthm} \label{thmB}
    Avec les m\^emes hypoth\`eses du th\'eor\`eme \ref{thmA}, soit $X$ un espace homog\`ene de $\SL_n$ de lien de Springer $\lien(F)$. Si $X(k) \neq \varnothing$, alors $X$ v\'erifie l'approximation faible.
\end{customthm}

Les th\'eor\`emes \ref{thmA} et \ref{thmB} seront d\'emontr\'es dans les paragraphes \ref{subsection53} et \ref{subsection54} respectifs du \S \ref{section5}. Leurs preuves reposent sur des calculs cohomologiques: comme $M = \Ind_{\Gamma_k}^{\Gamma_L} A$, on peut utiliser le lemme de Shapiro pour passer de $\H^r(k,M)$ \`a $\H^r(L,A)$. La compatibilit\'e avec l'isomorphisme de Shapiro de certaines applications seront donc n\'ecessaires, et nous allons en discuter dans le \S \ref{section2}.	
	
Pour toute extension $K$ de $k$, la fl\`eche connectante $\Delta: \H^1(K,M \oplus M) \to \H^2(K,Z)$ induite par l'extension centrale \eqref{eqExtension0} s'exprime en termes de cup-produits (voir lemme \ref{lemImageOfDelta}). Cette observation importante nous permet de r\'eduire le probl\`eme des points rationnels sur $X$ aux calculs cohomologiques. L'argument cl\'e sera le th\'eor\`eme \ref{thmKey}, qui est lui-m\^eme une application d'un \og lemme arithm\'etique \fg{}, les propositions \ref{propSerre} et \ref{propArithmetic} du paragraphe \ref{subsection41}. Ces propositions sont des g\'en\'eralisations de \cite[Chapitre I, \S 2.2, Th\'eor\`eme 4]{Serre}, qui d\'ecrivent une propri\'et\'e globale du symbole de Hilbert, et dont les d\'emonstrations reposent sur la dualit\'e de Poitou--Tate. 

Montrons que le stabilisateur des espaces homog\`enes de Borovoi--Kunyavski\u{\i} comme dans l'\'enonc\'e du th\'eor\`eme \ref{thmA} n'est pas toujours hyper-r\'esoluble au sens de Harpaz--Wittenberg.

\begin{lem} \label{lemNotSuperSolvabe}
    Soit $\gfrak$ le groupe cyclique $\pair{g\,|\,g^3 = 1}$ et soit $M = (\Zbb/2)[\gfrak] = \bigoplus_{i=0}^2 (\Zbb/2) e_i$. Alors le quotient $\ol{M} = M / (e_0 + e_1 + e_2)$ est un $\gfrak$-module simple mais pas un groupe cyclique.
\end{lem}
\begin{proof}
    Le groupe $\ol{M}$ est bien un $\gfrak$-module puisque $e_0 + e_1 + e_2 \in M$ est fix\'e par $\gfrak$. Il est d'ordre $4$ et d'exposant $2$ donc non cyclique. Montrons qu'il est un $\gfrak$-module simple. Soit $N \subseteq \ol{M}$ un sous-groupe $\gfrak$-stable et supposons $N \neq 0$. Soit $x \in N \setminus \{0\}$ et affirmons que $x$ n'est pas fix\'e par $\gfrak$. En effet, on peut \'ecrire $x = \ol{a_1e_1 + a_2e_2}$, o\`u $a_1,a_2 \in \Zbb/2$. Alors $\tensor[^g]{x}{} = \ol{a_1 e_0 + a_2 e_1}$. Si $\tensor[^g]{x}{} = x$, alors il existe $b \in \Zbb/2$ tel que $a_1 e_0 + a_2 e_1 = a_1e_1 + a_2e_2 + b(e_0 + e_1 + e_2)$, d'o\`u 
	\begin{equation*}
		a_1 = b, \quad a_2 = a_1 + b \quad \text{et} \quad a_2 + b = 0,
	\end{equation*}
    et donc $a_1 = a_2 = 0$, ou $x = 0$. Cette contradiction montre qu'en fait on a $\tensor[^g]{x}{} \neq x$. On en d\'eduit que $N$ contient 4 \'el\'ements deux \`a deux distincts, \`a savoir $0, x, \tensor[^g]{x}{}$ et $\tensor[^{g^2}]{x}{}$. Comme $\ol{M}$ est d'ordre 4, on voit que $N = \ol{M}$, ce qui montre que $\ol{M}$ n'a pas de sous-groupe $\gfrak$-stable non trivial.
\end{proof}

\begin{cor} \label{corNotSuperSolvabe}
    Soit $L/k$ une extension galoisienne de degr\'e 3 de corps de nombres et $M = \Ind_{\Gamma_k}^{\Gamma_L} \mu_2$. De $M$, on construit le $k$-groupe fini $F$ comme dans l'\'enonc\'e du th\'eor\`eme \ref{thmA}. Alors $F$ n'est pas hyper-r\'esoluble (en tant que groupe fini muni de l'action ext\'erieure de $\Gamma_k$).
\end{cor}
\begin{proof}
    En effet, si $F$ \'etait hyper-r\'esoluble, alors $M \oplus M = F/Z$ est un $\Gamma_k$-module hyper-r\'esoluble. On en d\'eduirait que ce serait aussi le cas pour $M$ et le quotient $\ol{M}$ comme dans le lemme \ref{lemNotSuperSolvabe}. Mais $\ol{M}$ est un $\Gamma_k$-module simple non cyclique, donc non hyper-r\'esoluble. Ainsi, $F$ n'est pas hyper-r\'esoluble.
\end{proof}
 
\section{Compatibilit\'e avec les isomorphismes de Shapiro} \label{section2}

Dans cette section, $G$ est un groupe profini et $H$ est un sous-groupe ouvert, distingu\'e de $G$.

Soit $r$ un entier. On rappelle que si $A$ est un group ab\'elien muni de l'action triviale de $H$, le groupe $\Z^r(H,A)$ des $r$-cocycles $H^r \to A$ est muni de l'action suivante de $G$:
    \begin{equation*}
        \forall \sigma \in G,\forall a \in \Z^r(H,A),\forall \tau_1,\ldots,\tau_r \in H, \quad (\tensor[^\sigma]{a)}{_{\tau_1,\ldots,\tau_r}} = a_{\sigma^{-1}\tau_1\sigma,\ldots,\sigma^{-1}\tau_r \sigma}.
    \end{equation*}
Cette action induit une action de $G$ sur le groupe $\H^r(H,A)$. Le groupe $H$ agit trivialement sur $\H^r(H,A)$ puisque les automorphismes int\'erieurs de $H$ induisent l'identit\'e. On obtient ainsi une action de $G/H$ sur $\H^r(H,A)$. Pour le reste de cette section, soit $M = \Ind_G^H A$ (voir \cite[Chapitre I, \S 2.5]{serre1994cohomologie} pour la notion de module induit). Alors $M$ est un $G$-module discret dont le groupe ab\'elien sous-jacent est celui des applications $G/H \to A$, et dont l'action de $G$ est donn\'ee par la formule
    \begin{equation*}
        \forall \sigma \in G,\forall f \in M,\forall g \in G/H, \quad (\tensor[^\sigma]{f)}{}(g) = f(g\ol{\sigma})
    \end{equation*}
Pour tout $r \ge 1$, on dispose d'un morphisme
    \begin{equation*}
        \sh: \Z^r(G,M) \to \Z^r(H,A)
    \end{equation*}
qui \`a chaque $r$-cocycle $a: G^r \to M$ associe le $r$-cocycle
    \begin{equation*}
        H^r \to A, \quad (\sigma_1,\ldots,\sigma_r) \mapsto a_{\sigma_1,\ldots,\sigma_r}(1).
    \end{equation*}
Le lemme de Shapiro affirme que $\sh$ induit un isomorphisme $\sh: \H^r(G,M) \to \H^r(H,A)$. Le but de cette section est de donner quelques r\'esultats de compatibilit\'e des isomorphismes de Shapiro avec certaines applications.

\subsection{Quasi-inverses aux applications de Shapiro} \label{subsection21}

On fixe une section (ensembliste) $u: G/H \to G$ avec $u(1) = 1$. Pour tous $g \in G/H$ et $\sigma \in G$, posons $\gamma(g,\sigma):=u(g) \sigma u(g\ol{\sigma})^{-1}$. Alors $\gamma(g,\sigma) \in H$ puisque son image dans $G/H$ est $g\ol{\sigma}(g\ol{\sigma})^{-1} = 1$. L'application $\gamma: G/H \times G \to H$ est continue. Elle v\'erifie une \og condition de cocycle\fg{}:
    \begin{equation} \label{eqCocycleConditionForGamma}
	\forall g \in G/H, \forall \sigma,\tau \in G, \quad \gamma(g,\sigma\tau) = \gamma(g,\sigma)\gamma(g\ol{\sigma},\tau).
    \end{equation}
Les r\'esultats suivants d\'ecrivent des morphismes $\Z^r(H,A) \to \Z^r(G,M)$ pour $r = 1, 2$, qui induisent des inverses de $\sh$ au niveau de la cohomologie. Ils sont inspir\'es par \cite[Lemma 2.1, Lemma 2.2]{Naidu}.

\begin{lem} \label{lemShapiro1}
    Soit $a \in \H^1(H,A)$, vu comme morphisme continu $H \to A$. Soit $x: G \to M$ l'application continue d\'efinie par
	\begin{equation*}
            \forall \sigma \in G, \forall g \in G/H, \quad x_\sigma(g) = a_{\gamma(g,\sigma)}.
        \end{equation*}
    Alors $x$ est un $1$-cocycle. De plus, pour tous $\sigma \in H$ et $g \in G/H$, on a $x_\sigma(g) = (\tensor[^{u(g)^{-1}}]{a)}{_\sigma}$. En particulier, $\sh(x) = a$.
\end{lem}
\begin{proof}
    Pour tous $g \in G/H$ et $\sigma,\tau\in G$, au vu de \eqref{eqCocycleConditionForGamma}, on a
	\begin{equation*}
            x_{\sigma\tau}(g) = a_{\gamma(g,\sigma\tau)} = a_{\gamma(g,\sigma)\gamma(g\ol{\sigma},\tau)} = a_{\gamma(g,\sigma)} + a_{\gamma(g\ol{\sigma},\tau)} = x_\sigma(g) + x_\tau(g\ol{\sigma}) = x_\sigma(g) + \tensor[^\sigma]{x}{_\tau}(g),
        \end{equation*}
    d'o\`u $x_{\sigma\tau} = x_\sigma + \tensor[^\sigma]{x}{_\tau}$, donc $x$ est bien un $1$-cocycle.
	
    Soient $\sigma \in H$ et $g \in G/H$. Alors $\gamma(g,\sigma) = u(g) \sigma u(g\ol{\sigma})^{-1} = u(g) \sigma u(g)^{-1}$, d'o\`u 
	\begin{equation*}
            x_\sigma(g) = a_{u(g) \sigma u(g)^{-1}} = (\tensor[^{u(g)^{-1}}]{a)}{_\sigma}.
        \end{equation*}
    En particulier, $x_\sigma(1) = a_\sigma$ puisque $u(1) = 1$ et donc $\sh(x) = a$.
\end{proof} 

\begin{lem} \label{lemShapiro2}
    Soit $a \in \Z^2(H,A)$. Soit $x: G \times G \to M$ l'application continue d\'efinie par
	\begin{equation*}
            \forall \sigma,\tau \in G, \forall g \in G/H, \quad x_{\sigma,\tau}(g) = a_{\gamma(g,\sigma),\gamma(g\ol{\sigma},\tau)}.
        \end{equation*}
    Alors $x$ est un $2$-cocycle. De plus, pour tous $\sigma,\tau \in H$ et $g \in G/H$, on a $x_{\sigma,\tau}(g) = (\tensor[^{u(g)^{-1}}]{a)}{_{\sigma,\tau}}$. En particulier, $\sh(x) = a$.	
\end{lem}
\begin{proof}
    Pour tous $g \in G/H$ et $\sigma,\tau,\upsilon \in G$, on a
	\begin{align*}
		& \tensor[^\sigma]{x}{_{\tau,\upsilon}}(g) - x_{\sigma\tau,\upsilon}(g) + x_{\sigma,\tau\upsilon}(g) - x_{\sigma,\tau}(g) \\
		= \ & x_{\tau,\upsilon}(g\ol{\sigma}) - x_{\sigma\tau,\upsilon}(g) + x_{\sigma,\tau\upsilon}(g) - x_{\sigma,\tau}(g) \\
		= \ & a_{\gamma(g\ol{\sigma},\tau),\gamma(g\ol{\sigma\tau},\upsilon)} - a_{\gamma(g,\sigma\tau),\gamma(g\ol{\sigma\tau},\upsilon)} + a_{\gamma(g,\sigma),\gamma(g\ol{\sigma},\tau\upsilon)} - a_{\gamma(g,\sigma),\gamma(g\ol{\sigma},\tau)} \\
		= \ & a_{\gamma(g\ol{\sigma},\tau),\gamma(g\ol{\sigma\tau},\upsilon)} - a_{\gamma(g,\sigma)\gamma(g\ol{\sigma},\tau),\gamma(g\ol{\sigma\tau},\upsilon)} + a_{\gamma(g,\sigma),\gamma(g\ol{\sigma},\tau)\gamma(g\ol{\sigma\tau},\upsilon)} - a_{\gamma(g,\sigma),\gamma(g\ol{\sigma},\tau)}, & \text{par \eqref{eqCocycleConditionForGamma},}\\
		= \ & 0, & \text{car } a \text{ est un cocycle,}
    \end{align*} 
    d'o\`u $\tensor[^\sigma]{x}{_{\tau,\upsilon}} - x_{\sigma\tau,\upsilon} + x_{\sigma,\tau\upsilon} - x_{\sigma,\tau} = 0$ et donc $x$ est bien un $2$-cocycle.
	
    Soient $\sigma,\tau \in H$ et $g \in G/H$. Alors $\gamma(g,\sigma) = u(g) \sigma u(g\ol{\sigma})^{-1} = u(g)\sigma u(g)^{-1}$ et $\gamma(g\ol{\sigma},\tau) =  \gamma(g,\tau) = u(g)\tau u(g)^{-1}$, d'o\`u 
	\begin{equation*}
            x_{\sigma,\tau}(g) = a_{u(g)\sigma u(g)^{-1},u(g)\tau u(g)^{-1}}= (\tensor[^{u(g)^{-1}}]{a)}{_{\sigma,\tau}}.
        \end{equation*}
    En particulier, $x_{\sigma,\tau}(1) = a_{\sigma,\tau}$ puisque $u(1) = 1$ et donc $\sh(x) = a$.
\end{proof}

\subsection{Description globale} \label{subsection22} 

\`A partir de ce paragraphe, on identifie $M \otimes M$ au groupe des applications $G/H \times G/H \to A \otimes A$, muni de l'action de $G$ d\'efinie par
    \begin{equation*}
        \forall \sigma \in G,\forall f \in M \otimes M,\forall g,h \in G/H, \quad (\tensor[^\sigma]{f)}{}(g,h) = f(g\ol{\sigma},h\ol{\sigma}).
    \end{equation*}
Alors pour tous $f_1,f_2 \in M$, l'\'el\'ement $f_1 \otimes f_2 \in M \otimes M$ est donn\'e par
    \begin{equation*}
        \forall g,h \in G/H, \quad (f_1 \otimes f_2)(g,h) = f_1(g) \otimes f_2(h).
    \end{equation*}
Pour tout $g \in G/H$, on d\'efinit le morphisme
    \begin{equation} \label{eqOmegaG}
	\omega_g: M \otimes M \to \Ind_G^H (A \otimes A) = \{G/H \to A \otimes A\}, \quad \omega_g(f)(h) = f(gh,h).
    \end{equation}
On v\'erifie sans peine que $\omega_g$ est $G$-\'equivariant. Soit $\omega = (\omega_g)_{g \in G/H}: M \otimes M \to (\Ind_G^H (A \otimes A))^{[G:H]}$. Explicitons la compos\'ee $\sh' = \sh^{[G:H]} \circ \omega_\ast: \Z^r(G, M \otimes M) \to \Z^r(H,A \otimes A)^{[G:H]}$ pour tout $r \ge 1$.

\begin{lem} \label{lemShapiroTensor}
    Soit $r \ge 1$.
        \begin{enumerate}
            \item $\omega$ est un isomorphisme de $G$-modules (donc $\sh': \H^r(G,M\otimes M) \to \H^r(H,A \otimes A)^{[G:H]}$ est un isomorphisme).
            
            \item Soit $x \in \Z^r(G, M \otimes M)$. Pour tout $g \in G/H$, soit $a_g \in \Z^r(H, A \otimes A)$ le $r$-cocycle $(\sigma_1,\ldots,\sigma_r) \mapsto x_{\sigma_1,\ldots,\sigma_r}(g,1)$. Alors $\sh'(x) = (a_g)_{g \in G/H}$.
	\end{enumerate}
\end{lem}
\begin{proof}
\begin{enumerate}
    \item Montrons l'injectivit\'e de $\omega$. Soit $f \in M \otimes M$ tel que $\omega(f) = 0$. Alors 
	\begin{equation*}
            f(g,h) = f(gh^{-1}h,h) = \omega_{gh^{-1}}(f)(h) = 0
        \end{equation*}
    pour tous $g,h \in G/H$, d'o\`u $f = 0$.
		
    Montrons maintenant que $\omega$ est surjectif. Soit alors $(f_g)_{g \in G/H}$ une famille d'applications $G/H \to A \otimes A$. On d\'efinit l'\'el\'ement $f \in M \otimes M$ par
	\begin{equation*}
            \forall g,h \in G/H, \quad f(g,h) = f_{gh^{-1}}(h).
        \end{equation*}
    Pour tous $g,h \in G/H$, on a $\omega_g(f)(h) = f(gh,h) = f_{ghh^{-1}}(h) = f_g(h)$, d'o\`u $\omega_g(f) = f_g$ et donc $\omega(f) = (f_g)_{g \in G/H}$.
		
    \item Pour tous $\sigma_1,\ldots,\sigma_r \in H$ et $g \in G/H$, on a
	\begin{equation*}
            ((\omega_g)_\ast x)_{\sigma_1,\ldots,\sigma_r}(1) = \omega_g(x_{\sigma_1,\ldots,\sigma_r})(1) = x_{\sigma_1,\ldots,\sigma_r}(g,1) = (a_g)_{\sigma_1,\ldots,\sigma_r},
        \end{equation*}
    donc $\sh((\omega_g)_\ast x) = a_g$, d'o\`u $\sh'(x) = (\sh((\omega_g)_\ast x))_{g \in G/H} =  (a_g)_{g \in G/H}$.
\end{enumerate} 
\end{proof}

D\'ecrivons l'application $\H^1(H,A) \times \H^1(H,A) \to \H^2(H,A \otimes A)^{[G:H]}$ associ\'ee au cup-produit 
	\begin{equation*}
            \cup: \H^1(G,M) \times \H^1(G,M) \to \H^2(G, M \otimes M).
        \end{equation*}

\begin{lem} \label{lemShapiroCupProduct}
    On a un diagramme commutatif, o\`u les fl\`eches verticales sont des isomorphismes:
	\begin{equation*}
        \xymatrix{
            \H^1(H,A) \times \H^1(H,A)  \ar[rrrr]^-{(a,b) \mapsto (\tensor[^{g^{-1}}]{a}{} \cup b)_{g \in G/H}} &&&& \H^2(H,A \otimes A)^{[G:H]} \\
            \H^1(G,M) \ar@<6ex>[u]_{\sh} \times \H^1(G,M) \ar@<-6ex>[u]_{\sh} \ar[rrrr]^-{\cup} &&&& \H^2(G, M \otimes M). \ar[u]_{\sh'}
	}
        \end{equation*}
\end{lem} 
\begin{proof}
    Soient $a,b \in \H^1(H,A)$, vus comme morphismes continus $H \to A$. Soient $x,y: G \to M$ les $1$-cocycles repr\'esentant les classes $\sh^{-1}(a), \sh^{-1}(b) \in \H^1(G,M)$ respectivement, qui sont construits dans le lemme \ref{lemShapiro1}. Alors l'on a $x_\sigma(g) = (\tensor[^{u(g)^{-1}}]{a)}{_\sigma}$ et $y_\sigma(g) = (\tensor[^{u(g)^{-1}}]{b)}{_\sigma}$ pour tous $\sigma \in H$ et $g \in G/H$ (o\`u $u(g) \in G$ d\'esigne un relev\'e de $g$).
	
    Regardons le cup-produit $x \cup y \in \Z^2(G, M \otimes M)$. Pour tous $g \in G/H$ et $\sigma,\tau \in H$, on a 
	\begin{equation*}
            (x \cup y)_{\sigma,\tau}(g,1) = (x_\sigma \otimes \tensor[^\sigma]{y}{_\tau})(g,1) = (x_\sigma \otimes y_\tau)(g,1) = x_\sigma(g) \otimes y_\tau(1) = (\tensor[^{u(g)^{-1}}]{a)}{_\sigma} \otimes b_\tau = (\tensor[^{u(g)^{-1}}]{a}{} \cup b)_{\sigma,\tau}.
        \end{equation*}
    Par le lemme \ref{lemShapiroTensor}, on a $\sh'(x \cup y) = (\tensor[^{u(g)^{-1}}]{a}{} \cup b)_{g \in G/H}$ dans $\Z^2(H,A \otimes A)^{[G:H]}$, ce qui implique que $\sh'([x] \cup [y]) = (\tensor[^{g^{-1}}]{a}{} \cup b)_{g \in G/H}$ dans $\H^2(G,A \otimes A)^{[G:H]}$, d'o\`u le lemme.
\end{proof}

\begin{lem} \label{lemShapiroJ}
    On munit $A$ de l'action triviale de $G$. Soit $j: A \otimes A \hookrightarrow M \otimes M$ l'inclusion $G$-\'equivariante qui \`a tout $m \in A \otimes A$ associe l'application
	\begin{equation*}
            G/H \times G/H \to A \otimes A, \quad (g,h) \mapsto m.
        \end{equation*}
    Soit $r \ge 1$ et soit $\res: \Z^r(G,A \otimes A) \to \Z^r(H,A \otimes A)$ le morphisme de restriction. Alors l'on dispose d'un diagramme commmutatif, o\`u la fl\`eche verticale est un isomorphisme:
	\begin{equation*}
            \xymatrix{
                \H^r(G,A \otimes A) \ar[rr]^-{(\res,\ldots,\res)} \ar[rrd]^{j_\ast} && \H^r(H,A \otimes A)^{[G:H]} \\
                && \H^r(G, M \otimes M). \ar[u]_{\sh'}
	    }
        \end{equation*}
\end{lem} 
\begin{proof}
    Soit $x \in \Z^r(G,A \otimes A)$. Alors $(j_\ast x)_{\sigma_1,\ldots,\sigma_r}(g,1) = j(x_{\sigma_1,\ldots,\sigma_r})(g,1) = x_{\sigma_1,\ldots,\sigma_r}$ pour tous $\sigma_1,\ldots,\sigma_r \in H$ et $g \in G/H$. Par le lemme \ref{lemShapiroTensor}, on a $\sh'(j_\ast x) = (\res(x),\ldots,\res(x))$ dans $\Z^r(H,A \otimes A)^{[G:H]}$ et donc $\sh'(j_\ast[x]) = (\res([x]),\ldots,\res([x]))$ dans $\H^r(H,A \otimes A)^{[G:H]}$.
\end{proof}

\subsection{Description locale} \label{subsection23}

Dans ce paragraphe et celui qui le suit, on prend $G = \Gamma_k$ et $H = \Gamma_L$, o\`u $L/k$ est une extension finie galoisienne de corps de nombres. Posons $\gfrak=\Gal(L/k) = \Gamma_k/\Gamma_L$. $A$ est toujours groupe ab\'elien muni de l'action triviale de $\Gamma_L$ et $M = \Ind_{\Gamma_k}^{\Gamma_L} A$. 

On fixe une place $v$ de $k$, une place $w$ de $L$ divisant $v$ et une extension de $w$ \`a $\ol{k}$ (d'o\`u une inclusion $\Gamma_{k_v} \subseteq \Gamma_k$ avec $\Gamma_{L_w} = \Gamma_{k_v} \cap \Gamma_L$). Soit $\gfrak_v = \Gal(L_w/k_v) = \Gamma_{k_v}/\Gamma_{L_w} \subseteq \gfrak$ le groupe de d\'ecomposition de $w|v$. On choisit un syst\`eme de repr\'esentants $\Escr_v$ des classes \`a gauche de $\gfrak$ suivantes $\gfrak_v$ et l'on pose $e_v=|\Escr_v| = [\gfrak:\gfrak_v]$. Alors toute place de $L$ divisant $v$ est de la forme $sw$ pour un unique $s \in \Escr_v$.

Posons $M_v = \Ind_{\Gamma_{k_v}}^{\Gamma_{L_w}} A = \{\gfrak_v \to A\}$. Pour tout $s \in \Escr_v$, soit $\varsigma_s: M \to M_v$ le morphisme $\Gamma_{k_v}$-\'equivariant d\'efini par
	\begin{equation*}
            \forall f \in M,\forall h \in \gfrak_v, \quad \varsigma_s(f)(h) = f(sh).
        \end{equation*}
Soit $\varsigma = (\varsigma_s)_{s \in \Escr_v}: M \to M_v^{e_v}$ et soit $\sh_v=\sh^{e_v} \circ \varsigma_\ast: \Z^r(k_v,M) \to \Z^r(L_w,A)^{e_v}$ pour tout $r \ge 1$.

\begin{lem} \label{lemShapiroLocalM}
    $\varsigma$ est un isomorphisme de $\Gamma_{k_v}$-modules (donc $\sh_v: \H^r(k_v,M) \to \H^r(L_w,A)^{e_v}$ est un isomorphisme pour tout $r \ge 1$).
\end{lem}
\begin{proof}
    Montrons l'injectivit\'e de $\varsigma$. Soit $f \in M$ tel que $\varsigma(f) = 0$. Pour tout $g \in \gfrak$, on peut \'ecrire $g = sh$, o\`u $s \in \Escr_v$ et $h \in \gfrak_v$. Alors $f(g) = f(sh) = \varsigma_s(f)(h) = 0$, d'o\`u $f = 0$.
	
    Montrons maintenant que $\varsigma$ est surjectif. Soit alors $(f_s)_{s \in \Escr_v}$ une famille d'\'el\'ements de $M_v$. On d\'efinit l'\'el\'ement $f \in M$ comme suit. Pour tout $g \in \gfrak$, il existe un unique $s \in \Escr_v$ et un unique $h \in \gfrak_v$ tel que $g = sh$. On d\'efinit $f(g) = f_s(h)$. Alors l'on a $\varsigma_s(f)(h) = f(sh) = f_s(h)$ pour tous $s \in \Escr_v$ et $h \in \gfrak$, d'o\`u $\varsigma_s(f) = f_s$ et donc $\varsigma(f) = (f_s)_{s \in \Escr_v}$.
\end{proof}

Par le lemme \ref{lemShapiroLocalM}, on dispose d'un isomorphisme
    \begin{equation*}
        \varsigma \otimes \varsigma = (\varsigma_s \otimes \varsigma_t)_{s,t \in \Escr_v}: M \otimes M \to M_v^{e_v} \otimes M_v^{e_v} = (M_v \otimes M_v)^{e_v^2}
    \end{equation*}
de $\Gamma_{k_v}$-modules. Notons que pour tous $s,t \in \Escr_v$, $f_1, f_2 \in M$ et $h_1,h_2 \in \gfrak_v$, on a 
    \begin{align*}
        (\varsigma_s \otimes \varsigma_t)(f_1 \otimes f_2)(h_1,h_2) & = (\varsigma_s(f_1) \otimes \varsigma_t(f_2))(h_1,h_2) \\
        & = \varsigma_s(f_1)(h_1) \otimes \varsigma_t(f_2)(h_2) \\
        & = f_1(sh_1) \otimes f_2(th_2) \\
        & = (f_1 \otimes f_2)(sh_1,th_2).
    \end{align*}
Mais comme $M \otimes M$ est engendr\'e (en tant que groupe ab\'elien) par les \'el\'ements de la forme $f_1 \otimes f_2$ (o\`u $f_1,f_2 \in M$), on voit que
    \begin{equation} \label{eqSigmaSTensorSigmaT}
	\forall s,t \in \Escr_v, \forall f \in M \otimes M, \forall h_1,h_2 \in \gfrak_v, \quad (\varsigma_s \otimes \varsigma_t)(f)(h_1,h_2) = f(sh_1,th_2).
    \end{equation}
Par le lemme \ref{lemShapiroTensor} (appliqu\'e aux groupes profinis $\Gamma_{L_w} \subseteq \Gamma_{k_v}$ et au $\Gamma_{k_v}$-module $M_v$), la compos\'ee 
    \begin{equation*}
        \sh'_v=(\sh')^{e_v^2} \circ (\varsigma \otimes \varsigma)_\ast: \Z^r(k_v,M \otimes M) \to \Z^r(L_w,A \otimes A)^{e_v^2 |\gfrak_v|}
    \end{equation*}
induit un isomorphisme $\H^r(k_v,M \otimes M) \to \H^r(L_w,A \otimes A)^{e_v^2 |\gfrak_v|}$ pour tout $r \ge 1$.

Explicitons maintenant l'application $\H^1(L_w,A)^{e_v} \times \H^1(L_w,A)^{e_v} \to \H^2(L_w,A \otimes A)^{e_v^2 |\gfrak_v|}$ associ\'ee au cup-produit $\cup: \H^1(k_v,M) \times \H^1(k_v,M) \to \H^2(k_v, M \otimes M).$

\begin{lem} \label{lemShapiroCupProductLocal}
    On a un diagramme commutatif, o\`u les fl\`eches verticales sont des isomorphismes:
	\begin{equation*}
            \xymatrix{
                \H^1(L_w,A)^{e_v} \times \H^1(L_w,A)^{e_v}  \ar[rrrrrr]^-{((a_s)_{s \in \Escr_v},(b_s)_{s \in \Escr_v}) \mapsto (\tensor[^{h^{-1}}]{a}{_s} \cup b_t)_{h \in \gfrak_v,s,t \in \Escr_v}} &&&&&& \H^2(L_w,A \otimes A)^{e_v^2|\gfrak_v|} \\
                \H^1(k_v,M) \ar@<6ex>[u]_{\sh_v} \times \H^1(k_v,M) \ar@<-6ex>[u]_{\sh_v} \ar[rrrrrr]^-{\cup} &&&&&& \H^2(k_v, M \otimes M). \ar[u]_{\sh'_v}
	    }
        \end{equation*}
\end{lem}
\begin{proof}
    Puisque $\sh_v = \sh^{e_v} \circ \varsigma_\ast$ et $\sh'_v = (\sh')^{e_v^2} \circ (\varsigma \otimes \varsigma)_\ast$, il suffit de montrer que les deux petits carr\'es du diagramme suivant sont commutatifs:
	\begin{equation} \label{eqlemShapiroCupProductLocal}
		\xymatrix{
				\H^1(L_w,A)^{e_v} \times \H^1(L_w,A)^{e_v}  \ar[rrrrrr]^-{((a_s)_{s \in \Escr_v},(b_s)_{s \in \Escr_v}) \mapsto (\tensor[^{h^{-1}}]{a}{_s} \cup b_t)_{h \in \gfrak_v,s,t \in \Escr_v}} &&&&&& \H^2(L_w,A \otimes A)^{e_v^2|\gfrak_v|} \\
				\H^1(k_v,M_v)^{e_v} \ar@<6ex>[u]_{\sh^{e_v}} \times \H^1(k_v,M_v)^{e_v} \ar@<-6ex>[u]_{\sh^{e_v}} \ar[rrrrrr]^-{((x_s)_{s \in \Escr_v},(y_s)_{s \in \Escr_v}) \mapsto (x_s \cup y_t)_{s,t \in \Escr_v}} &&&&&& \H^2(k_v, M_v \otimes M_v)^{e_v^2}. \ar[u]_{(\sh')^{e_v^2}}\\			
				\H^1(k_v,M) \ar@<6ex>[u]_{\varsigma_\ast} \times \H^1(k_v,M) \ar@<-6ex>[u]_{\varsigma_\ast} \ar[rrrrrr]^-{\cup} &&&&&& \H^2(k_v, M \otimes M). \ar[u]_{(\varsigma \otimes \varsigma)_{\ast}}
		}
	\end{equation}
    Par le lemme \ref{lemShapiroCupProduct}, le diagramme
	\begin{equation*}
            \xymatrix{
                \H^1(L_w,A) \times \H^1(L_w,A)  \ar[rrrr]^-{(a,b) \mapsto (\tensor[^{h^{-1}}]{a}{} \cup b)_{h \in \gfrak_v}} &&&& \H^2(L_w,A \otimes A)^{|\gfrak_v|} \\
                \H^1(k_v,M_v) \ar@<6ex>[u]_{\sh} \times \H^1(k_v,M_v) \ar@<-6ex>[u]_{\sh} \ar[rrrr]^-{\cup} &&&& \H^2(k_v, M_v \otimes M_v) \ar[u]_{\sh'}
		}
        \end{equation*}
    est commutatif, d'o\`u la commutativit\'e du carr\'e du haut de \eqref{eqlemShapiroCupProductLocal}. Le carr\'e du bas de \eqref{eqlemShapiroCupProductLocal} commute tout simplement par fonctorialit\'e du cup-produit.
\end{proof}

Le r\'esultat suivant est une version locale du lemme \ref{lemShapiroJ}.

\begin{lem} \label{lemShapiroJLocal}
    On munit $A$ de l'action triviale de $\Gamma_k$. Soit $j: A \otimes A \hookrightarrow M \otimes M$ l'inclusion $\Gamma_k$-\'equivariante qui \`a tout $m \in A \otimes A$ associe l'application
	\begin{equation*}
            \gfrak \times \gfrak \to A \otimes A, \quad (g,h) \mapsto m.
        \end{equation*}
    Soit $r \ge 1$ et soit $\res: \Z^r(k_v,A \otimes A) \to \Z^r(L_w,A \otimes A)$ le morphisme de restriction. Alors l'on dispose d'un diagramme commmutatif, o\`u la fl\`eche verticale est un isomorphisme:
	\begin{equation*}
            \xymatrix{
                \H^r(k_v,A \otimes A) \ar[rr]^-{(\res,\ldots,\res)} \ar[rrd]^{j_\ast} && \H^r(L_w,A \otimes A)^{e_v^2|\gfrak_v|} \\
                && \H^r(k_v, M \otimes M). \ar[u]_{\sh'_v}
	    }
        \end{equation*}
\end{lem} 
\begin{proof}
    Soit $x \in \Z^r(k_v,A \otimes A)$. Pour tous $s, t \in \Escr_v$, $\sigma_1,\ldots,\sigma_r \in \Gamma_{L_w}$ et $h \in \gfrak_v$, on a
	\begin{align*}
            ((\varsigma_s \otimes \varsigma_t)_\ast j_\ast x)_{\sigma_1,\ldots,\sigma_r}(h,1) & = (\varsigma_s \otimes \varsigma_t)(j(x_{\sigma_1,\ldots,\sigma_r}))(h,1) \\
            & = j(x_{\sigma_1,\ldots,\sigma_r})(sh,t), & \text{par \eqref{eqSigmaSTensorSigmaT}}, \\
            & = x_{\sigma_1,\ldots,\sigma_r} \\
            & = \res(x)_{\sigma_1,\ldots,\sigma_r},
	\end{align*}
    d'o\`u $\sh'((\varsigma_s \otimes \varsigma_t)_\ast j_\ast x) = (\res(x),\ldots,\res(x)) \in \Z^r(L_w,A \otimes A)^{|\gfrak_v|}$ pour tous $s,t \in \Escr_v$ en vertu du lemme \ref{lemShapiroTensor}. Ainsi, on a $\sh'_v(j_\ast x) = (\sh')^{e_v^2}((\varsigma \otimes \varsigma)_\ast j_\ast x) = (\res(x),\ldots,\res(x))$ dans $\Z^r(L_w,A \otimes A)^{e_v^2|\gfrak_v|}$, ou $\sh'_v(j_\ast [x]) = (\res([x]),\ldots,\res([x]))$ dans $\H^r(L_w,A \otimes A)^{e_v^2|\gfrak_v|}$, ce qui ach\`eve la d\'emonstration.
\end{proof}

\subsection{Description des localisations} \label{subsection24}

On garde les notations du paragraphe pr\'ec\'edent. Rappelons qu'on a des isomorphismes
    \begin{equation*}
        \sh: \H^1(k,M) \to \H^1(L,A)
    \end{equation*}
et
    \begin{equation*}
        \sh_v = \sh^{e_v} \otimes \varsigma_\ast: \H^1(k_v,M) \to \H^1(L_w,A)^{e_v}.
    \end{equation*}
D\'ecrivons le morphisme $\H^1(L,A) \to \H^1(L_w,A)^{e_v}$ associ\'e au morphisme de localisation 
    \begin{equation*}
        \loc_v: \H^1(k,M) \to \H^1(k_v,M).
    \end{equation*}

\begin{lem} \label{lemShapiroLocalization1}
    On a un diagramme commutatif, o\`u les fl\`eches verticales sont des isomorphismes:
	\begin{equation*}
        \xymatrix{
            \H^1(L,A) \ar[rrrr]^-{a \mapsto ((\tensor[^{s^{-1}}]{a)}{_w})_{s \in \Escr_v}} &&&& \H^1(L_w,A)^{e_v} \\
            \H^1(k,M) \ar[rrrr]^{\loc_v} \ar[u]_{\sh} &&&& \H^1(k_v,M). \ar[u]_{\sh_v}
	}
        \end{equation*}
\end{lem}
\begin{proof}
    Soit $a \in \H^1(L,A)$, vu comme morphisme continu $\Gamma_L \to A$. Soit $x: \Gamma_k \to M$ le $1$-cocycle repr\'esentant la classe $\sh^{-1}(a) \in \H^1(k,M)$ construit dans le lemme \ref{lemShapiro1}. Alors l'on a $x_\sigma(g) = (\tensor[^{u(g)^{-1}}]{a)}{_\sigma}$ pour tous $\sigma \in \Gamma_L$ et $g \in \gfrak$ (o\`u  $u(g) \in \Gamma_k$ d\'esigne un relev\'e de $g$).
	
    \'Etudions le localis\'e $x_v: \Gamma_{k_v} \to M$. Pour tous $s \in \Escr_v$ et $\sigma \in \Gamma_{L_w}$, on a
        \begin{equation*}
            ((\varsigma_s)_\ast x_v)_\sigma(1) = \varsigma_s((x_v)_\sigma)(1) = (x_v)_\sigma(s) = x_\sigma(s) = (\tensor[^{u(s)^{-1}}]{a)}{_\sigma} = ((\tensor[^{u(s)^{-1}}]{a)}{_w})_{\sigma}.
        \end{equation*}
    Ainsi, $\sh((\varsigma_s)_\ast x_v) = (\tensor[^{u(s)^{-1}}]{a)}{_w}$ dans $\Z^1(L_w,A)$, d'o\`u $\sh_v(x_v) = ((\varsigma_s)_\ast x_v)_{s \in \Escr_v} = ((\tensor[^{u(s)^{-1}}]{a)}{_w})_{s \in \Escr_v}$ dans $\Z^1(L_w,A)^{e_v}$, ou $\sh_v([x_v]) = ((\tensor[^{s^{-1}}]{a)}{_w})_{s \in \Escr_v}$ dans $\H^1(L_w,A)^{e_v}$.
\end{proof}

Rappelons maintenant qu'on a des isomorphismes
    \begin{equation*}
        \sh': \H^2(k,M \otimes M) \to \H^2(L,A \otimes A)^{|\gfrak|}
    \end{equation*}
et
    \begin{equation*}
        \sh'_v = (\sh')^{e_v^2} \circ (\varsigma \otimes \varsigma)_\ast: \H^2(k_v,M \otimes M) \to \H^2(L_w,A \otimes A)^{e_v^2|\gfrak_v|}.
    \end{equation*}
D\'ecrivons le morphisme $\H^2(L,A \otimes A)^{|\gfrak|} \to \H^2(L_w,A \otimes A)^{e_v^2|\gfrak_v|}$ associ\'e au morphisme de localisation 
    \begin{equation*}
        \loc_v: \H^2(k,M \otimes M) \to \H^2(k_v,M \otimes M).
    \end{equation*}

\begin{lem} \label{lemShapiroLocalization2}
	On a un diagramme commutatif, o\`u les fl\`eches verticales sont des isomorphismes:
	\begin{equation*}
        \xymatrix{
		\H^2(L,A \otimes A)^{|\gfrak|} \ar[rrrrrr]^-{(\alpha_g)_{g \in \gfrak} \mapsto ((\tensor[^{t^{-1}}]{\alpha}{_{sht^{-1}}})_w)_{h \in \gfrak_v,s,t \in \Escr_v}} &&&&&& \H^2(L_w,A \otimes A)^{e_v^2 |\gfrak_v|} \\
		\H^2(k,M \otimes M) \ar[rrrrrr]^{\loc_v} \ar[u]_{\sh'} &&&&&& \H^2(k_v,M \otimes M). \ar[u]_{\sh'_v}
	}
        \end{equation*}
\end{lem}
\begin{proof}
    Soit $(a_g)_{g \in \gfrak}$ une famille de $2$-cocycles $\Gamma_L \times \Gamma_L \to A \otimes A$. Pour tout $g \in \gfrak$, soit $x_g: \Gamma_k \times \Gamma_k \to \Ind_{\Gamma_k}^{\Gamma_L} (A \otimes A) = \{\gfrak \to A \otimes A\}$ le $2$-cocycle construit dans le lemme \ref{lemShapiro2} (appliqu\'e au $A \otimes A$ au lieu de $A$), de sorte que $\sh(x_g) = a_g$ et que
	\begin{equation} \label{eqlemShapiroLocalization2}
            \forall \sigma,\tau \in \Gamma_L,\forall h \in \gfrak, \quad (x_g)_{\sigma,\tau}(h) = (\tensor[^{u(h)^{-1}}]{a}{_g})_{\sigma,\tau},
	\end{equation}
    o\`u $u(h) \in \Gamma_k$ d\'esigne un relev\'e de $h$. 
	
    Reprenons les notations $\omega_g: M \otimes M \to \Ind_{\Gamma_k}^{\Gamma_L} (A \otimes A)$ et
	\begin{equation*}
            \omega = (\omega_g)_{g \in \gfrak}: M \otimes M \to (\Ind_{\Gamma_k}^{\Gamma_L} (A \otimes A))^{|\gfrak|}
        \end{equation*}
    au d\'ebut du paragraphe \ref{subsection22}. Par le lemme \ref{lemShapiroTensor}, il existe un $2$-cocycle $y \in \Z^2(k,M \otimes M)$ tel que $x_g = (\omega_g)_\ast y$ pour tout $g \in \gfrak$. Alors $\sh'(y) = (\sh((\omega_g)_\ast y))_{g \in \gfrak} = (\sh(x_g))_{g \in \gfrak} = (a_g)_{g \in \gfrak} \in \Z^2(L,A \otimes A)^{|\gfrak|}$.
	
    \'Etudions le localis\'e $y_v \in \Z^2(k_v,M \otimes M)$. Pour tous $s,t \in \Escr_v$, $h \in \gfrak_v$ et $\sigma,\tau \in \Gamma_{L_w}$, on a
	\begin{align*} 
            ((\varsigma_s \otimes \varsigma_t)_\ast y_v)_{\sigma,\tau}(h,1) & = (\varsigma_s \otimes \varsigma_t)((y_v)_{\sigma,\tau})(h,1) \\
            & = (\varsigma_s \otimes \varsigma_t)(y_{\sigma,\tau})(h,1) \\
            & = y_{\sigma,\tau}(sh,t), & \text{par \eqref{eqSigmaSTensorSigmaT}}, \\
            & = y_{\sigma,\tau}(sht^{-1}t,t) \\
            & = \omega_{sht^{-1}}(y_{\sigma,\tau})(t), & \text{par \eqref{eqOmegaG}},\\
            & = ((\omega_{sht^{-1}})_\ast y)_{\sigma,\tau}(t) \\
            & = (x_{sht^{-1}})_{\sigma,\tau}(t) \\
            & = (\tensor[^{u(t)^{-1}}]{a}{_{sht^{-1}}})_{\sigma,\tau}, & \text{par \eqref{eqlemShapiroLocalization2}}, \\
            & = ((\tensor[^{u(t)^{-1}}]{a}{_{sht^{-1}}})_w)_{\sigma,\tau},
	\end{align*}
    d'o\`u $\sh'((\varsigma_s \otimes \varsigma_t)_\ast y_v) = ((\tensor[^{u(t)^{-1}}]{a}{_{sht^{-1}}})_w)_{h \in \gfrak_v}$ dans $\Z^2(L_w,A \otimes A)^{|\gfrak_v|}$ par le lemme \ref{lemShapiroTensor}. Ainsi
	\begin{equation*}
            \sh'_v(y_v) = (\sh'((\varsigma_s \otimes \varsigma_t)_\ast y_v))_{s,t \in \Escr_v} = ((\tensor[^{u(t)^{-1}}]{a}{_{sht^{-1}}})_w)_{h \in \gfrak_v, s,t \in \Escr_v}
        \end{equation*}
    dans $\Z^2(L_w,A \otimes A)^{e_v^2|\gfrak_v|}$, ou $\sh'_v([y_v]) = ((\tensor[^{t^{-1}}]{[a}{_{sht^{-1}}}])_w)_{h \in \gfrak_v, s,t \in \Escr_v}$ dans $\H^2(L_w,A \otimes A)^{e_v^2|\gfrak_v|}$. Cette \'egalit\'e \'etablit le r\'esultat voulu.
\end{proof} 
\section{La construction de Borovoi--Kunyavski\u{\i}} \label{section3}

\subsection{$\H^2$ non ab\'elien et espaces homog\`enes} \label{subsection31}

Soit $K$ un corps de caract\'eristique nulle. Pour \'etudier les $K$-espaces homog\`enes de $\SL_n$ (qui peut ne pas avoir de $K$-point), il convient d'introduire la notion de $2$-cohomologie non ab\'elienne. Nous allons utiliser sa version concr\`ete en termes de cocycles \cite[\S 1]{FSS}. 

 Seuls les $K$-liens dont le $\ol{K}$-groupe sous-adjacent est fini (c'est donc simplement un groupe abstrait fini) seront pris en consid\'eration. Un {\em $K$-lien (fini)} $L = (F,\kappa)$ est alors la donn\'ee d'un groupe fini $F$ muni d'une action ext\'erieure $\kappa: \Gamma_K \to \Out(F)$. Si $F$ est un $K$-groupe ({\em i.e.} si une action (continue) $\rho: \Gamma_K \to \Aut(F)$ est donn\'ee), on lui associe son $K$-lien canonique $\lien(F)$.

Soit $L = (F,\kappa)$ un $K$-lien. Un {\em $2$-cocycle \`a coefficients dans $L$} est un couple $(\rho, u)$ o\`u $\rho: \Gamma_K \to \Aut(F)$ et $u: \Gamma_K \times \Gamma_K \to F$ sont des applications continues telles que
    \begin{enumerate}
	\item $\rho_\sigma$ rel\`eve $\kappa_\sigma$ pour tout $\sigma \in \Gamma_K$;
		 
	\item $\rho_{\sigma\tau} = \int(u_{\sigma,\tau}) \circ \rho_\sigma \circ \rho_\tau$ pour tous $\sigma,\tau \in \Gamma_K$;
		
	\item $u_{\sigma,\tau\upsilon}\rho_\sigma(u_{\tau,\upsilon}) = u_{\sigma\tau,\upsilon} u_{\sigma,\tau}$ pour tous $\sigma,\tau,\upsilon \in \Gamma_K$.
    \end{enumerate} 
Notons $\Z^2(K,L)$ l'ensemble des $2$-cocycles \`a coefficients dans $L$. Deux tels $2$-cocycles $(\rho,u)$ et $(\rho',u')$ sont dits {\em cohomologues} s'il existe une application continue $c: \Gamma_K \to F$ telle que

\begin{enumerate}
		\item $\rho'_\sigma = \int(c_\sigma) \circ  \rho_\sigma$ pour tout $\sigma \in \Gamma_K$;
		
		\item $u'_{\sigma,\tau} = c_{\sigma\tau} u_{\sigma,\tau} \rho_\sigma(c_\tau)^{-1} c_{\sigma}^{-1}$ pour tous $\sigma,\tau \in \Gamma_K$.
	\end{enumerate}
Alors \og \^etre cohomologues \fg{} est une relation d'\'equivalence sur $\Z^2(K,L)$. On appelle {\em ensemble de $2$-cohomologie galoisienne \`a coefficients dans $L$} le quotient $\Z^2(K,L)$ par cette relation, et on le note $\H^2(K,L)$. Si $F$ est un $K$-groupe, on note $\H^2(K,F) := \H^2(K,\lien(F))$. Dans le cas o\`u $F$ est un $\Gamma_K$-module, $\H^2(K,F)$ est le groupe de $2$-cohomologie galoisienne usuel.

Soit $L = (F,\kappa)$ un $K$-lien. Une classe $\eta \in \H^2(K,L)$ est dite {\em neutre} si elle est repr\'esent\'ee par un $2$-cocycle de la forme $(\rho,1)$. Dans ce cas, $\rho: \Gamma_K \to \Aut(F)$ est une action continue de $\Gamma_K$ sur $F$ relevant $\kappa$ et donc $\eta$ correspond \`a une $K$-forme $F'$; on note $\eta = n(F')$. L'ensemble $\H^2(K,L)$ peut ne pas avoir de classe neutre, et elle peut \'egalement en poss\'eder plusieurs. Si $F$ est un $K$-groupe et $\rho: \Gamma_K \to F$ d\'esigne l'action de $\Gamma_K$, l'ensemble $\H^2(K,F)$ a une classe neutre privil\'egi\'ee $\eta_0 = [(\rho,1)] = n(F)$. Si $F$ est un $\Gamma_K$-module, la seule classe neutre du groupe $\H^2(K,F)$ est $0$.

Soit $L = (F,\kappa)$ un $K$-lien et soit $Z = Z(F)$. Comme $Z$ est commutatif et caract\'eristique dans $F$, l'action ext\'erieure de $\Gamma_K$ sur $F$ induit une action sur $Z$, {\em i.e.} $Z$ est naturellement un $\Gamma_K$-module. Dans le cas o\`u l'ensemble $\H^2(K,L)$ est non vide, c'est un espace principal homog\`ene du groupe ab\'elien $\H^2(K,Z)$, l'action \'etant d\'efinie par
    \begin{equation*}
        [\beta] \cdot [(\rho,u)] := [(\rho,\beta u)],
    \end{equation*}
o\`u $\beta$ (resp. $(\rho,u)$) est un $2$-cocycle \`a coefficients dans $Z$ (resp. dans $L$) \cite[Lemma 1.9]{Borovoi93}.

Soit $X$ un $K$-espace homog\`ene d'un $K$-groupe alg\'ebrique $G$ \`a stabilisateur g\'eom\'etrique fini $F$. Alors $F$ n'est pas {\em a priori} muni d'une action, mais seulement une action ext\'erieure de $\Gamma_K$. On peut donc d\'efinir le {\em lien de Springer} $L_X$ de $X$. Plus concr\`etement, soit $x \in X(\ol{K})$ un point g\'eom\'etrique et soit $F \subseteq G(\ol{K})$ son stabilisateur. Pour  tout $\sigma \in \Gamma_K$, \'ecrivons $\tensor[^\sigma]{x}{} = x \cdot g_\sigma$, o\`u $g_\sigma \in \SL_n(\ol{K})$ est unique modulo multiplication \`a gauche par un \'el\'ement de $F$. On peut choisir les $g_\sigma$ de sorte que l'application $\sigma \mapsto g_\sigma$ est continue. On d\'efinit
    \begin{equation*}
        \rho_\sigma: F \to F, \quad \rho_\sigma(f) = g_\sigma \tensor[^\sigma]{f}{} g_\sigma^{-1}.
    \end{equation*}
Alors $\rho: \Gamma_K \to \Aut(F)$ est continue et la compos\'ee $\kappa: \Gamma_K \xrightarrow{\rho} \Aut(F) \to \Out(F)$ est une action ext\'erieure. On v\'erifie que le $K$-lien $L_X = (F, \kappa)$ ne d\'epend pas du choix de $x$ et des $g_\sigma$; c'est le lien de Springer de $X$. De plus, si l'on note $u_{\sigma,\tau}:=g_{\sigma\tau}\tensor[^\sigma]{g}{}_\tau^{-1}g_{\sigma}^{-1}$ pour tous $\sigma,\tau \in \Gamma_K$, alors $(\rho,u)$ est un $2$-cocycle \`a coefficients dans $L_X$. On d\'efinit la {\em classe de Springer} de $X$ comme \'etant $\eta_X := [(\rho,u)] \in \H^2(K,L_X)$. Elle est neutre si et seulement si $X$ est domin\'e par un $K$-torseur sous $G$ \cite[\S 7.6]{Borovoi93}. En particulier, si $\H^1(K,G) = 1$ (par exemple, c'est le cas pour $G = \GL_n$ ou $G = \SL_n$ par une variante du th\'eor\`eme 90 de Hilbert), alors $\eta_X$ est neutre si et seulement si $X(K) \neq \varnothing$.

\begin{lem} \label{lemExistenceOfHomogeneousSpace}
    Soit $L$ un $K$-lien et soit $\eta \in \H^2(K,L)$. Alors il existe un entier $n$ et un espace homog\`ene $X$ de $\SL_n$ de lien de Springer $L$ et de classe de Springer $\eta$. Deux tels espaces homog\`enes sont $K$-stablement birationnels.
\end{lem}
\begin{proof}
    On pourra consulter \cite[Corollaire 3.3, Corollaire 3.5]{demarche2019reduction}.
\end{proof} 

\begin{rmk}
    Soit $F$ un $K$-groupe fini et $Z = Z(F)$. \`A partir de chaque classe $\eta \in \H^2(K,Z)$, Borovoi et Kunyavski\u{\i} ont construit explicitement un espace homog\`ene $X$ de $\SL_n$ de lien de Springer $\lien(F)$ \cite[\S 2]{Bokun}. Notant $\eta_0 = n(F) \in \H^2(K,F)$ est la classe neutre privil\'egi\'ee, alors la classe de Springer de $X$ est $\eta \cdot \eta_0$ \cite[Lemma 5.3]{HS02}. Puisque $\H^2(K,Z)$ agit transitivement sur $\H^2(K,F)$, cette construction-l\`a donne tous les espaces homog\`enes de $\SL_n$ de lien de Springer $\lien(F)$, \`a \'equivalence birationnelle stable pr\`es.
\end{rmk}

\begin{prop} \label{propSpringerClass}
    Soient $F$ un $K$-groupe fini, $Z = Z(F)$ et $\eta_0 \in \H^2(K,F)$ la classe neutre privil\'egi\'ee. On note $\Delta: \H^1(K,F/Z) \to \H^2(K,Z)$ l'application connectante induite par l'extension centrale ({\em cf.} \cite[Chapitre VII, Appendice, Proposition 2]{Serre2004Locaux})
        \begin{equation*}
            1 \to Z \to F \to F/Z \to 1.
        \end{equation*}
    Soit $\beta \in \H^2(K,Z)$ et soit $X$ un espace homog\`ene de $\SL_n$ de lien de Springer $\lien(F)$ et de classe de Springer $\eta_X = \beta \cdot \eta_0 \in \H^2(K,F)$. Alors $X(K) \neq \varnothing$ si et seulement s'il existe $\alpha \in \H^1(K,F/Z)$ tel que $\beta = \Delta(\alpha)$. Dans ce cas, $X$ est $K$-isomorphe \`a $\tensor[_\afrak]{F}{} \backslash \SL_n$, o\`u $\afrak: \Gamma_K \to F/Z$ est n'importe quel $1$-cocycle repr\'esentant $\alpha$. Ici, $\tensor[_\afrak]{F}{}$ d\'esigne la $K$-forme de $F$ tordue par un $1$-cocycle $\afrak: \Gamma_K \to F/Z$ repr\'esentant $\alpha$, c'est-\`a-dire que son action de Galois $\cdot_\afrak: \Gamma_k \times F^{\afrak} \to F^{\afrak}$ est donn\'ee par       
        \begin{equation*}
            \forall \sigma \in \Gamma_k,\forall f \in F, \quad \sigma \cdot_\afrak f = \tilde{\afrak}_s \tensor[^\sigma]{f}{} \tilde{\afrak}_\sigma^{-1}
        \end{equation*} 
    o\`u $\tilde{\afrak}_\sigma \in F$ est un relev\'e quelconque de $\afrak_\sigma$.
\end{prop} 
\begin{proof}
    Par \cite[Lemma 2.4, Lemma 2.5]{Borovoi93}, $\eta$ est neutre si et seulement s'il existe un $1$-cocycle $\afrak: \Gamma_K \to F/Z$ tel que $\eta_X = n(\tensor[_\afrak]{F}{}) = \Delta([\afrak]) \cdot \eta_0$. Cette \'egalit\'e \'equivaut \`a $\Delta([\afrak]) = \beta$ puisque l'action de $\H^2(K,Z)$ sur $\H^2(K,F)$ est libre. Dans ce cas, $X$ est $k$-stablement birationnel \`a $\tensor[_\afrak]{F}{} \backslash \SL_{n}$, donc $X$ poss\`ede un $K$-point dont le stabilisateur est $K$-isomorphe \`a $\tensor[_\afrak]{F}{}$, {\em i.e.} $X \simeq \tensor[_\afrak]{F}{} \backslash \SL_n$.
\end{proof}

\begin{rmk}
    Lorsque $K = k$ et un corps de nombres, $F$ est un $k$-groupe fini et $Z = Z(F) = [F,F]$, Harari et Skorobogatov ont utilis\'e une \og th\'eorie de la descente non ab\'elienne \fg{} pour d\'emontrer le fait suivant: Si $X$ est un $k$-espace homog\`ene de $\SL_n$ de lien de Springer $\lien(F)$ et de classe de Springer $\eta \cdot \eta_0$, o\`u $\eta \in \Sha^2(k,Z)$, alors $X$ poss\`ede un point ad\'elique orthogonal \`a $\Br_1 X:= \Ker(\Br X \to \Br X_{\ol{K}})$ \cite[Proposition 5.5]{HS02}. Ils ont soulev\'e la question si l'on peut en d\'eduire l'existence des contres-exemples au principe de Hasse non expliqu\'es par l'obstruction de Brauer--Manin alg\'ebrique ({\em cf. loc. cit.}, Remark 5.6). Dans ce texte, nous y donnons une r\'eponse partielle: c'est impossible pour les donn\'ees de $k$ et de $F$ comme dans l'\'enonc\'e du th\'eor\`eme \ref{thmA}.
\end{rmk}

\subsection{Construction du stabilisateur} \label{subsection32}

Soit $K$ un corps de caract\'eristique nulle. Soient $M$ et $Z$ des $\Gamma_K$-modules finis et $\phi: M \otimes M \to Z$ un morphisme $\Gamma_K$-\'equivariant. On dira que $\phi$ est {\em non d\'eg\'en\'er\'e} s'il poss\`ede les propri\'et\'es suivantes.
\begin{enumerate}
	\item Si $x \in M$ est tel que $\phi(x \otimes y) = 0$ pour tout $y \in M$, alors $x = 0$.
	\item Si $y \in M$ est tel que $\phi(x \otimes y) = 0$ pour tout $x \in M$, alors $y = 0$.
\end{enumerate}
On munit $Z$ de l'action triviale de $M \oplus M$. Posons
\begin{equation} \label{eqBigPhi}
	\Phi: (M \oplus M) \times (M \oplus M) \to Z, \quad \Phi((x,y),(x',y')) = \phi(x \otimes y').
\end{equation}
Alors $\Phi$ est biadditive, donc c'est un $2$-cocycle {\em normalis\'e} ({\em i.e.} $\Phi(a,0) = \Phi(0,a) = 0$ pour tout $a \in M \oplus M$). Soit $F$ le produit crois\'e $Z \times_{\Phi} (M \oplus M)$, {\em i.e.} $F = Z \times (M \oplus M)$ ensemblistement et la loi de composition sur $F$ est donn\'ee par la formule
\begin{equation} \label{eqGroupLawOnF}
	\forall z,z' \in Z, \forall a,a' \in M \oplus M, \quad (z,a)(z',a') = (z + z' + \Phi(a,a'), a + a').
\end{equation}
En particulier, $(z,a)=(z,0)(0,a) = (0,a)(z,0)$. De plus, comme 
	\begin{equation*}
            (0,a)(0,-a) = (\Phi(a,-a),0) = (-\Phi(a,a),0) = (\Phi(a,a),0)^{-1},
        \end{equation*}
on a $(0,a)^{-1} = (\Phi(a,a),0)(0,-a) = (\Phi(a,a),-a)$ et donc
    \begin{equation} \label{eqInverseOnF}
	(z,a)^{-1} = (z,0)^{-1}(0,a)^{-1} = (-z,0)(\Phi(a,a),-a) = (\Phi(a,a)-z,-a).
    \end{equation}
On obtient une extension {\em centrale} de groupes {\em abstraits}:
\begin{equation} \label{eqExtension}
	0 \to Z \to F \to M \oplus M \to 0.
\end{equation}

\begin{lem} \label{lemCommutator}
    Soient $M,Z,\phi,\Phi,F$ comme ci-dessus. Soient $z,z' \in Z$ et $a,a' \in M \oplus M$.
        \begin{enumerate}
            \item On a $(z,a)(z',a')(z,a)^{-1} = (z' + \Phi(a,a') - \Phi(a',a),a')$.
		
            \item Le commutateur $[(z,a),(z',a')] = (\Phi(a,a') - \Phi(a',a),0)$.
		
            \item Si $\phi$ est surjectif, $Z = [F,F]$.
	\end{enumerate}
\end{lem}
\begin{proof} 
    \begin{enumerate}
	\item On a
		\begin{align*}
				& (z,a)(z',a')(z,a)^{-1} \\
				= \ & (z,0)(0,a)(z',a')(0,a)^{-1}(z,0)^{-1} \\
				= \ & (0,a)(z',a')(0,a)^{-1}, & \text{comme $Z$ est central dans $F$}, \\
				= \ & (z' + \Phi(a,a'), a + a')(\Phi(a,a),-a), & \text{par \eqref{eqGroupLawOnF} et \eqref{eqInverseOnF}}, \\
				= \ & (z' + \Phi(a,a') + \Phi(a,a) - \Phi(a+a',a), a + a' - a), & \text{par \eqref{eqGroupLawOnF}}, \\
				= \ & (z' + \Phi(a,a') - \Phi(a',a),a').
		\end{align*}
		
	\item On a
		\begin{align*}
				& [(z,a),(z',a')] \\
				= \ & [(z,0)(0,a),(z',0)(0,a')] \\
				= \ & [(0,a),(0,a')], & \text{comme $Z$ est central dans $F$},\\
				= \ & (0,a)(0,a')(0,a)^{-1} (0,a')^{-1} \\
				= \ & (\Phi(a,a') - \Phi(a',a),a') (\Phi(a',a'),-a'), & \text{par le premier point et \eqref{eqInverseOnF}}, \\
				= \ & (\Phi(a,a') - \Phi(a',a) + \Phi(a',a') - \Phi(a',a'), a' - a'), & \text{par \eqref{eqGroupLawOnF}}, \\
				= \ &  (\Phi(a,a') - \Phi(a',a),0).
		\end{align*}
		
	\item Comme $F/Z = M \oplus M$ est commutatif, on a toujours que $[F,F] \subseteq Z$. De plus,
		\begin{equation*}
                \phi(x \otimes y) = \phi(x \otimes y) - \phi(0,0) = \Phi((x,0),(0,y)) - \Phi((0,y),(x,0))
            \end{equation*}
	pour tous $x,y \in M$, d'o\`u
		\begin{equation*}
                (\phi(x \otimes y),0) = (\Phi((x,0),(0,y)) - \Phi((0,y),(x,0)), 0) = [(0,(x,0)),(0,(0,y))] \in [F,F]
            \end{equation*}
	par le deuxi\`eme point. Ainsi, on aura $Z = Z \times \{0\} \subseteq [F,F]$ d\`es que $\phi$ est surjectif.
    \end{enumerate}
\end{proof} 

\begin{lem} \label{lemCenterOfF}
    Soient $M,Z,\phi,\Phi,F$ comme ci-dessus. Si $\phi$ est non d\'eg\'en\'er\'e, $Z = Z(F)$.
\end{lem}
\begin{proof} 		
    Comme $Z$ est centrale dans $F$, il reste \`a montrer que $Z(F) \subseteq Z$. Supposons $(z,a) \in Z(F)$, o\`u $z \in Z$ et $a = (x,y) \in M \oplus M$. En vertu du lemme \ref{lemCommutator}, on a
	\begin{equation*}
            (0,0) = [(z,a),(0,a')] = (\Phi(a,a') - \Phi(a',a),0) = (\phi(x \otimes y') - \phi(x' \otimes y),0)
        \end{equation*}
    pour tout $a' = (x',y') \in M \oplus M$. En choisissant $x' = 0$, on obtient
	\begin{equation*}
            \forall y' \in M, \quad \phi(x \otimes y') = 0,
        \end{equation*}
    d'o\`u $x = 0$ puisque $\phi$ est non d\'eg\'en\'er\'e. De m\^eme, $y = 0$, d'o\`u $a = 0$ et donc $(z,a) = (z,0) \in Z$.
\end{proof} 

Finalement, on munit $F$ d'une action de $\Gamma_K$ qui est compatible avec celles sur $Z$ et sur $M \oplus M$; alors \eqref{eqExtension} devient une extenstion centrale de $K$-groupes finis. Dans ce texte, on se concentre principalement sur l'action {\em coordonn\'ees par coordonn\'ees}, c'est-\`a-dire
    \begin{equation*}
        \forall \sigma \in \Gamma_K,\forall z \in Z,\forall a \in M \oplus M, \quad \tensor[^\sigma]{(z,a)}{} = (\tensor[^\sigma]{z}{},\tensor[^\sigma]{a}{}).
    \end{equation*}
De \eqref{eqGroupLawOnF}, on voit que cette action est bien compatible avec la loi de composition sur $F$ puisque $\Phi$ est $\Gamma_K$-\'equivariant. 

\begin{defn}
    On appelle {\em espace homog\`ene de Borovoi--Kunyavski\u{\i}} tout espace homog\`ene de $\SL_n$ \`a stabilisateur g\'eom\'etrique fini $F$ une extension de la forme $\eqref{eqExtension}$ avec $Z = Z(F) = [F,F]$.
\end{defn}

Rappelons que la condition $Z = Z(F) = [F,F]$ se garantit lorsque $\phi$ est surjectif et non d\'eg\'en\'er\'e (lemmes \ref{lemCommutator} et \ref{lemCenterOfF}).

\subsection{Quelques calculs avec des cocycles} \label{subsection33} 

Soit $K$ un corps de caract\'eristique nulle. Soient $M$ et $Z$ des $\Gamma_K$-modules finis et $\phi: M \otimes M \to Z$ un morphisme $\Gamma_K$-\'equivariant. On d\'efinit l'application biadditive $\Phi: (M \oplus M) \times (M \oplus M) \to Z$ par \eqref{eqBigPhi} et l'on pose $F = Z \times_\Phi (M \oplus M)$, muni de l'action coordonn\'ees par coordonn\'ees de $\Gamma_K$.

Si $x,y : \Gamma_K \to M$ sont des $1$-cocha\^ines, on notera $x \otimes y: \Gamma_K \to M \otimes M$ la $1$-cocha\^ine $\sigma \mapsto x_\sigma \otimes y_\sigma$. Dans ce qui suit, on consid\'erera les cup-produits des cocha\^ines \`a coefficients dans $M$ induit par l'accouplement $\otimes: M \times M \to M \otimes M$. Pour la preuve du th\'eor\`eme \ref{thmHassePrinciple} (c'est-\`a-dire du th\'eor\`eme \ref{thmA}), au vu de la proposition \ref{propSpringerClass}, il conviendra de d\'ecrire l'application connectante $\Delta: \H^1(K,M \oplus M) \to \H^2(K,Z)$ induite par \eqref{eqExtension}.

\begin{lem} \label{lemImageOfDelta}
    Avec les notations ci-dessus, on a les r\'esultats suivants.
        \begin{enumerate}
            \item Soit $\afrak = (\xfrak, \yfrak): \Gamma_K \to M \oplus M$ un $1$-cocycle ({\em i.e.}, $\xfrak$ et $\yfrak$ sont des $1$-cocycles), et soit $\tensor[_\afrak]{F}{}$ la $K$-forme de $F$ tordue par $\afrak$ ({\em cf.} proposition \ref{propSpringerClass}). Soient $z: \Gamma_K \to Z$, $a = (x,y): \Gamma_K \to M \oplus M$ des $1$-cocha\^ines, et $f = (z,a): \Gamma_K \to \tensor[_\afrak]{F}{}$. Alors $f$ est un cocycle si et seulement si
		\begin{itemize}
                \item $a$ est un cocycle;
                \item $\diff z + \phi_\ast(\xfrak \cup y + x \cup \yfrak + x \cup y + \diff(x \otimes \yfrak)) = 0$.
		\end{itemize}
	
        \item Notons $\Delta: \H^1(K,M \oplus M) \to \H^2(K,Z)$ l'application connectante induite par \eqref{eqExtension} ({\em cf.} \cite[Chapitre VII, Appendice, Proposition 2]{Serre2004Locaux}). Alors $\Delta([a]) = \phi_\ast[x \cup y]$ pour tout $1$-cocycle $a = (x,y): \Gamma_K \to M \oplus M$.
    \end{enumerate}
\end{lem}
\begin{proof}
\begin{enumerate}
    \item Notons $\cdot_\afrak$ l'action de $\Gamma_K$ sur $\tensor[_\afrak]{F}{}$. Alors
        \begin{equation} \label{eqlemImageOfDelta}
            \sigma \cdot_\afrak (\zeta,\alpha )   = (0,\afrak_\sigma)  (\tensor[^\sigma]{\zeta}{},\tensor[^\sigma]{\alpha}{})(0,\afrak_\sigma)^{-1} = (\tensor[^\sigma]{\zeta}{} + \Phi(\afrak_\sigma,\tensor[^\sigma]{\alpha}{}) - \Phi(\tensor[^\sigma]{\alpha }{},\afrak_\sigma), \tensor[^\sigma]{\alpha}{})
	\end{equation}
    pour tous $(\zeta,\alpha) \in F$ et $\sigma \in \Gamma_K$, en vertu du lemme \ref{lemCommutator}.

    Le morphisme $\tensor[_\afrak]{F}{} \to M \oplus M$ \'etant la deuxi\`eme projection, on voit que $a$ est forc\'ement un $1$-cocycle d\`es que $f$ l'est. Soient $\sigma,\tau \in \Gamma_K$. Sous la condition que $a$ est un cocycle, on a
	\begin{align*}
            & (0,a_\sigma)(\sigma \cdot_\afrak (z_\tau,a_\tau))(0,a_{\sigma\tau})^{-1}\\
		= \ & (0,a_\sigma) (\sigma \cdot_\afrak (z_\tau,a_\tau)) (\Phi(a_{\sigma\tau},a_{\sigma\tau}),-a_{\sigma\tau}), & \text{par \eqref{eqInverseOnF}}, \\
		= \ & (0,a_\sigma) (\tensor[^\sigma]{z}{_\tau} + \Phi(\afrak_\sigma,\tensor[^\sigma]{a}{_\tau}) - \Phi(\tensor[^\sigma]{a}{_\tau},\afrak_\sigma), \tensor[^\sigma]{a}{_\tau}) (\Phi(a_{\sigma\tau},a_{\sigma\tau}),-a_{\sigma\tau}), & \text{par \eqref{eqlemImageOfDelta}}, \\
		= \ & (\tensor[^\sigma]{z}{_\tau} + \Phi(\afrak_\sigma,\tensor[^\sigma]{a}{_\tau}) - \Phi(\tensor[^\sigma]{a}{_\tau},\afrak_\sigma) + \Phi(a_\sigma,\tensor[^\sigma]{a}{_\tau}), a_\sigma + \tensor[^\sigma]{a}{_\tau}) (\Phi(a_{\sigma\tau},a_{\sigma\tau}),-a_{\sigma\tau}), & \text{par \eqref{eqGroupLawOnF}}, \\
		= \ & (\tensor[^\sigma]{z}{_\tau} + \Phi(\afrak_\sigma,\tensor[^\sigma]{a}{_\tau}) - \Phi(\tensor[^\sigma]{a}{_\tau},\afrak_\sigma) + \Phi(a_\sigma,\tensor[^\sigma]{a}{_\tau}), a_{\sigma\tau}) (\Phi(a_{\sigma\tau},a_{\sigma\tau}),-a_{\sigma\tau}) \\
		= \ & (\tensor[^\sigma]{z}{_\tau} +\Phi(\afrak_\sigma,\tensor[^\sigma]{a}{_\tau}) - \Phi(\tensor[^\sigma]{a}{_\tau},\afrak_\sigma) + \Phi(a_\sigma,\tensor[^\sigma]{a}{_\tau}) + \Phi(a_{\sigma\tau},a_{\sigma\tau}) - \Phi(a_{\sigma\tau},a_{\sigma\tau}), 0), & \text{par \eqref{eqGroupLawOnF}}, \\
		= \ & (\tensor[^\sigma]{z}{_\tau} + \Phi(\afrak_\sigma,\tensor[^\sigma]{a}{_\tau}) - \Phi(\tensor[^\sigma]{a}{_\tau},\afrak_\sigma) + \Phi(a_\sigma,\tensor[^\sigma]{a}{_\tau}), 0),
	\end{align*}
	d'o\`u
		\begin{align*}
			 f_\sigma (\sigma \cdot_\afrak f_\tau)f_{\sigma\tau}^{-1} & = (z_\sigma,a_\sigma)(\sigma \cdot_\afrak (z_\tau,a_\tau))(z_{\sigma\tau},a_{\sigma\tau})^{-1} \\
			 & = (z_\sigma,0)(0,a_\sigma)(\sigma \cdot_\afrak (z_\tau,a_\tau))(0,a_{\sigma\tau})^{-1}(z_{\sigma\tau},0)^{-1}  \\
			 & = (z_\sigma,0)(\tensor[^\sigma]{z}{_\tau} + \Phi(\afrak_\sigma,\tensor[^\sigma]{a}{_\tau}) - \Phi(\tensor[^\sigma]{a}{_\tau},\afrak_\sigma) + \Phi(a_\sigma,\tensor[^\sigma]{a}{_\tau}), 0)(-z_{\sigma\tau},0) \\ 
			 & = (z_\sigma + \tensor[^\sigma]{z}{_\tau} - z_{\sigma\tau} + \Phi(\afrak_\sigma,\tensor[^\sigma]{a}{_\tau}) - \Phi(\tensor[^\sigma]{a}{_\tau},\afrak_\sigma) + \Phi(a_\sigma,\tensor[^\sigma]{a}{_\tau}), 0) \\
			 & = (z_\sigma + \tensor[^\sigma]{z}{_\tau} - z_{\sigma\tau} + \phi(\xfrak_\sigma \otimes \tensor[^\sigma]{y}{_\tau} - \tensor[^\sigma]{x}{_\tau}\otimes \yfrak_\sigma + x_\sigma \otimes \tensor[^\sigma]{y}{_\tau}), 0) \\
			 & = ((\diff z)_{\sigma,\tau} + \phi((\xfrak \cup y)_{\sigma,\tau} - \tensor[^\sigma]{x}{_\tau}\otimes \yfrak_\sigma + (x \cup y)_{\sigma,\tau}),0).
		\end{align*}
	Notons de plus que
		\begin{equation*}
                x_{\sigma\tau} \otimes \yfrak_{\sigma\tau} = (x_\sigma + \tensor[^\sigma]{x}{_\tau}) \otimes (\yfrak_\sigma + \tensor[^\sigma]{\yfrak}{_\tau}) = x_\sigma \otimes \yfrak_\sigma + \tensor[^\sigma]{x}{_\tau} \otimes \yfrak_\sigma + x_\sigma \otimes \tensor[^\sigma]{\yfrak}{_\tau} + \tensor[^\sigma]{x}{_\tau} \otimes \tensor[^\sigma]{\yfrak}{_\tau},
            \end{equation*}
	d'o\`u $-\tensor[^\sigma]{x}{_\tau} \otimes \yfrak_\sigma = x_\sigma \otimes \tensor[^\sigma]{\yfrak}{_\tau} + (x_{\sigma} \otimes \yfrak_{\sigma} + \tensor[^\sigma]{(}{}x_{\tau} \otimes \yfrak_{\tau})  - x_{\sigma\tau} \otimes \yfrak_{\sigma\tau}) = (x \cup \yfrak)_{\sigma,\tau} + \diff(x \otimes \yfrak)_{\sigma,\tau}$, donc on a
		\begin{equation} \label{eqlemImageOfDelta1} 
			f_\sigma (\sigma \cdot_\afrak f_\tau)f_{\sigma\tau}^{-1} = ((\diff z)_{\sigma,\tau} + \phi((\xfrak \cup y)_{\sigma,\tau} + (x \cup \yfrak)_{\sigma,\tau} + \diff(x \otimes \yfrak)_{\sigma,\tau} + (x \cup y)_{\sigma,\tau}),0)
		\end{equation} 
    De \eqref{eqlemImageOfDelta1}, on voit que $f$ est un cocycle si et seulement si $a$ l'est et
	\begin{equation*}
            \diff z + \phi_\ast(\xfrak \cup y + x \cup \yfrak + x \cup y + \diff(x \otimes \yfrak)) = 0,
        \end{equation*}
    d'o\`u le premier point du lemme.

    \item On consid\`ere la $1$-cocha\^ine $f = (0,a): \Gamma_K \to F$ relevant $a$. Alors 
        \begin{equation*}
            f_\sigma \tensor[^\sigma]{f}{_\tau} f_{\sigma\tau}^{-1} = (\phi((x \cup y)_{\sigma,\tau}),0)
        \end{equation*} 
    pour tous $\sigma,\tau \in \Gamma_K$, au vu de \eqref{eqlemImageOfDelta1}. Par d\'efinition de $\Delta$, on a bien $\Delta([a]) = \phi_\ast[x \cup y]$.
\end{enumerate}
\end{proof}

Le lemme suivant sera utilis\'e dans la preuve du th\'eor\`eme \ref{thmWeakApproximation} (c'est-\`a-dire du th\'eor\`eme \ref{thmB}).

\begin{lem} \label{lemCohomologousOnF}
    Gardons les notations au d\'ebut du paragraphe, et soit $\afrak = (\xfrak, \yfrak): \Gamma_K \to M \oplus M$ un $1$-cocycle. Soient $f = (z,a)$ et $f' = (z',a')$ des $1$-cocycles $\Gamma_K \to \tensor[_\afrak]{F}{}$ ({\em cf.} proposition \ref{propSpringerClass}), o\`u $a = (x,y)$ et $a' = (x',y')$ sont des $1$-cocycles $\Gamma_K \to M \oplus M$ (voir lemme \ref{lemImageOfDelta}). On suppose que $\alpha = (\xi,\eta) \in M \oplus M$ est un \'el\'ement satisfaisant $a' = a + \diff \alpha$ et l'on consid\`ere la $1$-cocha\^ine
	\begin{equation*}
            c:=z' - z + \phi_\ast(-(x + \xfrak) \cup \eta + \xi \cup (y' + \yfrak) + \diff \xi \otimes \yfrak): \Gamma_K \to Z.
        \end{equation*}
	\begin{enumerate}
		\item $c$ est un cocycle.
		
		\item Si $c$ est un cobord, $f$ et $f'$ sont cohomologues.
		
		\item Supposons $\phi$ surjectif. \'Ecrivons $z = \phi_\ast \varepsilon$ et $z' = \phi_\ast \varepsilon'$, o\`u $\varepsilon,\varepsilon': \Gamma_K \to M \otimes M$ sont des $1$-cocha\^ines. Alors on peut \'ecrire
		\begin{equation*}
                \diff \varepsilon + \xfrak \cup y + x \cup \yfrak + x \cup y + \diff(x \otimes \yfrak) = j_\ast \lambda \quad \text{et} \quad \diff \varepsilon' + \xfrak \cup y' + x' \cup \yfrak + x' \cup y' + \diff(x' \otimes \yfrak) = j_\ast \lambda',
            \end{equation*}
		o\`u $j: \Ker \phi \hookrightarrow M \otimes M$ est l'inclusion et o\`u  $\lambda,\lambda': \Gamma_K \times \Gamma_K \to \Ker \phi$ sont des $2$-cocycles. De plus, si $\delta: \H^1(K,Z) \to \H^2(K, \Ker \phi)$ d\'esigne le morphisme connectant induit par la suite exacte  $0 \to \Ker \phi \xrightarrow{j} M \otimes M \xrightarrow{\phi} Z \to 0$, alors $\delta([c]) = [\lambda' - \lambda]$.
	\end{enumerate}
\end{lem}
\begin{proof}
\begin{enumerate}
	\item Par l'hypoth\`ese, $x' = x + \diff \xi$ et $y' = y + \diff \eta$. Par le lemme \ref{lemImageOfDelta}, on a
		\begin{equation*}
                \diff z = -\phi_\ast(\xfrak \cup y + x \cup \yfrak + x \cup y + \diff(x \otimes \yfrak))
            \end{equation*}
	et
	\begin{align*}
		\diff z' & = -\phi_\ast(\xfrak \cup y' + x' \cup \yfrak + x' \cup y' + \diff(x' \otimes \yfrak)) \\
		& = -\phi_\ast(\xfrak \cup (y + \diff \eta) + (x + \diff \xi) \cup \yfrak + (x + \diff \xi) \cup  (y + \diff \eta) + \diff((x + \diff \xi) \otimes \yfrak))
	\end{align*}
     puisque $f$ et $f'$ sont des cocycles. D'o\`u
	\begin{align*}
		\diff z' - \diff z & = -\phi_\ast(\xfrak \cup \diff \eta + \diff \xi \cup \yfrak + x \cup \diff \eta + \diff \xi \cup y + \diff \xi \cup \diff \eta + \diff(\diff \xi \otimes \yfrak)) \\
		& = -\phi_\ast((x + \xfrak) \cup \diff \eta + \diff \xi \cup (y + \yfrak + \diff \eta) + \diff (\diff \xi \otimes \yfrak)) \\
		& = -\phi_\ast((x + \xfrak) \cup \diff \eta + \diff \xi \cup (y' + \yfrak) + \diff (\diff \xi \otimes \yfrak)).
	\end{align*}
    Notant que $x,\xfrak, y'$ et $\yfrak$ sont des cocycles, on en d\'eduit que
	\begin{align*}
		\diff c & =  \diff (z' - z + \phi_\ast(-(x + \xfrak) \cup \eta + \xi \cup (y' + \yfrak) + \diff \xi \otimes \yfrak)) \\
		& = \diff z' - \diff z + \phi_\ast((x + \xfrak) \cup \diff \eta + \diff \xi \cup (y' + \yfrak) + \diff(\diff \xi \otimes \eta)) \\
		& = 0,
	\end{align*}
	donc $c$ est bien un cocycle.
	
	\item \'Ecrivons $c = \diff \zeta$ avec $\zeta \in 
	Z$. Pour tout $\sigma \in \Gamma_K$, on a
		\begin{equation} \label{eqlemCohomologousOnF1}
			a'_\sigma + \alpha = a_\sigma + \tensor[^\sigma]{\alpha}{}
		\end{equation}
	puisque $a' = a + \diff\alpha$, et
		\begin{equation} \label{eqlemCohomologousOnF2}
			c_\sigma = z'_\sigma - z_\sigma + \phi(-(x_\sigma+\xfrak_\sigma) \otimes \tensor[^\sigma]{\eta}{} + \xi \otimes (y'_\sigma + \yfrak_\sigma) + (\tensor[^\sigma]{\xi}{} - \xi) \otimes \yfrak_\sigma)
		\end{equation}
	par d\'efinition de $c$. Calculons
	\begin{align*}
		& (\zeta,\alpha)^{-1} f_\sigma (\sigma \cdot_\afrak (\zeta,\alpha)) \\
		= \ & (\Phi(\alpha,\alpha) -\zeta,-\alpha) (z_\sigma,a_\sigma) (\sigma \cdot_\afrak (\zeta,\alpha)), & \text{par \eqref{eqInverseOnF}}, \\
		= \ & (\Phi(\alpha,\alpha) -\zeta,-\alpha) (z_\sigma,a_\sigma)  (\tensor[^\sigma]{\zeta}{} + \Phi(\afrak_\sigma,\tensor[^\sigma]{\alpha}{}) - \Phi(\tensor[^\sigma]{\alpha }{},\afrak_\sigma), \tensor[^\sigma]{\alpha}{}), & \text{par \eqref{eqlemImageOfDelta}}, \\
		= \ & (\Phi(\alpha,\alpha) -\zeta,-\alpha)  (z_\sigma + \tensor[^\sigma]{\zeta}{} + \Phi(\afrak_\sigma,\tensor[^\sigma]{\alpha}{}) - \Phi(\tensor[^\sigma]{\alpha }{},\afrak_\sigma)+ \Phi(a_\sigma,\tensor[^\sigma]{\alpha}{}), a_\sigma + \tensor[^\sigma]{\alpha}{}), & \text{par \eqref{eqGroupLawOnF}}, \\
		= \ & (\Phi(\alpha,\alpha)-\zeta,-\alpha) (z_\sigma + \tensor[^\sigma]{\zeta}{} + \Phi(a_\sigma+\afrak_\sigma,\tensor[^\sigma]{\alpha}{}) - \Phi(\tensor[^\sigma]{\alpha }{},\afrak_\sigma), a'_\sigma + \alpha), & \text{par \eqref{eqlemCohomologousOnF1}}, \\
		= \ & (\Phi(\alpha,\alpha) -\zeta + z_\sigma + \tensor[^\sigma]{\zeta}{} + \Phi(a_\sigma+\afrak_\sigma,\tensor[^\sigma]{\alpha}{}) - \Phi(\tensor[^\sigma]{\alpha }{},\afrak_\sigma) - \Phi(\alpha,a'_\sigma + \alpha), a'_\sigma), & \text{par \eqref{eqGroupLawOnF}}, \\
		= \ & (\tensor[^\sigma]{\zeta}{} - \zeta + z_\sigma + \Phi(a_\sigma+\afrak_\sigma,\tensor[^\sigma]{\alpha}{}) - \Phi(\tensor[^\sigma]{\alpha }{},\afrak_\sigma) - \Phi(\alpha,a'_\sigma), a'_\sigma) \\
		= \ & (c_\sigma + z_\sigma + \phi((x_\sigma+\xfrak_\sigma) \otimes \tensor[^\sigma]{\eta}{} - \tensor[^\sigma]{\xi }{} \otimes  \yfrak_\sigma - \xi \otimes  y'_\sigma), a'_\sigma), & \text{car } c = \diff \zeta,\\
		= \ & (c_\sigma + z_\sigma + \phi((x_\sigma+\xfrak_\sigma) \otimes \tensor[^\sigma]{\eta}{} - \xi \otimes (y'_\sigma + \yfrak_\sigma) - (\tensor[^\sigma]{\xi}{} - \xi) \otimes \yfrak_\sigma), a'_\sigma) \\
		= \ & (z'_\sigma,a'_\sigma), & \text{par \eqref{eqlemCohomologousOnF2}}, \\
		= \ & f'_\sigma.
	\end{align*}
	Cette \'egalit\'e signifie que $f$ et $f'$ sont cohomologues.
	
	\item Comme $f$ est un cocycle, on a
		\begin{equation*}
                \phi_\ast(\diff \varepsilon + \xfrak \cup y + x \cup \yfrak + x \cup y + \diff(x \otimes \yfrak)) = \diff z + \phi_\ast(\xfrak \cup y + x \cup \yfrak + x \cup y + \diff(x \otimes \yfrak)) = 0
            \end{equation*}
	par le lemme \ref{lemImageOfDelta}. Donc on peut \'ecrire
		\begin{equation} \label{eqlemCohomologousOnF3}
			\diff \varepsilon + \xfrak \cup y + x \cup \yfrak + x \cup y + \diff(x \otimes \yfrak) = j_\ast \lambda
		\end{equation}
	pour une $2$-cocha\^ine $\lambda: \Gamma_K \times \Gamma_K \to \Ker \phi$. Le membre du gauche de \eqref{eqlemCohomologousOnF3} \'etant un cocycle, il en va de m\^eme pour $\lambda$ puisque $j$ est injectif. De m\^eme, on a 
		\begin{equation} \label{eqlemCohomologousOnF4}
			\diff \varepsilon' + \xfrak \cup y' + x' \cup \yfrak + x' \cup y' + \diff(x' \otimes \yfrak) = j_\ast \lambda'
		\end{equation}
	pour un $2$-cocycle $\lambda': \Gamma_K \times \Gamma_K \to \Ker \phi$. Consid\'erons maintenant la $1$-cocha\^ine
		\begin{equation*}
                e:=\varepsilon' - \varepsilon - (x + \xfrak) \cup \eta + \xi \cup (y' + \yfrak) + \diff \xi \otimes \yfrak: \Gamma_K \to M \otimes M.
            \end{equation*}
	Alors $\phi_\ast e = c$. Calculons le cobord
		\begin{align*}
			\diff e & = \diff \varepsilon' - \diff \varepsilon + (x + \xfrak) \cup \diff \eta + \diff \xi \cup (y' + \yfrak) + \diff(\diff \xi \otimes \yfrak) \\
			& = \diff \varepsilon' - \diff \varepsilon + (x + \xfrak) \cup (y' - y) + (x' - x) \cup (y' + \yfrak) + \diff((x' - x) \otimes \yfrak) \\
			& = (\diff \varepsilon' + \xfrak \cup y' + x' \cup \yfrak + x' \cup y' + \diff(x' \otimes \yfrak)) - (\diff \varepsilon + \xfrak \cup y + x \cup \yfrak + x \cup y + \diff(x \otimes \yfrak))\\
			& = j_\ast (\lambda' - \lambda),
		\end{align*}
	o\`u la derni\`ere \'egalit\'e vient de \eqref{eqlemCohomologousOnF3} et \eqref{eqlemCohomologousOnF4}. Par d\'efinition de $\delta$, on a $\delta([c]) = [\lambda' - \lambda]$.
\end{enumerate}
\end{proof}  
\section{Calculs des groupes de Brauer} \label{section4}

\subsection{G\'en\'eralit\'es} \label{subsection41}

Soit $K$ un corps de caract\'eristique nulle et soit $X$ un $K$-espace homog\`ene de $\SL_n$. Rappelons d'abord quelques groupes associ\'es \`a $X$.
\begin{itemize}
	\item Soit $\ol{X} = X \times_K \ol{K}$, alors $\Br \ol{X}$ est le {\em groupe de Brauer g\'eom\'etrique} de $X$.
	\item $\Br_1 X = \Ker(\Br X \to \Br \ol{X})$ est le {\em groupe de Brauer alg\'ebrique} de $X$ et 
	$$\Br_{\nr,1} X = \Br_1 X \cap \Br_{\nr} X$$
	est sa {\em partie non ramifi\'ee}.
	\item $\Br_0 X = \Img(\Br K \to \Br X) \subseteq \Br_{\nr, 1} X$ est le {\em sous-groupe des \'el\'ements constants} de $\Br X$.
	\item $\Br_a X = (\Br_1 X)/(\Br_0 X)$ est le {\em groupe de Brauer arithm\'etique} de $X$ et 		
	$$\Br_{\nr,a} X = (\Br_{\nr,1})/(\Br_0 X)$$
	est sa partie non ramifi\'ee.
	\item Lorsque $K = k$ est un corps de nombres, on note $X_v = X \times_k k_v$ pour toute place $v$ de $k$. On d\'efinit les groupes
	$$\Be(X) = \Ker\left(\Br_a X\to \prod_{v \in \Omega_k} \Br_a X_v \right) = \{\alpha \in \Br_a X: \forall v \in \Omega_k,\alpha_v = 0 \}$$
	et 
	$$\Be_\omega(X) = \{\alpha \in \Br_a X: \alpha_v = 0 \text{ pour presque toute } v \in \Omega_k\}.$$
\end{itemize}
Notons $F$ le stabilisateur g\'eom\'etrique de $X$, qu'on suppose fini. Alors $F$ est munie d'une action ext\'erieure de $\Gamma_K$ ({\em cf.} paragraphe \ref{subsection21}), ce qui induit une action de $\Gamma_K$ sur $F^{\ab} = F/[F,F]$. Les groupes $\Br_a X$, $\Be(X)$ et $\Be_\omega(X)$ s'expriment en termes de $F^{\ab}$; c'est un argument classique qui se trouve par exemple dans \cite[Theorem 4.1.1]{skorobogatov2001torsors}.

Rappelons que $\hat{-}$ d\'esigne le foncteur de dual de Cartier.

\begin{prop} \label{propBrauerGroups}
    Soit $K$ un corps de caract\'eristique nulle satisfaisant $\H^3(K,\ol{K}^\times) = 0$ et soit $X$ un $K$-espace homog\`ene de $\SL_n$ \`a stabilisateur g\'eom\'etrique $F$. On munit $\hat{F} = \Hom(F,\ol{K}^\times)$ de l'action de $\Gamma_K$ induite par son action ext\'erieure sur $F$.
	\begin{enumerate}
		\item On a $\Br_a X = \H^1(K,\hat{F})$.
		\item Si $K = k$ est un corps de nombres, $\Be(X) = \Sha^1(k,\hat{F})$ et $\Be_\omega(X) = \Sha^1_\omega(k,\hat{F})$.
	\end{enumerate}
\end{prop}

Appliquons maintenant la proposition \ref{propBrauerGroups} aux espaces homog\`enes de Borovoi--Kunyavski\u{\i}.

\begin{cor} \label{corBe}
    Soient $k$ un corps de nombres, $M$ et $Z$ des $\Gamma_k$-modules finis, et $\phi: M \otimes M \to Z$ un morphisme $\Gamma_k$-\'equivariant. \`A partir de $\phi$, on construit l'extension
	\begin{equation*}
	    0 \to Z \to F \to M \oplus M \to 0
	\end{equation*}
    de groupes abstraits comme dans le paragraphe \ref{subsection22}. Supposons $Z = Z(F) = [F,F]$. Si $X$ est un espace homog\`ene de $\SL_n$ \`a stabilisateur g\'eom\'etrique $F$ (muni de n'importe quelle action ext\'erieure de $\Gamma_k$ compatible avec celles sur $Z$ et sur $M \oplus M$, pas n\'ecessairement celle induite par l'action coordonn\'ees par coordonn\'ees), alors $\Be(X) = \Sha^1(k,\hat{M})^2$ et $\Be_\omega(X) = \Sha^1_\omega(k,\hat{M})^2$.
\end{cor}
\begin{proof}
    Puisque $Z = [F,F]$, $F^{\ab} = F/Z = M \oplus M$ et donc $\hat{F} = \hat{F}^{\ab} = \hat{M} \oplus \hat{M}$. Ainsi, la proposition \ref{propBrauerGroups} donne $\Be(X) = \Sha^1(k,\hat{F}) = \Sha^1(k,\hat{M})^2$ et $\Be_\omega(X) = \Sha_\omega^1(k,\hat{F}) = \Sha^1_\omega(k,\hat{M})^2$.
\end{proof}

\subsection{\'Etude de la fl\`eche $\H^1(K,F) \to \H^1(K,F^{\ab})$} \label{subsection42} 

Soit $X$ un espace homog\`ene de Borovoi--Kunyavski\u{\i} sur un corps de nombres $k$. Afin de calculer le sous groupe $\Br_{\nr,a} X$ de $\Br_a X$, on va utiliser la formule de Demarche. Il faudra \'etudier l'image de la fl\`eche $\H^1(k_v,F) \to \H^1(k_v,F^{\ab})$ (cette image n'est pas forc\'ement un sous-groupe de $\H^1(k_v,F^{\ab})$) pour presque toute place $v$ de $k$. 

On se donne alors un corps $p$-adique $K$ et un $K$-groupe fini $F$. Dans le cas o\`u $F$ est un $K$-groupe constant, Demarche \cite[\S 3]{Demarche} a donn\'e une description assez explicite du sous-groupe de $\H^1(K,F^{\ab})$ engendr\'e par l'image de la fl\`eche $\H^1(K,F) \to \H^1(K,F^{\ab})$. Nous allons adapter son argument pour le cas o\`u $F$ est un $K$-groupe {\em non ramifi\'e}, c'est-\`a-dire que l'action du sous-groupe d'inertie de $\Gamma_K$ sur $F$ est trivial. 

On note $K_{\nr}$ (resp. $K_{\mr}$) l'extension maximale non ramifi\'ee (resp. mod\'er\'ement ramifi\'ee) de $K$. Soit $q$ le cardinal du corps r\'esiduel de $K$. Notons $\pi: F \to F^{\ab}$ la projection. Posons $\Gamma = \Gal(K_{\mr}/K) = \ol{\pair{\sigma,\tau | \sigma \tau \sigma^{-1} = \tau^q}}$. Par l'hypoth\`ese que $F$ est non ramifi\'e, $\tau$ agit trivialement sur $F$. \'Etudions la fl\`eche $\pi_\ast: \H^1(\Gamma,F) \to \H^1(\Gamma,F^{\ab})$.

Soit $f: \Gamma \to F$ un cocycle. Comme $\tau$ agit trivialement sur $F$, on a
    \begin{equation*}
        f(\tau)^q = f(\tau^q) = f(\sigma \tau\sigma^{-1}) = f(\sigma) \tensor[^\sigma]{f(\tau)}{} \tensor[^{\sigma\tau}]{f(\sigma^{-1})}{} = f(\sigma) \tensor[^\sigma]{f(\tau)}{} \tensor[^{\sigma}]{f(\sigma^{-1})}{}.
    \end{equation*}
Comme $1 = f(1) = f(\sigma\sigma^{-1}) = f(\sigma)\tensor[^{\sigma}]{f(\sigma^{-1})}{}$, on obtient $\tensor[^{\sigma}]{f(\sigma^{-1})}{} = f(\sigma)^{-1}$ et donc
    \begin{equation*}
        f(\tau)^q = f(\sigma) \tensor[^\sigma]{f(\tau)}{}f(\sigma)^{-1}.
    \end{equation*}
Inversement, soient $a,b \in F$ tels que $a\tensor[^\sigma]{b}{}a^{-1} = b^q$. Alors il existe un unique cocycle $f: \Gamma \to F$ tel que $f(\sigma) = a$ et $f(\tau) = b$.

\begin{defn}
    On appelle {\em $q$-relevable} tout \'el\'ement $\ol{b} \in F^{\ab}$ ayant un relev\'e $b \in F$ tel que $b^q$ soit conjugu\'e \`a $\tensor[^\sigma]{b}{}$ (en particulier, $\tensor[^\sigma]{\ol{b}}{} = \ol{b}^q$).
\end{defn}

L'inverse d'un \'el\'ement $q$-relevable est encore $q$-relevable. Si $f,g: \Gamma \to F^{\ab}$ sont deux cocycles cohomologues, alors $f(\tau) = g(\tau)$ (car l'action de $\tau$ sur $F$ est triviale), donc $f(\tau)$ est $q$-relevable si et seulement si $g(\tau)$ l'est aussi.

On note $I(F)$ l'image de la fl\`eche $\pi_\ast: \H^1(\Gamma,F) \to \H^1(\Gamma,F^{\ab})$, et $J(F)$ le sous-ensemble de $\H^1(\Gamma,F^{\ab})$ form\'e des classes des cocycles $f$ tels que $f(\tau)$ soit $q$-relevable. C'est \'evident que $I(F) \subseteq J(F)$, et que $J(F)$ est stable par inversion.

\begin{lem} \label{lemIF}
    Les sous-groupes de $\H^1(\Gamma,F^{\ab})$ engendr\'es par $I(F)$ et par $J(F)$ co\"incident.
\end{lem}
\begin{proof}
    Soit $f: \Gamma \to F^{\ab}$ un cocycle tel que $[f] \in J(F)$. Il existe un relev\'e $b \in F$ de $f(\tau)$ et un \'el\'ement $a \in F$ tel que $a \tensor[^\sigma]{b}{} a^{-1} = b^q$. Soit $\tilde{f}_1: \Gamma \to F^{\ab}$ le cocycle d\'etermin\'e par $\tilde{f}_1(\sigma) = a$ et $\tilde{f}_1(\tau) = b$. Posons $f_1 = \pi \circ \tilde{f}_1$, alors $f_1(\sigma) = \pi(a)$ et $f_1(\tau) = f(\tau)$. 
	
    Soit $a' \in F$ un relev\'e de $\pi(a)^{-1}f(\sigma)$. Soit $\tilde{f}_2: \Gamma \to F$ le cocycle d\'etermin\'e par $\tilde{f}_2(\sigma) = a'$ et $\tilde{f}_2(\tau) = 1$. Posons $f_2 = \pi \circ \tilde{f}_2$, alors $f_1(\sigma)f_2(\sigma) = f(\sigma)$ et $f_1(\tau)f_2(\tau) = f(\tau)$, donc $[f] = [f_1] + [f_2] \in \H^1(\Gamma,F^{\ab})$. Comme $[f_1]$, $[f_2] \in I(F)$, $[f]$ appartient au sous-groupe engendr\'e par $I(F)$.
\end{proof}

\begin{prop} \label{propIF2}
    Soit $K$ un corps $p$-adique et notons $q$ le cardinal du corps r\'esiduel de $K$. Soit $\sigma \in \Gal(K_{\mr}/K)$ un relev\'e de l'automorphisme de Frobenius. Soit $F$ un $K$-groupe fini non ramifi\'e de cardinal $|F| < p$. Alors on a \'equivalence entre
	\begin{enumerate}
            \item $\H^1(K,F^{\ab})$ est engendr\'e par l'image de la fl\`eche $\H^1(K,F) \to \H^1(K,F^{\ab})$;
            \item le sous-groupe $\{a \in F^{\ab} : \tensor[^\sigma]{a}{} = a^q\}$ de $F^{\ab}$ est engendr\'e par les \'elements $q$-relevables.
	\end{enumerate}
\end{prop}
\begin{proof}
    Pour tout cocycle $f: \Gamma_K \to F$, l'ensemble $\{\upsilon \in \Gamma_K: f(\upsilon) = 1\}$ est un sous-groupe ouvert de $\Gamma_K$. Il correspond \`a une extension finie $L/K$. On v\'erifie sans peine que pour tous $\upsilon,\upsilon' \in \Gamma_K$, $f(\upsilon') = f(\upsilon)$ \'equivaut \`a $\upsilon' \in \upsilon \Gamma_L$. En particulier, $[L: K] = [\Gamma_K : \Gamma_L] \le |F| < p$, donc $L/K$ est mod\'er\'ement ramifi\'ee. Ainsi, tout cocycle $\Gamma_K \to F$ se factorise par un cocycle $\Gal(K_{\mr}/K) \to F$ et il en va de m\^eme pour $F^{\ab}$. Donc le premier point \'equivaut \`a dire que $\H^1(\Gal(K_{\mr}/K),F^{\ab})$ est engendr\'e par l'image $I(F)$ de la fl\`eche $\H^1(\Gal(K_{\mr}/K),F) \to \H^1(\Gal(K_{\mr}/K),F^{\ab})$. 
	
    Soit $\tau$ un g\'en\'er\'ateur topologique de $\Gal(K_{\mr}/K_{\nr})$, de sorte que $\sigma \tau \sigma^{-1} = \tau^q$.
	
    Supposons 1. Soit $a \in F^{\ab}$ satisfaisant $\tensor[^\sigma]{a}{} = a^q$. Alors on peut d\'efinir un cocycle 
        \begin{equation*}
            f: \Gal(K_{\mr}/K) \to F^{\ab}
        \end{equation*}
    par $f(\sigma) = 1$ et $f(\tau) = a$. Comme $I(F) \subseteq J(F)$, on peut \'ecrire $[f] = [f_1] + \cdots + [f_r]$, o\`u $[f_i] \in J(F)$, c'est-\`a-dire que chacun des $f_i(\tau)$ est $q$-relevable. Or $\tau$ agit trivialement sur $F^{\ab}$, donc $a = f(\tau) = f_1(\tau)\cdots f_r(\tau)$, d'o\`u 2.
	
    Supposons 2. Soit $f: \Gal(K_{\mr}/K) \to F^{\ab}$ un cocycle. Alors $\tensor[^\sigma]{f(\tau)}{} = f(\tau)^q$ et donc $f(\tau) = \ol{b}_1  \cdots \ol{b}_r$, o\`u chaque $\ol{b}_i \in F^{\ab}$ est $q$-relevable. D\'efinissons les cocycles $f_i: \Gal(K_{\mr}/K) \to F^{\ab}$ par
        \begin{equation*}
            f_1(\sigma) = f(\sigma), \quad f_1(\tau) = \ol{b}_1,
        \end{equation*}
        \begin{equation*}
            f_i(\sigma) = 1, \quad f_i(\tau) = \ol{b}_i, \quad i=2,\ldots,r.
        \end{equation*}
    Alors $[f_i] \in J(F)$ et $[f] = [f_1] +  \cdots + [f_r]$, donc $[f]$ appartient au sous-groupe engendr\'e par $J(F)$. Par le lemme \ref{lemIF}, $[f]$ appartient au sous-groupe engendr\'e par l'image de $\H^1(\Gal(K_{\mr}/K),F) \to \H^1(\Gal(K_{\mr}/K),F^{\ab})$, d'o\`u 1.
\end{proof}

\subsection{Groupe de Brauer arithm\'etique} \label{subsection43}

Nous allons maintenant combiner l'analyse du paragraphe \ref{subsection32} avec la formule de Demarche pour calculer le groupe de Brauer arithm\'etique non ramifi\'e d'un espace homog\`ene de Borovoi--Kunyavski\u{\i}.

Rappelons d'abord la dualit\'e locale de Tate (voir par exemple \cite[Chapitre 10]{harari2017cohomologie}). Soit $K$ un corps $p$-adique $K$ et soit $A$ un $\Gamma_K$-module fini, alors l'accouplement
    \begin{equation*}
        \H^1(K,A) \times \H^1(K,\hat{A}) \to \Qbb/\Zbb, \quad (a,\alpha) \mapsto \inv_K(a \cup \alpha)
    \end{equation*}
est une dualit\'e parfaite de groupes ab\'eliens finis.

\begin{prop} \label{propBrNr}
    Soient $k$ un corps de nombres, $M$ et $Z$ des $\Gamma_k$-modules finis, et $\phi: M \otimes M \to Z$ un morphisme $\Gamma_k$-\'equivariant. \`A partir de $\phi$, on construit l'extension
	\begin{equation*}
            0 \to Z \to F \to M \oplus M \to 0
	\end{equation*}
    de groupes abstraits comme dans le paragraphe \ref{subsection22}. Supposons $Z = Z(F) = [F,F]$.
	\begin{enumerate}
            \item On munit $F$ de l'action coordonn\'ees par coordonn\'ees de $\Gamma_k$. Alors le groupe $\H^1(k_v,F^{\ab})$ est engendr\'e par l'image de la fl\`eche $\H^1(k_v, F) \to \H^1(k_v,F^{\ab})$ pour presque toute $v \in \Omega_k$.
            
            \item Soit $X$ est un espace homog\`ene de $\SL_n$ \`a stabilisateur g\'eom\'etrique $F$, muni de l'action ext\'erieure de $\Gamma_k$ induite par l'action coordonn\'ees par coordonn\'ees. Alors $\Br_{\nr,a} X = \Be_{\omega} (X) = \Sha^1_\omega(k,\hat{M})^2$.
	\end{enumerate}
\end{prop}
\begin{proof}
    \begin{enumerate}
        \item Soit $L/k$ une extension finie galoisienne d\'eployant $F$. Soit $v$ une place de $k$ qui est non ramifi\'ee dans $L/k$, et qui divise un nombre premier impair $p > |F|$. On va d\'emontrer que $\H^1(k_v,F^{\ab})$ est engendr\'e par l'image de la fl\`eche $\H^1(k_v, F) \to \H^1(k_v,F^{\ab})$. Au vu de la proposition \ref{propIF2}, il suffit de d\'emontrer que tout \'el\'ement $a \in F^{\ab}$ satisfaisant $\tensor[^\sigma]{a}{} = qa$ est $q$-relevable, o\`u $q$ est le cardinal du corps r\'esiduel de $k_v$ et o\`u $\sigma \in \Gal(k_{v,\mr}/k)$ est un relev\'e de l'automorphisme de Frobenius. Notons que $F^{\ab} = F/[F,F] = F/Z = M \oplus M$.
	
        Rappelons de \eqref{eqBigPhi} qu'on dispose d'une application biadditive 
            \begin{equation*}
                \Phi: (M \oplus M) \times (M \oplus M) \to Z, \quad ((x,y),(x',y')) \mapsto \phi(x \otimes y').
            \end{equation*}
        Soit $a = (x,y) \in M \oplus M$ tel que $\tensor[^\sigma]{a}{} = qa$ et montrons que $a$ est $q$-relevable. Consid\'erons alors le relev\'e $(0,a) \in F$ de $a$. On va d\'emontrer que $(0,a)^q$ est conjug\'e \`a $\tensor[^\sigma]{(0,a)}{} = (0,qa)$.
        
        \`A l'aide de \eqref{eqGroupLawOnF}, on peut calculer
            \begin{equation} \label{eqpropBrNr1}
                (0,a)^q = \left(\tfrac{q(q-1)}{2}\Phi(a,a),qa\right)
		\end{equation}
        par r\'ecurrence. Soit $a' = \left(\tfrac{q+1}{2}x, y\right)$. Alors
            \begin{equation} \label{eqpropBrNr2}
                \Phi(qa,a') = \phi(qx \otimes y) = q\Phi(a,a) \quad \text{et} \quad  \Phi(a',qa) = \phi\left(\tfrac{q+1}{2}x \otimes qy\right) = \tfrac{q(q+1)}{2}\Phi(a,a).
		\end{equation}
	On a
            \begin{align*}
                (0,a')(0,qa)(0,a')^{-1} & = (\Phi(a',qa) - \Phi(qa,a'),qa), & \text{par le lemme \ref{lemCommutator}},\\
                & = (\tfrac{q(q+1)}{2}\Phi(a,a) - q\Phi(a,a),qa), & \text{par \eqref{eqpropBrNr2}},\\
                & = (\tfrac{q(q-1)}{2}\Phi(a,a),qa) \\
                & = (0,a)^q, & \text{par \eqref{eqpropBrNr1}},
		\end{align*}
        donc $(0,a)^q$ est bien conjugu\'e \`a $(0,qa) = \tensor[^\sigma]{(0,a)}{}$, d'o\`u $a$ est $q$-relevable. Cela nous permet de conclure que $\H^1(k_v,F^{\ab})$ est engendr\'e par l'image de $\H^1(k_v,F) \to \H^1(k_v,F^{\ab})$ pour presque toute place $v$ de $k$.
	
    \item Rappelons que $\Br_a X = \H^1(k,\hat{F}) = \H^1(k,\hat{M})^2$ et que $\Br_a X_v = \H^1(k_v,\hat{F}) = \H^1(k_v,\hat{M})^2$ pour toute place $v$ de $k$, {\em cf.} proposition \ref{propBrauerGroups}. Soit $\eta_0 \in \H^2(k,F)$ la classe neutre privil\'egi\'ee (ou $F$ est muni de l'action coordonn\'ees par coordonn\'es de $\Gamma_k$). Alors il existe une unique classe $\beta \in \H^2(k,Z)$ telle que la classe de Springer de $X$ vaut $\eta_X = \beta \cdot \eta_0$ ({\em cf.} paragraphe \ref{subsection31}). Pour presque toute $v \in \Omega_k$, on a $\beta_v = 0 \in \H^2(k_v,Z)$, donc $\loc_v(\eta_X) = \loc_v(\eta_0)$, ou $X_v$ est $k_v$-isomorphe \`a $F \backslash \SL_n$ (voir proposition \ref{propSpringerClass}). Par la formule de Demarche \cite[Th\'eor\`eme 2.1]{Demarche}, un \'el\'ement $\alpha \in \Br_a X$ appartient \`a $\Br_{\nr,a} X$ si et seulement si pour presque toute $v \in \Omega_k$, l'image de $\alpha_v$ dans $\H^1(k_v,\hat{F})$ est orthogonal \`a l'image de $\H^1(k_v,F) \to \H^1(k_v,F^{\ab})$. Par la dualit\'e locale de Tate, ce dernier siginifie que $\alpha_v = 0$ pour presque toute $v \in \Omega_k$, {\em i.e.} que $\alpha \in \Be_\omega(X)$. Ainsi, $\Br_{\nr,a} X = \Be_\omega(X) = \Sha^1_\omega (k,\hat{F}) =  \Sha^1_\omega(k,\hat{M})^2$.
	\end{enumerate}
\end{proof}

\subsection{Groupe de Brauer g\'eom\'etrique} \label{subsection44}

Dans ce paragraphe, soit $K$ un corps {\em alg\'ebriquement clos} de caract\'eristique nulle. Pour calculer le groupe de Brauer g\'eom\'etrique non ramifi\'e des espaces homog\`enes de Borovoi--Kunyavski\u{\i}, on rappelle la formule suivante de Bogomolov \cite[\S 3]{Bogomolov}.

\begin{prop} \label{propBogomolovFormula}
    Soit $F$ un groupe fini, vu comme $K$-groupe. On choisit un plongement $F \hookrightarrow \SL_n$ de $K$-groupes et l'on pose $X = F \backslash \SL_n$. Alors 
        \begin{equation*}
            \Br_{\nr} X = B_0(F):=\Ker\left(\H^2(F,\Qbb/\Zbb) \to \prod_{A} \H^2(A,\Qbb/\Zbb) \right),
        \end{equation*}
    o\`u $A$ parcourt les sous-groupes ab\'eliens de $F$.
\end{prop}

Lorsque $F$ est nilpotent de classe 2, on dispose d'un morphisme
    \begin{equation} \label{eqLambda}
        \lambda_F: \bigwedge^2 F^{\ab} \to Z, \qquad \lambda(a \wedge b) = [\tilde{a},\tilde{b}],
    \end{equation}
o\`u $Z = Z(F)$ et o\`u $\tilde{a},\tilde{b} \in F$ sont des relev\'es respectifs de $a,b \in F/Z$. Si de plus $Z = [F,F]$ (par exemple, pour les espaces homog\`enes de Borovoi--Kunyavski\u{\i}), {\em i.e.} $\lambda_F$ est surjectif, le groupe $B_0(F)$ a la description explicite suivante \cite[Lemma 5.1]{Bogomolov} (en fait, l'\'enonc\'e dans {\em loc. cit.} demande que $F$ soit un $p$-groupe, mais sa d\'emonstration vaut pour $F$ de cardinal quelconque). Un \'enonc\'e similaire se trouve dans \cite[\S 5]{Moravec}. On rappelera la preuve pour la commodit\'e du lecteur. 

\begin{lem} \label{lemBogomolovFormula}
    Soit $F$ un groupe fini v\'erifiant $Z(F) = [F,F]$. Notons $Z = Z(F)$ et $\lambda = \lambda_F$ le morphisme surjectif d\'efini par \eqref{eqLambda}. Alors $B_0(F) = \Hom(S/S_\lambda,\Qbb/\Zbb)$, o\`u $S = \Ker \lambda$ et o\`u $S_\lambda$ est le sous-groupe de $S$ engendr\'e par $S \cap \{a \wedge b: a,b \in F^{\ab}\}$.
\end{lem}
\begin{proof}
   La premi\`ere \'etape est de montrer que $B_0(F)$ est inclus dans l'image du morphism d'inflation $\pi^\ast: \H^2(F^{\ab},\Qbb/\Zbb) \to \H^2(F,\Qbb/\Zbb)$ (o\`u $\pi: F \to F^{\ab}$ d\'esigne la projection). C'est le lemme 3.5 dans {\em loc. cit.}: dans la suite spectrale de Hochschild--Serre $\H^p(F^{\ab},\H^q(Z,\Qbb/\Zbb)) \Rightarrow \H^{p+q}(F,\Qbb/\Zbb)$, le groupe $\H^2(F,\Qbb/\Zbb)$ admet une filtration $E_0 \subseteq E_1 \subseteq E_2$, o\`u $E_0 \subseteq \H^2(F,\Qbb/\Zbb)$ est l'image de $\pi^\ast$, o\`u $E_1/E_0 \subseteq \H^1(F^{\ab},\H^1(Z,\Qbb/\Zbb))$, et o\`u $E_2/E_1 \subseteq \H^2(Z,\Qbb/\Zbb))$. Pour tout $\alpha \in B_0(F)$, son image (par restriction) dans $\H^2(Z,\Qbb/\Zbb)$ est nulle, on peut alors consid\'erer son image $\beta$ dans $\H^1(F^{\ab},\H^1(Z,\Qbb/\Zbb)) = \Hom(F^{\ab},\Hom(Z,\Qbb/\Zbb))$. Il suffit de montrer que $\beta = 0$. En effet, soit $C \subseteq F^{\ab}$ n'importe quel sous-groupe cyclique. On dispose de la suite spectrale de Hochschild--Serre $\H^p(C,\H^q(Z,\Qbb/\Zbb)) \Rightarrow \H^{p+q}(\pi^{-1}(C),\Qbb/\Zbb)$, d'o\`u une filtration $\tilde{E_0} \subseteq \tilde{E_1} \subseteq \tilde{E_2}$ compatible avec celle mention\'ee ci-dessus. Le sous-groupe $\pi^{-1}(C) \subseteq F$ est ab\'elien (c'est une extension centrale du groupe cyclique $C$). La classe $\alpha$ \'etant dans $B_0(F)$, elle se restreint \`a $0 \in \H^2(\pi^{-1}(C),\Qbb/\Zbb)$, donc on a $\beta|_C = 0 \in \H^1(C,\H^1(Z,\Qbb/\Zbb))$. Cela implique que $\beta = 0$, ou $\alpha \in E_0 = \Img \pi^\ast$.
    
    La deuxi\`eme \'etape est de d\'ecrire $\Img \pi^\ast$. Il est connu qu'on a un isomorphisme
        \begin{equation} \label{eqBogomolov1}
            \H^2(F^{\ab},\Qbb/\Zbb) \simeq \Hom\tuple{\bigwedge^2 F^{\ab},\Qbb/\Zbb}
        \end{equation}
    qui \`a la classe $[E]$ de toute extension centrale $0 \to \Qbb/\Zbb \to E \to F^{\ab} \to 0$ associe le morphisme $\lambda_E: \bigwedge^2 F^{\ab} \to \Qbb/\Zbb$ d\'efini par \eqref{eqLambda}. Si $\pi^\ast[E] = 0 \in \H^2(F,\Qbb/\Zbb)$, la projection $\pi: F \to F^{\ab}$ se rel\`eve en un morphisme $\rho: F \to E$. Dans ce cas, $\rho$ induit un morphisme $\chi: Z \to \Qbb/\Zbb$ satisfaisant $\lambda_E = \chi \circ \lambda$. Inversement, supposons qu'il existe un morphisme $\chi: Z \to \Qbb/\Zbb$ tel que $\lambda_E = \chi \circ \lambda$, alors $[E] = \chi_\ast[F]$, o\`u $\chi_\ast: \H^2(F^{\ab},Z) \to \H^2(F^{\ab},\Qbb/\Zbb)$ est le morphism induit par $\chi$. Or $\pi^\ast[F] = 0 \in \H^2(F,Z)$, donc $\pi^\ast[E] = 0 \in \H^2(F,\Qbb/\Zbb)$ puisque les op\'erateurs $\pi^\ast$ et $\chi_\ast$ commutent. On en d\'eduit que sous l'identification de \eqref{eqBogomolov1}, $\Ker \pi^\ast$ correspond au sous-groupe $\Hom(Z,\Qbb/\Zbb) \subseteq \Hom\tuple{\bigwedge^2 F^{\ab},\Qbb/\Zbb}$ (l'inclusion \'etant induit par la surjection $\lambda$). On conclut que $\Img \pi^\ast = \Hom(S,\Qbb/\Zbb)$. 
    
    Soit maintenant $0 \to \Qbb/\Zbb \to E \to F^{\ab} \to 0$ une extension centrale et $\chi:S \to \Qbb/\Zbb$ la restriction de $\lambda_E$. Soient $a,b \in F^{\ab}$ tels que $a \wedge b \in S$ et les relevons respectivement en $\tilde{a},\tilde{b} \in F$. Comme $a \wedge b \in S$, la formule \eqref{eqLambda} implique que $\tilde{a}$ et $\tilde{b}$ commutent, donc le sous-groupe $\pair{\tilde{a},\tilde{b}} \subseteq F$ est ab\'elien. Si $[E] \in B_0(F)$, la compos\'ee $\bigwedge^2 \pair{\tilde{a},\tilde{b}} \to \bigwedge^2 F^{\ab} \xrightarrow{\lambda_E} \Qbb/\Zbb$ est nulle, donc en particulier $\chi(a \wedge b) = \lambda_E(a \wedge b) = 0$. Ainsi, $[E] \in B_0(F)$ implique $\chi|_{S_\lambda} = 0$. Inversement, supposons $\chi|_{S_\lambda} = 0$. Soit $A \subseteq F$ un sous-groupe ab\'elien. Pour tous $\tilde{a},\tilde{b} \in A$, on a $\lambda(a \wedge b) = 0$ (o\`u $a = \pi(\tilde{a})$ et $b = \pi(\tilde{b})$), donc $a \wedge b \in S_{\lambda}$, d'o\`u $\lambda_E(a \wedge b) = \chi(a \wedge b) = 0$. La compos\'ee $\bigwedge^2 A \to \bigwedge^2 F^{\ab} \xrightarrow{\lambda_E} \Qbb/\Zbb$ est nulle, donc $[E] \mapsto 0 \in \H^2(A,\Qbb/\Zbb)$. On en d\'eduit que $[E] \in B_0(F)$. Ainsi
        \begin{equation*}
            B_0(F) = \Ker(\Hom(S,\Qbb/\Zbb) \to \Hom(S_\lambda,\Qbb/\Zbb)) = \Hom(S/S_\lambda,\Qbb/\Zbb).
        \end{equation*}
\end{proof}

\begin{prop} \label{propGeometricBr}
    Soient $M$ et $Z$ des groupes ab\'eliens, $\phi: M \otimes M \to Z$ un morphisme, \`a partir desquels on construit une extension
        \begin{equation*}
            0 \to Z \to F \to M \oplus M \to 0
        \end{equation*}
    de groupes abstraits comme dans le paragraphe \ref{subsection22}, vu comme extension de $K$-groupes finis. Supposons $Z = Z(F) = [F,F]$, soit $F \hookrightarrow \SL_n$ un plongement de $K$-groupes et notons $X = F \backslash \SL_n$. Alors $\Br_{\nr} X = \Hom\left((\Ker \phi) / H, \Qbb/\Zbb\right)$, o\`u $H$ est le sous-groupe $\pair{x \otimes y: x,y \in M, \phi(x \otimes y) = 0}$ de $M \otimes M$.
\end{prop}
\begin{proof}
    Sous l'identification $\bigwedge^2(M \oplus M) = (\bigwedge^2 M) \oplus (M \otimes M) \oplus (\bigwedge^2 M)$, on a 
        \begin{equation*}
            (x,y) \wedge (x',y') = (x \wedge x', x \otimes y', y \wedge y')
        \end{equation*}
    pour tous $x,y,x',y' \in M$. Soit $\lambda = \lambda_F: \bigwedge^2(M \oplus M) \to Z$ le morphisme d\'efini par \eqref{eqLambda}. On note $\Phi: (M \oplus M) \times (M \oplus M) \to Z$ l'application biadditive d\'efinie par \eqref{eqBigPhi}. Par le lemme \ref{lemCommutator}, on a
	\begin{align*}
            \lambda(x \wedge y, 0 , 0) & = \lambda((x,0) \wedge (y,0)) = \Phi((x,0),(y,0)) - \Phi((y,0),(x,0)) = 0,\\ 
            \lambda(0,x \otimes y, 0) & = \lambda((x,0) \wedge (0,y)) = \Phi((x,0),(0,y)) - \Phi((0,y),(x,0)) = \phi(x \otimes y),\\
            \lambda(0 , 0, x \wedge y) & = \lambda((0,x) \wedge (0,y)) = \Phi((0,x),(0,y)) - \Phi((0,y),(0,x)) = 0.
	   \end{align*}
	On en d\'eduit que $\lambda(\gamma_1,\gamma_2,\gamma_3) = \phi(\gamma_2)$ pour tous $\gamma_1,\gamma_3 \in \bigwedge^2 M$ et $\gamma_2 \in M \otimes M$, d'o\`u $S = \Ker \lambda = \bigwedge^2(M \oplus M) = (\bigwedge^2 M) \oplus (\Ker \phi) \oplus (\bigwedge^2 M)$, au vu des notations du lemme \ref{lemBogomolovFormula}.
	
    Affirmons que $S_\lambda = \bigwedge^2(M \oplus M) = (\bigwedge^2 M) \oplus H \oplus (\bigwedge^2 M)$. En effet, soient $x,y,x',y' \in M$. Alors $\lambda((x,y) \wedge (x',y')) = \lambda'(x \wedge x', x \otimes y', y \wedge y') = \phi(x \otimes y')$, donc $(x,y) \wedge (x',y') \in S$ si et seulement si $\phi(x \otimes y') = 0$. Dans ce cas, $x \otimes y' \in H$ par d\'efinition, d'o\`u $S_\lambda \subseteq (\bigwedge^2 M) \oplus H \oplus (\bigwedge^2 M)$. Inversement,  
	\begin{itemize}
            \item $(\bigwedge^2 M) \oplus 0 \oplus 0 \subseteq S_\lambda$ puisque $(x \wedge y, 0,0) = (x,0) \wedge (y,0) \in S$ pour tous $x,y \in M$;
            
            \item $0 \oplus 0 \oplus (\bigwedge^2 M) \subseteq S_\lambda$ puisque $(0,0,x \wedge y) = (0,x) \wedge (0,y) \in S$ pour tous $x,y \in M$;
            
            \item $0 \oplus H \oplus 0 \subseteq S_\lambda$ puisque $(0,x \otimes y, 0) = (x,0) \wedge (0,y) \in S$ pour tous $x,y \in M$ tels que $\phi(x \otimes y) = 0$;
	\end{itemize}
    d'o\`u $S_\lambda = (\bigwedge^2 M) \oplus H \oplus (\bigwedge^2 M)$ comme voulu. Finalement, $S/S_\lambda = (\Ker \phi) / H$ et donc la proposition \ref{propBogomolovFormula} donne $\Br_{\nr} X = B_0(F) = \Hom(S/S_\lambda,\Qbb/\Zbb) = \Hom((\Ker \phi) / H,\Qbb/\Zbb)$.
\end{proof}
  
\section{Les principaux r\'esultats} \label{section5}

\subsection{Un lemme arithm\'etique} \label{subsection51}

Pour \'etablir les principaux r\'esultats de ce texte, une propri\'et\'e globale du symbole de Hilbert est n\'ecessaire. On propose la g\'en\'eralisation suivante de \cite[Chapitre I, \S 2.2, Th\'eor\`eme 4]{Serre}.

\begin{prop} [\og Lemme arithm\'etique \fg{} pour les groupes cycliques] \label{propSerre}
    Soit $k$ un corps de nombres, soit $n$ un entier et soit $(a_i)_{i \in I}$ une famille finie d'\'el\'ements de $\H^1(k,\Zbb/n)$. Soit $(\lambda_{i})_{i\in I} \in (\Br k)^{|I|}$ remplissant la condition suivante: pour toute $v \in \Omega_k$, il existe $c_v \in \H^1(k_v,\mu_n)$ tel que $(a_i)_v \cup c_v = (\lambda_{i})_v$ pour tout $i \in I$. Alors, pour tout sous-ensemble fini $S \subseteq \Omega_k$, il existe $b \in \H^1(k,\mu_n)$ tel que
	\begin{enumerate}
            \item $a_i \cup b = \lambda_i$ pour tout $i \in I$;
            \item $b_v = c_v$ pour toute $v \in S$.
	\end{enumerate}
\end{prop}

La preuve da la proposition \ref{propSerre} repose sur des th\'eor\`emes de dualit\'e arithm\'etique. Pour les d\'etails de la th\'eorie de la dualit\'e arithm\'etique, on pourra consulter \cite[Chapitre 17]{harari2017cohomologie}. Rappelons alors quelques notations. Soit $k$ un corps de nombres et soit $A$ un $\Gamma_k$-module fini. Pour toute place finie $v$ de $k$ telle que $A$ soit un $\Gamma_{k_v}$-module non ramifi\'e, on note $\H^1_{\nr}(k_v,A)$ l'image de la fl\`eche d'inflation $\H^1(\Gal(k_{v,\nr}/k_v), A) \to \H^1(k_v,A)$. Notons $\Pbb^1(k,A)$ le produit restreint $\prod_{v \in \Omega_k}' \H^1(k_v,A)$ par rapport aux sous-groupes $\H^1_{\nr}(k_v,A)$. On dispose d'une application diagonale de $\H^1(k,A)$ dans $\Pbb^1(k,A)$. En rassemblent les dualit\'es locales de Tate entre $\H^1(k_v,A)$ et $\H^1(k_v,\hat{A})$, on obtient un accouplement parfait de groupes ab\'eliens localement compacts
    \begin{equation*}
        \Pbb^1(k,A) \times \Pbb^1(k,\hat{A}) \to \Qbb/\Zbb, \quad ((a_v)_{v \in \Omega_k}, (b_v)_{v \in \Omega_k}) \mapsto \sum_{v \in \Omega_k} \inv_v(a_v \cup b_v).
    \end{equation*}

\begin{proof} [D\'emonstration de la proposition \ref{propSerre}]
    Calculons le sous-groupe 	
    \begin{equation*}
        P = \Img(\H^1(k,\Zbb/n) \to \Pbb^1(k,\Zbb/n)) \cap \left(\prod_{v \in S} \H^1(k_v,\Zbb/n) \times \prod_{v \in \Omega_k \setminus S} \pair{(a_i)_v: i \in I} \right)
    \end{equation*}
    de $\Pbb^1(k,\Zbb/n)$. Soit alors $a \in \H^1(k,\Zbb/n)$ tel que pour toute $v \notin S$, $a_v$ soit une combinaison lin\'eaire des $(a_i)_v$, $i \in I$. On voit les $a_i$ et $a$ comme des morphismes continus $\Gamma_k \to \Zbb/n$; alors il existe un quotient $G$ de $\Gamma_k$, qui est fini, ab\'elien et de $n$-torsion, par lequel ces morphismes se factorisent. Notons $a_i',a' \in \Hom(G,\Zbb/n)$ les morphismes induits respectifs. La condition sur $a$ et le th\'eor\`eme de Chebotarev implique que $a'|_H \in \pair{a'_i|_H: i \in I}$ pour tout sous-groupe cyclique $H$ de $G$. En appliquant cette condition \`a chaque sous-groupe cyclique de $\bigcap_{i \in I} \Ker(a'_i)$, on voit que $a'$ s'annule sur $\bigcap_{i \in I} \Ker(a'_i)$. Au vu de l'accouplement parfait entre $G$ et $\Hom(G,\Zbb/n)$, $a'$ est orthogonal \`a $\bigcap_{i \in I} \Ker(a'_i)$, qui est l'orthogonal de $\pair{a'_i: i \in I}$, ainsi $a' \in \pair{a'_i: i \in I}$. Il s'ensuit que $a \in \pair{a_i: i \in I}$ et on conclut que $P = \pair{((a_i)_v)_{v \in \Omega_k}: i \in I}$.
	
    Lorsque $v \in \Omega_k \setminus S$ est une place telle que $(\lambda_i)_v = 0$ pour tour $i \in I$, on peut supposer que $c_v = 0$. Cela implique que $(c_v)_{v \in \Omega_k} \in \bigoplus_{v \in \Omega_k} \H^1(k_v,\mu_n) \subseteq \Pbb^1(k,\mu_n)$. Au vu de l'accouplement parfait entre $\Pbb^1(k,\Zbb/n)$ et $\Pbb^1(k,\mu_n)$, la loi de r\'eciprocit\'e globale implique que $(c_v)_{v \in \Omega_k}$ est orthogonal \`a $((a_i)_v)_{v \in \Omega_k}$ pour tout $i \in I$, {\em i.e.} il est orthogonal \`a $P$. Or l'orthogonal de $\Img(\H^1(k,\Zbb/n) \to \Pbb^1(k,\Zbb/n))$ est $\Img(\H^1(k,\mu_n) \to \Pbb^1(k,\mu_n))$ (c'est l'exactitude au 5$^e$ terme de la suite exacte \`a 9 termes de Poitou--Tate, {\em cf.} \cite[Th\'eor\`eme 17.13]{harari2017cohomologie}), donc 
	\begin{equation*}
            P^{\perp} = \Img(\H^1(k,\mu_n) \to \Pbb^1(k,\mu_n)) + \prod_{v \in S} \{0\} \times \prod_{v \in \Omega_k \setminus S}' \pair{(a_i)_v: i \in I}^{\perp},
        \end{equation*}
    donc il existe $b \in \H^1(k,\mu_n)$ satisfaisant les conditions suivantes.
	\begin{enumerate}
            \item Pour toute $v \in S$, $b_v - c_v = 0$ pour toute $v \in S$.
            \item Pour toute $v \notin S$, $b_v - c_v$ est orthogonal aux $(a_i)_v$, $i \in I$.
	\end{enumerate}
    On a ainsi que $(a_i \cup b)_v = (a_i)_v \cup c_v = (\lambda_{i})_v$ pour tous $i \in I$ et $v \in \Omega_k$ (d'o\`u $a_i \cup b = \lambda_i$ par la loi de r\'eciprocit\' e globale) et que $b_v = c_v$ pour toute $v \in S$. La proposition est finalement d\'emontr\'ee.
\end{proof}

La proposition \ref{propSerre} se g\'en\'eralise en la proposition \ref{propArithmetic} ci-dessous. Afin de la d\'emontrer, remarquons le fait suivant: Soient $m$ et $n$ deux entiers, soit $d = \PGCD(m,n)$, et soit $K$ un corps de caract\'eristique nulle. Alors l'on dispose d'un diagramme commutatif de $\Gamma_K$-modules, dont les fl\`eches horizontales sont des isomorphismes:
    \begin{equation*}
        \xymatrix{
		\Zbb/m \otimes \mu_n \ar[d]^{\iota^m_n \otimes \id} \ar[rr]^{\simeq} && \mu_d  \ar@{^{(}->}[d] \\
		\Zbb/n \otimes \mu_n \ar[rr]^{\simeq} && \mu_n,
	}
    \end{equation*}
et dont le morphisme $\iota^m_n$ est la compos\'ee 
    \begin{equation} \label{eqIotaMN}
        \Zbb/m \twoheadrightarrow \Zbb/d \hookrightarrow \Zbb/n.
    \end{equation}
Il s'ensuit qu'il y a un diagramme commutatif
    \begin{equation} \label{eqPairingsMN}
        \xymatrix{
            \H^1(K,\Zbb/m) \times \H^1(K,\mu_n) \ar@<-6.5ex>[d]^{(\iota^m_n)_\ast} \ar@{=}@<+6.5ex>[d] \ar[rr]^-{\cup} && (\Br K)[d] \ar@{^{(}->}[d] \\
		\H^1(K,\Zbb/n) \times \H^1(K,\mu_n)  \ar[rr]^-{\cup} && (\Br K)[n].
        }
    \end{equation}
D'ailleurs, dans le cas o\`u $K$ contient $\mu_{\PPCM(m,n)}$ (un g\'en\'erateur duquel sera fix\'e), on dispose d'un diagramme commutatif de $\Gamma_K$-modules, \`a lignes exactes:
    \begin{equation*}
        \xymatrix{
        1 \ar[r] & \Zbb/m \ar[rr] \ar@{->>}[d] && \ol{K}^\times  \ar[rr]^{(-)^{m}} \ar[d]^{(-)^{m/d}} && \ol{K}^\times \ar[r] \ar@{=}[d] & 1 \\
        1 \ar[r] & \Zbb/d \ar[rr] \ar@{^{(}->}[d] && \ol{K}^\times \ar[rr]^{(-)^{d}} \ar@{=}[d] && \ol{K}^\times \ar[r] \ar[d]^{(-)^{n/d}} & 1\\
        1 \ar[r] & \Zbb/n \ar[rr] && \ol{K}^\times \ar[rr]^{(-)^{n}} && \ol{K}^\times \ar[r] & 1,
    }
    \end{equation*}
d'o\`u un diagramme commutatif
    \begin{equation} \label{eqKummerMN}
        \xymatrix{
            K^\times \ar[rr] \ar[d]^{(-)^{n/d}} && \H^1(K,\Zbb/m) \ar[d]^{(\iota^m_n)_\ast} \\
             K^\times \ar[rr] && \H^1(K,\Zbb/n),
        }
    \end{equation}
dont les fl\`eches horizontales sont des morphismes de la th\'eorie de Kummer.

\begin{prop} [\og Lemme arithm\'etique \fg{} \label{propArithmetic} pour les groupes ab\'eliens finis]
    Soient $k$ un corps de nombres, $A$ un groupe ab\'elien fini muni de l'action triviale de $\Gamma_k$, et $(a_i)_{i \in I}$ une famille finie d'\'el\'ements de $\H^1(k,A)$. Soit $(\lambda_{i})_{i\in I} \in \H^2(k,A \otimes \hat{A})^{|I|}$ remplissant la condition suivante: pour toute $v \in \Omega_k$, il existe $c_v \in \H^1(k_v,\hat{A})$ tel que $(a_i)_v \cup c_v = (\lambda_{i})_v$ pour tout $i \in I$. Alors pour tout sous-ensemble fini $S \subseteq \Omega_k$, il existe $b \in \H^1(k,\hat{A})$ tel que
	\begin{enumerate}
            \item $a_i \cup b = \lambda_i$ pour tout $i \in I$;
            \item $b_v = c_v$ pour toute $v \in S$.
	\end{enumerate}
\end{prop}
\begin{proof}
    \'Ecrivons $A = \prod_{p \in J} \Zbb/{n_p}$ (alors $\hat{A} = \prod_{q \in J} \mu_{n_q}$), et $d_{p,q} = \PGCD(n_p,n_q)$ pour tous $p,q \in J$. De plus, pour tous $i \in I$ et $v \in \Omega_k$, \'ecrivons
        \begin{align*}
            a_i & = (a_i^p)_{p \in J} \in \H^1(k,A) = \prod_{p \in J}\H^1(k,\Zbb/n_p),\\
            c_v & = (c_v^q)_{q \in J} \in \H^1(k_v,\hat{A}) = \prod_{q \in J}\H^1(k_v,\mu_{n_q}),\\
            \lambda_i & = (\lambda_i^{p,q})_{p,q \in J} \in \H^2(k,A \otimes \hat{A}) = \prod_{p,q \in J}(\Br k)[d_{p,q}].
        \end{align*}
    La condition $(a_i)_v \cup c_v = (\lambda_i)_v$ se r\'e\'ecrit sous la forme
        \begin{equation} \label{eqpropArithmetic}
            \forall p,q \in J, \quad ((\iota^{n_p}_{n_q})_\ast a_i^p)_v \cup c_v^{q} = (a_i^p)_v \cup c_v^{q} = (\lambda_i^{p,q})_v,
        \end{equation}
    par biadditivit\'e des cup-produits et au vu de \eqref{eqPairingsMN}, o\`u $\iota^{n_p}_{n_q}: \Zbb/n_p \to \Zbb/n_q$ est d\'efini comme la compos\'ee \eqref{eqIotaMN}. Fixons $q \in J$ et appliquons la proposition \ref{propSerre} \`a la famille $((\iota^{n_p}_{n_q})_\ast a_i^p)_{i \in I, p \in J}$ pour trouver un $b^{q} \in \H^1(k,\mu_{n_q})$ tel que
        \begin{enumerate}
            \item $a_i^p \cup b^{q} = \lambda_i^{p,q}$ pour tous $i \in I$ et $p \in J$;
            \item $b^{q}_v = c^{q}_v$ pour toute $v \in S$.
        \end{enumerate}
    Posons finalement $b = (b^{q})_{q \in J} \in \H^1(k,\hat{A})$. Alors $a_i \cup b = \lambda_i$ (par biadditivit\'e des cup-produits), et $b_v = c_v$ pour toute $v \in S$, ce qui ach\`eve la d\'emonstration.
\end{proof}

Afin d'appliquer la proposition \ref{propArithmetic}, il convient de prouver l'\'enonc\'e suivant.

\begin{lem} \label{lemHilbertSymbol}
    Soient $(K,v)$ un corps $p$-adique et $n$ un entier tels que $K$ contienne $\mu_n$. Soit $d$ un diviseur de $n$ et soit $\tilde{a} \in K^\times$ tel que $v(\tilde{a})$ soit premier \`a $d$. Soit $a$ son image dans $\H^1(k,\mu_n)$ par la th\'eorie de Kummer. Alors pour tout $r \in \frac{1}{d}\Zbb/\Zbb$, il existe $b \in \H^1(K,\Zbb/n)$ tel que $\inv_K(a \cup b) = r$.
\end{lem}
\begin{proof}
    Posons $A:=\{\inv_K(a \cup b): b \in \H^1(K,\Zbb/n)\} \subseteq \Qbb/\Zbb$ et montrons que $A$ contient $\frac{1}{d}$. Commen\c cons par le cas o\`u $d = n$. Il suffit de montrer que pour tout nombre premier $\ell$ divisant $n$, $A$ contient $\frac{1}{\ell^{v_\ell(n)}}$. En effet, comme $\ell$ ne divise pas $v(\tilde{a})$, on a $(\tilde{a})^{\frac{n}{\ell}} \notin K^{\times n}$ et donc $\frac{n}{\ell}a \neq 0$. Par dualit\'e locale de Tate, il existe $b \in \H^1(k,\Zbb/n)$ tel que $\inv_K(\frac{n}{\ell}a \cup b) \neq 0$ dans $\Qbb/\Zbb$. \'Ecrivons $\inv_K(a \cup b) = \frac{t}{n}$ avec $t \in \Zbb$, alors $\inv_K(\frac{n}{\ell}a \cup b) = \frac{n}{\ell}\inv_K(a \cup b) = \frac{t}{\ell}$ est non nul dans $\Qbb/\Zbb$, d'o\`u $\ell$ ne divise pas $t$. Or $A$ contient $\inv_K(a \cup \frac{n}{\ell^{v_\ell(n)}}b) = \frac{n}{\ell^{v_\ell(n)}} \inv_K(a \cup b) = \frac{t}{\ell^{v_\ell(n)}}$, donc il contient $\frac{1}{\ell^{v_\ell(n)}}$. Ainsi, $A$ contient $\frac{1}{n}$ comme voulu.
	
    Revenons au cas g\'en\'eral. On a un diagramme commutatif de $\Gamma_K$-modules
	\begin{equation*}
            \xymatrix{
                1 \ar[r] & \mu_n \ar[r] \ar[d]^{(-)^{n/d}} & \ol{K}^\times \ar[r]^{(-)^n} \ar[d]^{(-)^{n/d}} & \ol{K}^\times \ar[r] \ar@{=}[d] & 1 \\
                1 \ar[r] & \mu_d \ar[r] & \ol{K}^\times \ar[r]^{(-)^d} & \ol{K}^\times \ar[r] & 1,
            }
        \end{equation*}
    dont les deux lignes sont exactes. D'o\`u un diagramme commutatif
	\begin{equation*}
            \xymatrix{
                K^\times \ar[r] \ar@{=}[d] & \H^1(K,\mu_n) \ar[d] \\
                K^\times \ar[r] & \H^1(K,\mu_d).
            }
        \end{equation*}
    D'o\`u l'image $a'$ de $a \in \H^1(K,\mu_n)$ par le morphisme induit par $\mu_n \xrightarrow{(-)^{n/d}} \mu_d$ co\"incide avec l'image de $\tilde{a}$ dans $\H^1(K,\mu_d)$ par la th\'eorie de Kummer. Comme $v(\tilde{a})$ est premier \`a $d$, par le r\'esultat dans le cas pr\'ec\'edent, il existe $b' \in \H^1(K,\Zbb/d)$ tel que $\inv_K(a' \cup b') = \frac{1}{d}$. On a un diagramme commutatif
	\begin{equation*}
            \xymatrix{
                \mu_n \times \Zbb/n \ar@<-2ex>[d]_{(-)^{n/d}} \ar[rr] && \ol{K}^\times \\
                \mu_d \times \Zbb/d \ar@{_{(}->}[u]<-2ex>  \ar[rru].
            }
        \end{equation*}
    Notant $b \in \H^1(K,\Zbb/n)$ l'image de $b'$ par le morphisme induit par $\Zbb/d \hookrightarrow \Zbb/n$, on a $a \cup b = a' \cup b'$, d'o\`u $A$ contient $\inv_K(a \cup b) = \inv_K(a' \cup b') = \frac{1}{d}$. Le lemme est finalement d\'emontr\'e.
\end{proof}

\subsection{Construction du morphisme $\phi: M \otimes M \to Z$} \label{subsection52}

Comme on a vu dans le paragraphe \ref{subsection22}, on peut construire une extension 
    \begin{equation*}
        0 \to Z \to F \to M \oplus M \to 0
    \end{equation*}
\`a partir d'un morphisme $\phi: M \otimes M \to Z$. Nous allons maintenant choisir $M$, $Z$ et $\phi$. Les notations suivantes seront fix\'ees jusqu'\`a la fin du texte.

\begin{itemize}
	\item $A$ est un groupe ab\'elien fini et $k$ est un corps de nombres contenant $\mu_{\exp(A)}$. On munit $A$ de l'action triviale de $\Gamma_k$ et l'on fixe un g\'en\'erateur de $\mu_{\exp(A)}$ (ce qui d\'efinit un isomorphisme $A = \hat{A}$). $L/k$ est une extension finie galoisienne et $\gfrak = \Gal(L/k)$ (dans la construction originelle de \cite{Bokun}, $A = \Zbb/p^2$ et $\gfrak = \Zbb/p \times \Zbb/p$, o\`u $p$ est un nombre premier).
	
	\item Pour toute place $v$ de $k$, une place $w_v$ de $L$ divisant $v$ est choisie et $\gfrak_v \subseteq \gfrak$ est le groupe de d\'ecomposition de $w_v|v$. On fixe \'egalement un syst\`eme de repr\'esentants $\Escr_v$ des classes \`a gauche de $\gfrak$ suivantes $\gfrak_v$ et l'on pose $e_v = |\Escr_v| = [\gfrak:\gfrak_v]$. Alors pour toute place $w$ de $L$ divisant $v$, il existe un unique $s \in \Escr_v$ tel que $w = sw_v$. De plus, si $\sigma \in \Gamma_k$ est un relev\'e de $s$, alors $\Gamma_{L_w} = \sigma \Gamma_{L_{w_v}} \sigma^{-1}$. De plus, pour tout $\Gamma_k$-module $B$ sur lequel $\Gamma_L$ agit trivialement et tout $r \ge 1$, on dispose d'un isomorphisme
	\begin{equation*}
            \Z^r(L_{w_v}, B) \to \Z^r(L_w, B),
        \end{equation*}
	qui \`a chaque $r$-cocycle $c: \Gamma_{L_{w_v}}^r \to B$ associe le $r$-cocycle
	\begin{equation*}
            \Gamma_{L_w}^r \to B, \quad (\tau_1,\ldots,\tau_r) \mapsto c(\sigma^{-1}\tau_1\sigma,\ldots,\sigma^{-1}\tau_r\sigma).
        \end{equation*}
    Cet isomorphisme-l\`a induit un isomorphisme $\H^r(L_{w_v}, B) \to \H^r(L_w, B)$, qui est compatible avec les fl\`eches de restriction de $\H^r(k_v,B)$. On va noter l'image de chaque classe $\gamma \in \H^r(L_{w_v},B)$ par $\tensor[^s]{\gamma}{} \in \H^r(L_w,B)$. Si $\gamma \in \H^r(L,B)$, alors $\tensor[^s]{(\gamma}{_{w_v}}) =  (\tensor[^s]{\gamma}{})_{_{w}} \in \H^r(L_w,B)$.
	
	\item $M = \Ind_{\Gamma_k}^{\Gamma_L} A$ et $j: A \otimes A \hookrightarrow M \otimes M$ est l'inclusion canonique, {\em cf.} lemme \ref{lemShapiroJ}. 
	
	\item $Z$ est le conoyau de $j$ et $\phi: M \otimes M \to Z$ est la projection. $\Phi: (M \oplus M) \times (M \oplus M) \to Z$ est l'application biadditive d\'efinie par \eqref{eqBigPhi}, {\em i.e.}
	\begin{equation*}
            \forall x,y,x',y' \in M, \quad \Phi((x,y),(x',y')) = \phi(x \otimes y').
        \end{equation*}
	Soit $F$ le produit crois\'e $Z \times_\Phi (M \oplus M)$. On verra que $Z = Z(F) = [F,F]$ (proposition \ref{propBorovoiKunyavskii}). On munit $F$ d'une action de $\Gamma_k$ compatible avec celles sur $Z$ et sur $M \otimes M$.
\end{itemize}

\begin{prop} \label{propBorovoiKunyavskii}
    Avec les notations ci-dessus, on a les r\'esultats suivants.
	\begin{enumerate}
		\item $Z = Z(F) = [F,F]$.
            \item Soit $X$ un espace homog\`ene de $\SL_n$ \`a stabilisateur g\'eom\'etrique $F$. Alors $\Be(X) = \Be_\omega (X) = 0$ et $\Br_{\nr} X = \Br_{\nr,1} X$. Si de plus l'action ext\'erieure de $\Gamma_k$ sur $F$ est induite par l'action coordonn\'ees par coordonn\'ees, alors $\Br_{\nr} X = \Br_0 X$.
	\end{enumerate}
\end{prop}
\begin{proof}
    On rappelle que $M = \{\gfrak \to A\}$ (le $\Gamma_k$-module des applications $\gfrak \to A$) et que $M \otimes M = \{\gfrak \times \gfrak \to A \otimes A\}$. 
    \begin{enumerate}
        \item Par le lemme \ref{lemCommutator}, $Z = [F,F]$ puisque $\phi$ est surjectif. Montrons que $\phi$ est non d\'eg\'en\'er\'e.

        Notons d'abord que si $a \in A$ est tel que $a \otimes b = 0$ pour tout $b \in A$, alors $a = 0$. En effet, en choisissant un isomorphisme $A \simeq \Hom(A,\Zbb/\exp(A))$, on obtient un accouplement parfait $A \otimes A \to \Zbb/\exp(A)$.

        Supposons maintenant que $x: \gfrak \to A$ est tel que $\phi(x \otimes y) = 0$ pour tout $y: \gfrak \to A$. Pour tout $a \in A$, soit $y_a: \gfrak \to A$ d\'efini par $y_a(1) = a$ et $y_a(g) = 0$ pour tout $g \neq 1$. Alors $\phi(x \otimes y_a) = 0$, donc il existe $m_a \in A \otimes A$ tel que $x \otimes y_a = j(m_a)$, c'est-\`a-dire que
            \begin{equation*}
                \forall g,h \in \gfrak, \quad x(g) \otimes y_a(h) = m_a,
            \end{equation*}
        d'o\`u $x(g) \otimes a = m_a = x(g) \otimes 0 = 0$ pour tous $g \in \gfrak$ et $a \in A$, donc $x(g) = 0$, ainsi $x = 0$. De m\^eme, si $y \in M$ est tel que $\phi(x \otimes y) = 0$ pour tout $x \in M$, on v\'erifie sans peine que $y = 0$. Ainsi, $\phi$ est bien non d\'eg\'en\'er\'e et donc $Z = Z(F)$ au vu du lemme \ref{lemCenterOfF}.

        \item Calculons $\Br_{\nr} \ol{X}$. Au vu de la proposition \ref{propGeometricBr}, ce groupe est $\Hom((\Ker \phi) / H, \Qbb/\Zbb)$, o\`u $H$ est le sous-groupe $\pair{x \otimes y: x,y \in M, \phi(x \otimes y) = 0}$ de $M \otimes M$. Tout \'el\'ement de $\Ker \phi$ est de la forme $j(m)$ pour un $m \in A \otimes A$. \'Ecrivons $m = \sum_{i \in I} a_i \otimes b_i$, o\`u $I$ est fini et $a_i,b_i \in A$. On d\'efinit $x_i,y_i: \gfrak \to A$ par
            \begin{equation*}
                x_i(g) = a_i, \quad y_i(g) = b_i
            \end{equation*}
        pour tous $g \in \gfrak$ et $i \in I$. Alors $j(a_i \otimes b_i) = x_i \otimes y_i$, d'o\`u $\phi(x_i \otimes y_i) = 0$. On en d\'eduit que 
            \begin{equation*}
                j(m) = \sum_{i \in I} j(a_i \otimes b_i) = \sum_i x_i \otimes y_i \in H,
            \end{equation*}
        d'o\`u $\Ker \phi = H$, donc $\Br_{\nr} \ol{X} = 0$. Il s'ensuit que $\Br_{\nr,1} X = \Br_{\nr} X$ et $\Br_{\nr,a} X = (\Br_{\nr} X)/(\Br_0 X)$.
		
        Calculons $\Be_\omega(X)$. Par le lemme de Shapiro, $\Sha^1_\omega(k,\hat{M}) = \Sha^1_\omega(L,\hat{A}) = \Sha^1_\omega(L,A)$. Or $\Sha^1_\omega(L,A) = 0$ par une application du th\'eor\`eme de Chebotarev ({\it cf.} \cite[Corollaire 18.4]{harari2017cohomologie}), donc $\Sha^1_\omega(k,\hat{M}) = 0$. Par la proposition \ref{propBrauerGroups}, on a $\Be_\omega(X) = 0$ et {\it a fortiori} $\Be(X) = 0$.
		
        Finalement, si l'action ext\'erieure de $\Gamma_k$ sur $F$ est induite par l'action coordonn\'ees par coordonn\'ees, alors  $\Br_{\nr,a} X = \Be_\omega(X) = 0$ par la proposition \ref{propBrNr}, d'o\`u $\Br_{\nr} X = \Br_0 X$.
    \end{enumerate}
\end{proof}

Pour les principaux th\'eor\`emes, on consid\`ere seulement l'action coordonn\'ees par coordonn\'ees de $\Gamma_k$ sur $F$. On suppose ceci pour le reste du texte.

\subsection{Principe de Hasse} \label{subsection53}

Le th\'eor\`eme suivant sera le c\oe ur des preuves de nos principaux r\'esultats. Son deuxi\`eme point sera utile pour l'\'etude de l'approximation faible. Les lecteurs souhaitant s'habituer \`a l'id\'ee de sa preuve sont conseill\'es de se restreindre au cas o\`u $A = \Zbb/n$.

\begin{thm} \label{thmKey}
	Soient $A$ un groupe ab\'elien fini, $k$ un corps de nombres contenant $\mu_{\exp(A)}$, $L/k$ une extension finie galoisienne, $M = \Ind_{\Gamma_k}^{\Gamma_L} A$, $j: A \otimes A \hookrightarrow M \otimes M$ l'inclusion canonique, $Z$ le conoyau de $j$ et $\phi: M \otimes M \to Z$ la projection. De $\phi$ on d\'efinit une application biadditive $\Phi: (M \oplus M) \times (M \oplus M) \to Z$ par \eqref{eqBigPhi} et l'on pose $F = Z \times_\Phi (M \oplus M)$, muni de l'action coordonn\'ees par coordonn\'ees de $\Gamma_k$. Soit $u \in \H^2(k, M \otimes M)$. Alors, \'etant donn\'e:
		\begin{itemize}
			\item un sous-ensemble fini $S \subseteq \Omega_k$ tel que $u_v = 0$ pour toute $v \notin S$;
			\item pour toute $v \in S$, des \'el\'ements $x'_v,y'_v \in \H^1(k_v,M)$ et $\lambda'_v \in \H^2(k_v,A \otimes A)$ tels que $x'_v \cup y'_v = u_v + j_\ast \lambda'_v$ dans $\H^2(k_v, M \otimes M)$;
		\end{itemize}
	il existe $x,y \in \H^1(k,M)$ et $\lambda \in \H^2(k,A \otimes A)$ tels que les conditions suivantes soient satisfaites.
	\begin{enumerate}
		\item $x \cup y = u + j_\ast \lambda$ dans $\H^2(k,M \otimes M)$.
		\item $x_v = x'_v$, $y_v = y_v'$ et $\lambda_v = \lambda'_v$ pour toute $v \in S$.
	\end{enumerate}
\end{thm}
\begin{proof}
	Gardons les notations au d\'ebut du paragraphe \ref{subsection42}. Par le lemme \ref{lemShapiroTensor}, on dispose d'un isomorphisme $\sh': \H^2(k, M \otimes M) \to \H^2(L, A \otimes A)^{|\gfrak|}$; \'ecrivons 
	   \begin{equation*}
                \sh'(u) = (\gamma_g)_{g \in \gfrak} \in \H^2(L, A \otimes A)^{|\gfrak|}.
            \end{equation*}
	Par le lemme \ref{lemShapiroLocalization2}, on a
	    \begin{equation*}
                \forall v \notin S, \forall h \in \gfrak_v,\forall s,t \in \Escr_v, \quad (\tensor[^{t^{-1}}]{\gamma}{_{sht^{-1}}})_{w_v} = 0 \quad \text{dans } \H^2(L_{w_v}, A \otimes A),
            \end{equation*}
	ou $(\gamma_{sht^{-1}})_{tw_v} = 0$ dans $\H^2(L_{tw_v}, A \otimes A)$. Pour tous $g \in \gfrak$ et $w|v$, on peut \'ecrire $w = tw_v$ pour un $t \in \Escr_v$, puis $gt = sh$ pour un $s \in \Escr_v$ et un $h \in \gfrak_v$. On a ainsi
	\begin{equation} \label{eqthmKey1}
		\forall v \notin S, \forall g \in \gfrak,\forall w|v, \quad (\gamma_g)_w = 0 \quad \text{dans } \H^2(L_{w}, A \otimes A).
	\end{equation}
	On proc\`ede en plusieurs \'etapes.
	
	{\bf \'Etape 1.} {\em Approchons les $\lambda'_v$.} \'Ecrivons $A = \prod_{p \in J} \Zbb/n_p = \prod_{p \in J} \mu_{n_p}$, et $d_{p,q} = \PGCD(n_p,n_q)$ pour tous $p,q \in J$. Par le th\'eor\`eme de Chebotarev, il existe $|J|^2$ places (deux \`a deux distinctes) $v_{p,q} \notin S$ de $k$ ($p,q \in J$) qui sont finies et totalement d\'ecompos\'ees dans $L$ (donc $\Escr_{v_{p,q}} = \gfrak$). Pour toute $v \in S$, \'ecrivons 
	\begin{equation*}
            \lambda'_v = (\lambda'_v{}^{p,q})_{p,q \in J} \in \H^2(k_v, A \otimes A) = \prod_{p,q \in J} (\Br k_v)[d_{p,q}].
        \end{equation*}
	Pour tous $p,q \in J$, par la loi de r\'eciprocit\'e globale, il existe $\lambda^{p,q} \in (\Br k)[d_{p,q}]$ tel que 
	\begin{enumerate}
		\item $\lambda_v^{p,q} = \lambda'_v{}^{p,q}$ pour toute $v \in S$;
		\item $\inv_{v_{p,q}}(\lambda_{v_{p,q}}^{p,q}) = -\sum_{v \in S} \inv_v(\lambda'_v{}^{p,q})$;
		\item $\lambda_v^{p,q} = 0$ pour toute $v \notin S \cup \{v_{p,q}\}$.
	\end{enumerate}
	Posons $\lambda = (\lambda^{p,q})_{p,q \in J} \in  \prod_{p,q \in J} (\Br k)[d_{p,q}] = \H^2(k, A \otimes A).$ Alors $\lambda_v = \lambda'_v$ pour toute $v \in S$ et $\lambda_v = 0$ pour toute $v \notin S \cup \{v_{p,q}: p,q \in J\}$. En particulier, on a
	\begin{equation} \label{eqthmKey2}
		\forall v \in S, \quad x_v' \cup y_v' = u_v + j_\ast \lambda_v \quad \text{dans } \H^2(k_v, M \otimes M).
	\end{equation}
	
	{\bf \'Etape 2.} {\em Approchons les $x'_v$.} Pour toute $v \in S$, on dispose d'un isomorphisme 
		\begin{equation*}
                \sh_v: \H^1(k_v,M) \xrightarrow{\simeq} \H^1(L_{w_v},A)^{e_v}
            \end{equation*}
	(lemme \ref{lemShapiroLocalM}). \'Ecrivons $\sh_v(x'_v) = (a'_{v,s})_{s \in \Escr_v} \in \H^1(L_{w_v},A)^{e_v}$ et $\sh_v(y'_v) = (b'_{v,s})_{s \in \Escr_v} \in \H^1(L_{w_v},A)^{e_v}$. 
	
	Pour tous $p,q \in J$, on choisit une uniformisante $\varpi_{p,q} \in k^\times$ de $k_{v_{p,q}}$. Pour tout $p \in J$, le th\'eor\`eme des restes chinois donne un \'el\'ement $\tilde{a}^p \in L^\times$ tel que
	    \begin{enumerate}
	        \item $\tilde{a}^p_w = \varpi_{p,q} \pmod {L_w^{\times n_p}}$ pour toute place $w$ de $L$ divisant une place $v_{p,q}$, o\`u $q \in J$;
	        
	        \item $\tilde{a}^p_w \in L_w^{\times n_p}$ pour toute place $w$ de $L$ divisant une place $v_{p',q}$, o\`u $p' \in J \setminus \{p\}$ et $q \in J$.
	    \end{enumerate}
	On note $a'{}^p \in \H^1(L,\mu_{n_p}) = \H^1(L,\Zbb/n_p)$ l'image de $\tilde{a}^p$ par la th\'eorie de Kummer. Comme $\Sha^1_\omega(L,A) = 0$ par le th\'eor\`eme de Chebotarev, le $\Gamma_L$-module $A$ v\'erifie l'approximation faible au sens du lemme \ref{lemWA} ({\em cf.} \cite[Exercice 17.5(a), Corollaire 18.4]{harari2017cohomologie}), donc il existe un $a \in \H^1(L,A)$ tel que
	    \begin{enumerate}
	        \item $a_{sw_v} = \tensor[^s]{a}{}'_{v,s} \in \H^1(L_{sw_v},A)$ pour tous $v \in S$ et $s \in \Escr_v$;
	        
	        \item $a_w = (a'{}^p_w)_{p \in J} \in \prod_{p \in J} \H^1(L_w,\Zbb/n_p) = \H^1(L_w,A)$ pour toute place $w$ de $L$ divisant une place $v_{p,q}$, o\`u $p,q \in J$.
	    \end{enumerate}
	Soit $x = \sh^{-1}(a) \in \H^1(k,M)$. Pour tous $v \in S$ et $s \in \Escr_v$, on a $a_{sw_v} = \tensor[^s]{a}{}'_{v,s} \in \H^1(L_{sw_v}, A)$, d'o\`u $(\tensor[^{s^{-1}}]{a)}{_{w_v}} = a'_{s,v} \in \H^1(L_{w_v},A)$. Par le lemme \ref{lemShapiroLocalization1}, on a $x_v = x'_v$ pour toute $v \in S$. Ainsi \eqref{eqthmKey2} devient
		\begin{equation*}
                \forall v \in S, \quad x_v \cup y_v' = u_v + j_\ast \lambda_v \quad \text{dans } \H^2(k_v, M \otimes M).
            \end{equation*}
	Au vu des lemmes \ref{lemShapiroCupProductLocal}, \ref{lemShapiroJLocal}, \ref{lemShapiroLocalization1} et \ref{lemShapiroLocalization2}, on a
	   \begin{equation*}
                \forall v \in S,\forall h \in \gfrak_v, \forall s,t \in \Escr_v, \quad \tensor[^{h^{-1}s^{-1}}]{a}{_{w_v}} \cup b_{v,t}' = (\tensor[^{t^{-1}}]{\gamma}{_{sht^{-1}}})_{w_v} + \res(\lambda)_{w_v} \quad \text{dans } \H^2(L_{w_v},A \otimes A),
            \end{equation*}
	ou $\tensor[^{th^{-1}s^{-1}}]{a}{_{tw_v}} \cup \tensor[^t]{b}{}_{v,t}' = (\gamma_{sht^{-1}})_{tw_v} + \res(\lambda)_{tw_v}$ dans $\H^2(L_{tw_v},A \otimes A)$. Pour tous $g \in \gfrak$ et $t \in \Escr_v$, on peut \'ecrire $gt = sh$ pour un $s \in \Escr_v$ et un $h \in \gfrak_v$. On a ainsi
	\begin{equation} \label{eqthmKey3}
		\forall v \in S,\forall g \in \gfrak, \forall t \in \Escr_v, \quad \tensor[^{g^{-1}}]{a}{_{tw_v}} \cup \tensor[^t]{b}{}_{v,t}' = (\gamma_g)_{tw_v} + \res(\lambda)_{tw_v} \quad \text{dans } \H^2(L_{tw_v},A \otimes A).
	\end{equation}	
	
	{\bf \'Etape 3.} {\em V\'erifions l'hypoth\`ese du lemme arithm\'etique.} On va appliquer la proposition \ref{propArithmetic}
	\begin{itemize}
		\item au corps de nombres $L$;
		\item \`a la famille $(\tensor[^{g^{-1}}]{a}{})_{g \in \gfrak} \in \H^1(L,A)$;
		\item \`a la famille $(\gamma_g + \res(\lambda))_{g \in \gfrak} \in \H^2(L, A \otimes A)^{|\gfrak|}$;
		\item et \`a l'ensemble $\{w \in \Omega_L: \exists v \in S \cup \{v_{p,q}: p,q \in J\}, w|v\}$.
	\end{itemize} 
     Soit $w$ une place de $L$ et montrons qu'il existe $c_w \in \H^1(L_w,A)$ tel que
	\begin{equation*}
            \forall g \in \gfrak, \quad  \tensor[^{g^{-1}}]{a}{_w} \cup c_w = (\gamma_g)_w + \res(\lambda)_w \quad \text{dans } \H^2(L_w, A \otimes A).
        \end{equation*}
	Notons $v$ la place de $k$ au-dessous de $w$. On distingue trois cas.
		\begin{itemize}
			\item Si $v \in S$: Au vu de \eqref{eqthmKey3}, il suffit de prendre $c_w = \tensor[^t]{b}{}'_{v,t}$, o\`u $t \in \Escr_v$ est tel que $w = tw_v$. 
			
			\item Si $v = v_{p',q'}$, o\`u $p',q' \in J$: Rappelons que $v$ est totalement d\'ecompos\'ee dans $L$ et que $\varpi := \varpi_{p',q'} \in k^\times$ est une uniformisante de $k_v$. D'une part, $\tilde{a}^{p'}_w = \varpi \pmod{L_w^{\times n_{p'}}}$, donc l'image $a{}^{p'}_w \in \H^1(L_w,\mu_{n_{p'}})$ de $\tilde{a}^{p'}_w$ (par la th\'eorie de Kummer) est $\gfrak$-\'equivariante. Or $\tilde{a}^{p}_{w} \in L_w^{\times n_{p}}$ pour toute $p \in J \setminus \{p'\}$, donc l'image de $\tilde{a}^{p}_w$ dans  $\H^1(L_w,\mu_{n_p})$ (par la th\'eorie de Kummer) est $a{}^p_w  = 0$. En particulier, $a_w = (a^p_w)_{p \in J} \in \H^1(L_w,A)$ est $\gfrak$-\'equivariant. D'autre part, $(\gamma_g)_w = 0$ pour tous $g \in \gfrak$ puisque $u_v = 0$. Ainsi, il faut chercher $c_w = (c_w^{q})_{q \in J} \in \prod_{q \in J} \H^1(L_w,\mu_{n_q})$ tel que
			    \begin{equation*}
                    a_w \cup c_w = \res(\lambda)_w \quad \text{dans } \H^2(L_w, A \otimes A),
                \end{equation*}
			{\em i.e.}
			    \begin{equation} \label{eqthmKey4}
			        \forall p,q \in J, \quad  a{}^{p}_{w} \cup c_w^{q} = \res(\lambda^{p,q})_w \quad \text{dans } (\Br L_w)[d_{p,q}].
			     \end{equation}
		On choisit $c_w^{q} = 0$ pour tout $q \in J \setminus \{q'\}$. Pour le $c_w^{q'}$, on note $\alpha \in \H^1(L_w,\mu_{n_{q'}}) =  \H^1(L_w,\Zbb/n_{q'})$ l'image de $\tilde{a}^{p'}_w$ par la th\'eorie de Kummer, et l'on choisit $r \in (\Br L_w)[n_{q'}]$ tel que $\frac{n_{q'}}{d_{p',q'}}r = \res(\lambda^{p',q'})_w$. Comme $w(\tilde{a}^{p'}) = v(\tilde{a}^{p'}) = 1$, le lemme \ref{lemHilbertSymbol} assure qu'il existe $c_w^{q'} \in \H^1(L_w, \Zbb/n_{q'}) = \H^1(L_w, \mu_{n_{q'}})$ tel que $\alpha \cup c_w^{q'} = r$, d'o\`u $\frac{n_{q'}}{d_{p',q'}}
			\alpha \cup c_w^{q'} = \res(\lambda^{p',q'})_w$. Au vu de \eqref{eqKummerMN}, on a $(\iota^{n_{p'}}_{n_{q'}})_\ast a_w^{p'} = \frac{n_q'}{d_{p',q'}} \alpha \in \H^1(L_w,\Zbb/n_{q'})$, d'o\`u $a_w^{p'} \cup c_w^{q'} = \res(\lambda^{p',q'})_w$ en vertu de \eqref{eqPairingsMN}.
   
            On prend finalement $c_w = (c_w^q)_{q \in J}$. On a d\'ej\`a vu que cet \'el\'ement v\'erifie \eqref{eqthmKey4} pour $(p,q) = (p',q')$. Or, pour tout $(p,q) \in (J \times J) \setminus \{(p',q')\}$, on a soit $a^p_w = 0$ (si $p \neq p'$) ou $c^q_w = 0$ (si $q \neq q'$), et de plus $\lambda^{p,q}_v = 0$ par notre choix de $\lambda^{p,q}$, d'o\`u $\res(\lambda^{p,q})_w = 0 = a_w^p \cup c_w^q$, ce qui \'etablit \eqref{eqthmKey4} pour tous $p,q \in J$.

            \item Si $v \notin S \cup \{v_{p,q}: p,q \in J\}$. Dans ce cas, on a $u_v = 0$ et  $\lambda_v = 0$, d'o\`u $(\gamma_g)_w = \res(\lambda)_w = 0$ pour tous $g \in \gfrak$, donc il suffit de prendre $c_w = 0$.
		\end{itemize}
	L'hypoth\`ese de la proposition \ref{propArithmetic} est alors v\'erifi\'ee.
	
	{\bf \'Etape 4.} {\em Utilisons le lemme arithm\'etique pour approcher les $y'_v$.} La proposition \ref{propArithmetic} nous permet de fabriquer un \'el\'ement $b \in \H^1(L,A)$ satisfaisant les conditions suivantes.
	\begin{enumerate}
		\item $\tensor[^{g^{-1}}]{a}{} \cup b = \gamma_g + \res(\lambda)$ dans $\H^2(L, A \otimes A)$ pour tout $g \in \gfrak$.
		
		\item $b_{tw_v} = \tensor[^t]{b}{}'_{v,t}$ dans $\H^1(L_{tw_v},A)$ pour tous $v \in S$ et $t \in \Escr_v$. 
	\end{enumerate}  
	On note $y = \sh^{-1}(b) \in \H^1(k,M)$. Regardons les deux conditions ci-dessus.
	\begin{enumerate}
		\item Pour la premi\`ere condition, par les lemmes \ref{lemShapiroCupProduct} et \ref{lemShapiroJ}, on a $x \cup y = u + j_\ast \lambda$ dans $\H^2(k, M \otimes M)$.
		
		\item Soit $v \in S$. Pour tout $t \in \Escr_v$, la deuxi\`eme condition implique que $\tensor[^{t^{-1}}]{b}{_{w_v}} = b'_{t,v}$ dans $\H^1(L_{w_v},A)$. En vertu du lemme \ref{lemShapiroLocalization1}, on a $y_v = y'_v$. Rappelons qu'on a d\'ej\`a $x_v = x'_v$ et $\lambda_v = \lambda'_v$ pour toute $v \in S$.
	\end{enumerate}
	Le th\'eor\`eme est finalement d\'emontr\'e.
\end{proof}

Nous sommes maintenant en mesure d'\'etablir le th\'eor\`eme \ref{thmA}.

\begin{thm} \label{thmHassePrinciple}
	Soient $A$ un groupe ab\'elien fini, $k$ un corps de nombres contenant $\mu_{\exp(A)}$, $L/k$ une extension finie galoisienne, $M = \Ind_{\Gamma_k}^{\Gamma_L} A$, $j: A \otimes A \hookrightarrow M \otimes M$ l'inclusion canonique, $Z$ le conoyau de $j$ et $\phi: M \otimes M \to Z$ la projection. De $\phi$ on d\'efinit une application biadditive $\Phi: (M \oplus M) \times (M \oplus M) \to Z$ par \eqref{eqBigPhi} et l'on pose $F = Z \times_\Phi (M \oplus M)$. On munit $F$ de l'action coordonn\'ees par coordonn\'ees de $\Gamma_k$. Si $X$ est un espace homog\`ene de $\SL_n$ de lien de Springer $\lien(F)$ tel que $X(k_v) \neq \varnothing$ pour toute place $v$ de $k$, alors $X(k) \neq \varnothing$.
\end{thm}
\begin{proof}
	Notons $\Delta: \H^2(k,M \oplus M) \to \H^2(k,Z)$ l'application connectante induite par l'extension centrale $0 \to Z \to F \to M \oplus M \to 0$. Rappelons que $Z = Z(F) = [F,F]$ par la proposition \ref{propBorovoiKunyavskii}. Soit $\eta_0 \in \H^2(k,F)$ la classe neutre privil\'egi\'ee. Alors il existe une unique classe $\beta \in \H^2(k,Z)$ telle que la classe de Springer de $X$ vaut $\eta_X = \beta \cdot \eta_0$ ({\em cf.} paragraphe \ref{subsection31}). Soit $v$ une place de $\Omega_k$. Si $X(k_v) \neq \varnothing$, par la proposition \ref{propSpringerClass} et le lemme \ref{lemImageOfDelta}, il existera $(x'_v,y'_v) \in \H^1(k_v, M \oplus M)$ tel que $\beta_v = \Delta(x'_v,y'_v) = \phi_\ast(x'_v \cup y'_v)$ dans $\H^2(k_v,Z)$, donc l'image de $\beta_v$ dans $\H^3(k_v,A \otimes A)$ sera nulle. Ainsi, lorsque $X(k_v) \neq \varnothing$ pour toute $v$, l'image de $\beta$ dans $\H^3(k,A \otimes A)$ sera nulle partout localement. Or $\Sha^3(k,A \otimes A) = 0$ (c'est une cons\'equence d'une version simple de la dualit\'e de Poitou--Tate, {\em cf.} \cite[Th\'eor\`eme 17.13]{harari2017cohomologie}), donc $\beta$ s'enverra sur $0 \in \H^3(k,A \otimes A)$, d'o\`u il existera $u \in \H^2(k, M \otimes M)$ tel que $\beta = \phi_\ast(u)$.
	
	Soit $S \subseteq \Omega_k$ l'ensemble des places $v$ telles que $u_v \neq 0$. Pour toute $v \in S$, comme $\phi_\ast(x'_v \cup y'_v) = \beta_v = \phi_\ast(u_v)$ dans $\H^2(k_v, M \otimes M)$, il existe $\lambda'_v \in \H^2(k_v, A \otimes A)$ tel que $x'_v \cup y'_v = u_v + j_\ast \lambda'_v $ dans $\H^2(k_v, M \otimes M)$. Par le th\'eor\`eme \ref{thmKey}, on peut trouver $x,y \in \H^1(k,M)$ et $\lambda \in \H^2(k, A \otimes A)$ tels que $x \cup y = u + j_\ast \lambda$ dans $\H^2(k, M \otimes M)$. En particulier, $\beta = \phi_\ast(u) = \phi_\ast(x \cup y) = \Delta(x,y)$ par le lemme \ref{lemImageOfDelta}. Au vu de la proposition \ref{propSpringerClass}, on peut conclure que $X(k) \neq \varnothing$.
\end{proof}

\subsection{Approximation faible} \label{subsection54}

Gardons les notations au d\'ebut du paragraphe \ref{subsection42}. On a une suite exacte 
\begin{equation}\label{eqExactSequence}
	0 \to A \otimes A \xrightarrow{j} M \otimes M \xrightarrow{\phi} Z \to 0
\end{equation}
de $\Gamma_k$-modules. En dualisant, on obtient une suite exacte
\begin{equation} \label{eqDualizedExactSequence}
	0 \to \hat{Z} \xrightarrow{\hat{\phi}} \widehat{M \otimes M} \xrightarrow{\hat{j}} \widehat{A \otimes A}  \to 0.
\end{equation}
Notons que $\Gamma_k$ agit trivalement sur $A \otimes A$ et sur $\widehat{ A \otimes A}$.

\begin{lem} \label{lemShaOfZ}
	Le sous-groupe $\Sha^1_\omega(k,\hat{Z}) \subseteq \H^1(k,\hat{Z})$ est inclus dans l'image du morphisme connectant $\hat{\delta}: \widehat{ A \otimes A} \to \H^1(k,\hat{Z})$ induit par \eqref{eqDualizedExactSequence}.
\end{lem}
\begin{proof}
	Comme $\eqref{eqDualizedExactSequence}$ est d\'eploy\'e par $L$, on peut identifier $\H^1(\gfrak,\hat{Z})$ \`a un sous-groupe de $\H^1(k,\hat{Z})$ (par la suite exacte d'inflation-restriction). Une application du  th\'eor\`eme de Chebotarev donne $\Sha^1_\omega(L,\hat{Z}) = 0$ (voir \cite[Corollaire 18.4]{harari2017cohomologie}), donc on a $\Sha^1_\omega(k,\hat{Z}) = \Sha^1_\omega(\gfrak,\hat{Z})  \subseteq \H^1(\gfrak,\hat{Z})$. Maintenant, $M \otimes M$ est un $\gfrak$-module induit par le lemme \ref{lemShapiroTensor}, donc il en va de m\^eme pour $\widehat{M \otimes M}$. Par le lemme de Shapiro, $\H^1(\gfrak,\widehat{M \otimes M}) = 0$ et donc le morphisme connectant $\widehat{A \otimes A} \to \H^1(\gfrak,\hat{Z})$ est surjectif. D'o\`u le lemme.
\end{proof}

On d\'emontre maintenant le th\'eor\`eme \ref{thmB}.

\begin{thm} \label{thmWeakApproximation}
	Soient $A$ un groupe ab\'elien fini, $k$ un corps de nombres contenant $\mu_{\exp(A)}$, $L/k$ une extension finie galoisienne, $M = \Ind_{\Gamma_k}^{\Gamma_L} A$, $j: A \otimes A \hookrightarrow M \otimes M$ l'inclusion canonique, $Z$ le conoyau de $j$ et $\phi: M \otimes M \to Z$ la projection. De $\phi$ on d\'efinit une application biadditive $\Phi: (M \oplus M) \times (M \oplus M) \to Z$ par \eqref{eqBigPhi} et l'on pose $F = Z \times_\Phi (M \oplus M)$. On munit $F$ de l'action coordonn\'ees par coordonn\'ees de $\Gamma_k$. Si $X$ est un espace homog\`ene de $\SL_n$ de lien de Springer $\lien(F)$ tel que $X(k) \neq \varnothing$, alors $X$ v\'erifie l'approximation faible.
\end{thm}
\begin{proof}
	Par la proposition \ref{propSpringerClass},  $X \simeq \tensor[_\afrak]{F}{} \setminus \SL_n$ pour un certain $1$-cocycle $\afrak = (\xfrak,\yfrak): \Gamma_k \to M \oplus M$. Par le lemme \ref{lemWA}, on s'est ramen\'e \`a la question d'approximation faible pour $\tensor[_\afrak]{F}{}$, {\em i.e.} il faut d\'emontrer que pour tout sous-ensemble fini $S \subseteq \Omega_k$, la restriction
	\begin{equation*}
            \H^1(k,\tensor[_\afrak]{F}{}) \to \prod_{v \in S} \H^1(k_v,\tensor[_\afrak]{F}{})
        \end{equation*}
	est surjective. Soit alors $f_v' = (z_v',a_v'): \Gamma_{k_v} \to \tensor[_\afrak]{F}{}$ une famille de $1$-cocycles, $v \in S$, o\`u $a_v' = (x_v',y_v')$. Pour toute $v \in S$, par le lemme \ref{lemImageOfDelta}, $a'_v$ est un cocycle et de plus
	\begin{equation*}
            \diff z'_v + \phi_\ast(\xfrak_v \cup y'_v + x'_v \cup \yfrak_v + x'_v \cup y'_v + \diff(x'_v \otimes \yfrak_v)) = 0 \quad \text{dans } \Z^2(k_v, Z).
        \end{equation*}
	\'Ecrivons $z'_v = \phi_\ast \varepsilon'_v$ pour une $1$-cocha\^ine $\varepsilon'_v: \Gamma_{k_v} \to M \otimes M$. Il existe alors un $2$-cocycle $\lambda'_v: \Gamma_{k_v} \times \Gamma_{k_v} \to A \otimes A$ tel que
		\begin{equation} \label{eqthmWeakApproximation1}
			\diff \varepsilon'_v + \xfrak_v \cup y'_v + x'_v \cup \yfrak_v + x'_v \cup y'_v + \diff(x'_v \otimes \yfrak_v) = j_\ast \lambda'_v \quad \text{dans } \Z^2(k_v, M \otimes M).
		\end{equation}
Consid\'erons l'\'el\'ement $u = [\xfrak \cup \yfrak] \in \H^2(k, M \otimes M)$. Quitte \`a agrandir $S$, on peut supposer que $u_v = 0$ pour toute $v \notin S$. De \eqref{eqthmWeakApproximation1}, on a
	\begin{equation*}
            \forall v \in S, \quad [(x'_v + \xfrak_v) \cup (y'_v + \yfrak_v)] = [\xfrak_v \cup \yfrak_v] + j_\ast [\lambda'_v] \quad \text{dans } \H^2(k_v, M \otimes M).
        \end{equation*}
	Par le th\'eor\`eme \ref{thmKey}, il existe des $1$-cocycles $\tilde{x},\tilde{y}: \Gamma_k \to M$ et un $2$-cocycle $\lambda: \Gamma_k \times \Gamma_k \to A \otimes A$ tels que les conditions suivantes soient remplies.
		\begin{enumerate}
			\item $[\tilde{x} \cup \tilde{y}] = [\xfrak \cup \yfrak] + j_\ast [\lambda]$ dans $\H^2(k,M \otimes M)$.
			\item pour toute $v \in S$, $[\tilde{x}_v] = [x'_v + \xfrak_v]$, $[\tilde{y}_v] = [y'_v + \yfrak_v]$ et $[\lambda_v] = [\lambda'_v]$.
		\end{enumerate} 
	Soient $x = \tilde{x} - \xfrak$, $y = \tilde{y} - \yfrak$ et $a = (x,y)$. Alors $[x_v] = [x'_v]$ et $[y_v] = [y'_v]$ pour toute $v \in S$. De $[\tilde{x} \cup \tilde{y}] = [\xfrak \cup \yfrak] + j_\ast [\lambda]$, on voit qu'il existe une $1$-cocha\^ine $\varepsilon: \Gamma_k \to M \otimes M$ telle que 
		\begin{equation} \label{eqthmWeakApproximation2}
			\diff \varepsilon + \xfrak \cup y + x \cup \yfrak + x \cup y + \diff(x \otimes \yfrak) = j_\ast \lambda.
		\end{equation} 
	Posons $z := \phi_\ast \varepsilon$, alors \eqref{eqthmWeakApproximation2} implique que $\diff z + \phi_\ast(\xfrak \cup y + x \cup \yfrak + x \cup y + \diff(x \otimes \yfrak)) = 0$ et donc $(z,a): \Gamma_k \to F$ est un cocycle par le lemme \ref{lemImageOfDelta}.
	
	Soit $v \in S$. Puisque $[x'_v] = [x_v]$ et $[y'_v] = [y_v]$ dans $\H^1(k_v,M)$, on peut \'erire $x'_v = x_v + \diff \xi_v$ et $y'_v = y_v + \diff \eta_v$, o\`u $\xi_v,\eta_v \in M$. On consid\`ere le $1$-cocycle
	\begin{equation*}
            c_v := z'_v - z_v + \phi_\ast(-(x_v + \xfrak_v) \cup \eta_v + \xi_v \cup (y'_v + \yfrak_v) + \diff\xi_v \otimes \yfrak_v): \Gamma_{k_v} \to Z
        \end{equation*}
	comme dans l'\'enonc\'e du lemme \ref{lemCohomologousOnF}. Au vu de \eqref{eqthmWeakApproximation1} et \eqref{eqthmWeakApproximation2}, le troisi\`eme point de ce lemme-l\`a donne $\delta([c_v]) = [\lambda'_v - \lambda_v] = 0$, o\`u $\delta: \H^1(k_v,Z) \to \H^2(k_v,A \otimes A)$ d\'esigne le morphisme connectant induit par \eqref{eqExactSequence}. Or, en vertu du lemme \ref{lemShaOfZ} tout \'el\'ement de $\Sha^1_S(k,\hat{Z})$ appartient \`a l'image du morphisme connectant $\hat{\delta}: \widehat{A \otimes A} \to \H^1(k,\hat{Z})$ induite par \eqref{eqDualizedExactSequence}. Puisque les cup-produits sont compatibles avec les morphismes connectants, on voit que la famille $([c_v])_{v \in S}$ est orthogonal \`a $\Sha^1_S(k,\hat{Z})$ par rapport \`a l'accouplement
	\begin{equation*}
            \prod_{v \in S} \H^1(k_v,Z) \times \H^1(k,\hat{Z}) \to \Qbb/\Zbb,
        \end{equation*}
	obtenu en rassemblant les dualit\'es locales de Tate. Par \cite[Exercice 17.5(a)]{harari2017cohomologie}, il existe un $1$-cocycle $\tilde{z}: \Gamma_{k} \to Z$ tel que $[\tilde{z}_v] = [c_v]$ dans $\H^1(k_v,Z)$ pour toute $v \in S$. On pose finalement
	\begin{equation*}
            f := (\tilde{z},0)(z,a) = (z + \tilde{z},a): \Gamma_k \to \tensor[_\afrak]{F}{},
        \end{equation*}
    qui est un $1$-cocycle puisque $(z,a)$ et $(\tilde{z},0)$ le sont (notons que $Z$ est central dans $\tensor[_\afrak]{F}{}$). Soit $v \in S$, il reste \`a montrer que $[f_v] = [f'_v]$ dans $\H^1(k_v,\tensor[_\afrak]{F}{})$. En effet, la $1$-cocha\^ine
	\begin{equation*}
            z'_v - (z_v + \tilde{z}_v) + \phi_\ast(-(x_v + \xfrak_v) \cup \eta_v + \xi_v \cup (y'_v + \yfrak_v) + \diff\xi_v \otimes \yfrak_v) = c_v - \tilde{z}_v: \Gamma_k \to Z
        \end{equation*}
	est un cobord, donc $f_v$ est cohomologue \`a $f'_v$ par le lemme \ref{lemCohomologousOnF}. On a donc trouv\'e une classe $[f] \in \H^1(k,\tensor[_a]{F}{})$ qui se restreint \`a $([f'_v])_{v \in S} \in \prod_{v \in S} \H^1(k_v,\tensor[_\afrak]{F}{})$. Ainsi $\tensor[_a]{F}{}$ v\'erifie bien l'approximation faible, ce qui ach\`eve la d\'emonstration du th\'eor\`eme.
\end{proof}   

\begin{remercie}
    Ce travail est financ\'e par un \og Contrat doctoral sp\'ecifique normalien\fg{} de l'\'Ecole normale sup\'erieure de Paris. L'auteur est reconnaissant \`a David Harari pour ses commentaires pr\'ecieux ainsi que pour son soutien. L'auteur remercie Alexei Skorobogatov, Jean-Louis Colliot-Th\'el\`ene et Mikhail Borovoi pour les discussions pertinentes. L'auteur remercie les \'evaluateurs d'avoir lu attentivement le texte et pour ses remarques.
\end{remercie}
	
\bibliographystyle{alpha-fr}
\selectlanguage{french}
\bibliography{ref}

\begin{thebibliography}{DLA19}
\expandafter\ifx\csname fonteauteurs\endcsname\relax
\def\fonteauteurs{\scshape}\fi

\bibitem[BK97]{Bokun}
Mikhail \bgroup\fonteauteurs\bgroup Borovoi\egroup\egroup{} et Boris \bgroup\fonteauteurs\bgroup Kunyavskii\egroup\egroup{} :
\newblock On the {H}asse principle for homogeneous spaces with finite stabilizers.
\newblock {\em Annales de la Facult\'e des sciences de Toulouse : Math\'ematiques}, Ser. 6, 6(3)\string:\penalty500\relax 481--497, 1997.

\bibitem[BK01]{BokunErratum}
Mikhail \bgroup\fonteauteurs\bgroup Borovoi\egroup\egroup{} et Boris \bgroup\fonteauteurs\bgroup Kunyavskii\egroup\egroup{} :
\newblock Erratum : On the {H}asse principle for homogeneous spaces with finite stabilizers.
\newblock {\em Annales de la Facult\'e des sciences de Toulouse : Math\'ematiques}, Ser. 6, 10(4)\string:\penalty500\relax 779--779, 2001.

\bibitem[Bog88]{Bogomolov}
Fiodor~A. \bgroup\fonteauteurs\bgroup Bogomolov\egroup\egroup{} :
\newblock {The Brauer group of quotient spaces by linear group actions}.
\newblock {\em Mathematics of the {USSR}-Izvestiya}, 30(3)\string:\penalty500\relax 455--485, 1988.

\bibitem[Bor93]{Borovoi93}
Mikhail \bgroup\fonteauteurs\bgroup Borovoi\egroup\egroup{} :
\newblock {Abelianization of the second nonabelian Galois cohomology}.
\newblock {\em Duke Mathematical Journal}, 72(1)\string:\penalty500\relax 217--239, 1993.

\bibitem[Bor96]{bor96}
Mikhail \bgroup\fonteauteurs\bgroup Borovoi\egroup\egroup{} :
\newblock {The Brauer-Manin obstructions for homogeneous spaces with connected or abelian stabilizer}.
\newblock {\em Journal f{\"u}r die reine und angewandte Mathematik}, 473\string:\penalty500\relax 181--194, 1996.

\bibitem[CTS21]{colliot2021brauer}
Jean-Louis \bgroup\fonteauteurs\bgroup Colliot-Th\'el\`ene\egroup\egroup{} et Alexei~N. \bgroup\fonteauteurs\bgroup Skorobogatov\egroup\egroup{} :
\newblock {\em The Brauer--Grothendieck Group}, volume~71 de {\em Ergebnisse der Mathematik und ihrer Grenzgebiete. 3. Folge / A Series of Modern Surveys in Mathematics}.
\newblock Springer International Publishing, 2021.

\bibitem[Dem10]{Demarche}
Cyril \bgroup\fonteauteurs\bgroup Demarche\egroup\egroup{} :
\newblock {Groupe de Brauer non ramifi{\'e} d'espaces homog{\`e}nes {\`a} stabilisateurs finis}.
\newblock {\em {Mathematische Annalen}}, 346(4)\string:\penalty500\relax 949--968, 2010.

\bibitem[DLA19]{demarche2019reduction}
Cyril \bgroup\fonteauteurs\bgroup Demarche\egroup\egroup{} et Giancarlo \bgroup\fonteauteurs\bgroup Lucchini~Arteche\egroup\egroup{} :
\newblock {Le principe de Hasse pour les espaces homog{\`e}nes : r{\'e}duction au cas des stabilisateurs finis}.
\newblock {\em {Compositio Mathematica}}, 155(8)\string:\penalty500\relax 1568--1593, 2019.

\bibitem[FSS98]{FSS}
Yuval~F. \bgroup\fonteauteurs\bgroup Flicker\egroup\egroup{}, Claus \bgroup\fonteauteurs\bgroup Scheiderer\egroup\egroup{} et Ramdorai \bgroup\fonteauteurs\bgroup Sujatha\egroup\egroup{} :
\newblock {Grothendieck's theorem on non-abelian $H^2$ and local-global principles}.
\newblock {\em Journal of the American Mathematical Society}, 11(3)\string:\penalty500\relax 731--750, 1998.

\bibitem[Har07]{harari_affbm}
David \bgroup\fonteauteurs\bgroup Harari\egroup\egroup{} :
\newblock Quelques propri\'et\'es d'approximation reli\'ees \`a la cohomologie galoisienne d'un groupe alg\'ebrique fini.
\newblock {\em Bulletin de la Soci\'et\'e Math\'ematique de France}, 135(4)\string:\penalty500\relax 549--564, 2007.

\bibitem[Har17]{harari2017cohomologie}
David \bgroup\fonteauteurs\bgroup Harari\egroup\egroup{} :
\newblock {\em Cohomologie galoisienne et th{\'e}orie du corps de classes}.
\newblock Savoirs Actuels. EDP Sciences, 2017.

\bibitem[HS02]{HS02}
David \bgroup\fonteauteurs\bgroup Harari\egroup\egroup{} et Alexei~N. \bgroup\fonteauteurs\bgroup Skorobogatov\egroup\egroup{} :
\newblock Non-abelian cohomology and rational points.
\newblock {\em Compositio Mathematica}, 130(3)\string:\penalty500\relax 241--273, 2002.

\bibitem[HW20]{Harpaz2020}
Yonatan \bgroup\fonteauteurs\bgroup Harpaz\egroup\egroup{} et Olivier \bgroup\fonteauteurs\bgroup Wittenberg\egroup\egroup{} :
\newblock {Z\'ero-cycles sur les espaces homog\`enes et probl\`eme de Galois inverse}.
\newblock {\em Journal of the American Mathematical Society}, 33(3)\string:\penalty500\relax 775–805, 2020.

\bibitem[Man71]{manin1971}
Yuri~I. \bgroup\fonteauteurs\bgroup Manin\egroup\egroup{} :
\newblock {Le groupe de Brauer--Grothendieck en g\'eom\'etrie diophantienne}.
\newblock \emph{In} {\em Actes du Congr{\`e}s international des Math{\'e}maticiens (Nice, 1970)}, pages 401--411. Gauthier-Villars, Paris, 1971.

\bibitem[Mor12]{Moravec}
Primo$\check{\textrm{z}}$ \bgroup\fonteauteurs\bgroup Moravec\egroup\egroup{} :
\newblock {Unramified Brauer groups of finite and infinite groups}.
\newblock {\em American Journal of Mathematics}, 134(6)\string:\penalty500\relax 1679--1704, 2012.

\bibitem[Nai07]{Naidu}
Deepak \bgroup\fonteauteurs\bgroup Naidu\egroup\egroup{} :
\newblock {Categorical Morita equivalence for group-theoretical categories}.
\newblock {\em Communications in Algebra}, 35(11)\string:\penalty500\relax 3544--3565, 2007.

\bibitem[Ser77]{Serre}
Jean-Pierre \bgroup\fonteauteurs\bgroup Serre\egroup\egroup{} :
\newblock {\em Cours d'arithm{\'e}tique}, volume~2 de {\em Collection SUP.: Math{\'e}maticien}.
\newblock Presses universitaires de France, 1977.

\bibitem[Ser94]{serre1994cohomologie}
Jean-Pierre \bgroup\fonteauteurs\bgroup Serre\egroup\egroup{} :
\newblock {\em {Cohomologie Galoisienne: Cinqui\`eme \'edition, r\'evis\'ee et compl\'et\'ee}}, volume~5 de {\em Lecture Notes in Mathematics}.
\newblock Springer-Verlag Berlin Heidelberg, 1994.

\bibitem[Ser04]{Serre2004Locaux}
Jean-Pierre \bgroup\fonteauteurs\bgroup Serre\egroup\egroup{} :
\newblock {\em Corps Locaux}, volume 1296 de {\em Actualit{\'e}s scientifiques et industrielles}.
\newblock Hermann, 4th \'edition, 2004.

\bibitem[Sko01]{skorobogatov2001torsors}
Alexei~N. \bgroup\fonteauteurs\bgroup Skorobogatov\egroup\egroup{} :
\newblock {\em {Torsors and Rational Points}}, volume 144 de {\em Cambridge Tracts in Mathematics}.
\newblock Cambridge University Press, 2001.

\end{thebibliography}
	
\end{document}